\documentclass[11pt]{amsart}

\usepackage{subfiles}
\usepackage[a4paper,top=3cm,bottom=3cm,inner=3cm,outer=3cm]{geometry}
\usepackage[foot]{amsaddr}
\usepackage{mathtools}
\usepackage[T1]{fontenc}
\usepackage{inputenc}
\usepackage{microtype}
\usepackage[usenames,dvipsnames]{xcolor}
\usepackage{enumitem}
\usepackage{booktabs}
\usepackage{array}
\usepackage{pbox}
\usepackage{tikz}
\usepackage[textsize=footnotesize]{todonotes}
\usepackage{multirow,comment}
\usepackage{hyperref}
\hypersetup{colorlinks,linkcolor={blue!50!black},citecolor={red!40!black},urlcolor={blue!75!black},linktoc=page}
\usepackage{amssymb}
\usepackage{amsthm}
\usepackage{mathrsfs}
\usepackage{paralist}
\usepackage{bussproofs}
\EnableBpAbbreviations
\usepackage[theoremfont]{newpxtext}
\usepackage[varg,bigdelims,nosymbolsc,noamssymbols]{newpxmath}
\usepackage[cal=euler]{mathalfa}

\usepackage{bm}
\usepackage[capitalize]{cleveref}

\crefname{equation}{}{}
\crefname{condition}{condition}{conditions}
\crefname{item}{}{}
\crefname{section}{Chapter}{Chapters}
\crefname{subsection}{Section}{Sections}

\usetikzlibrary{positioning,decorations.pathreplacing,decorations.markings,arrows.meta,calc,fit,quotes,cd,math,arrows,backgrounds,shapes.geometric}
\hypersetup{final}

\setlist{nosep,leftmargin=*,labelsep=.5em}
\setlist[1]{labelindent=\parindent}

\theoremstyle{plain}
\newtheorem{theorem}{Theorem}[subsection] 
\newtheorem{proposition}[theorem]{Proposition}
\newtheorem{corollary}[theorem]{Corollary}
\newtheorem{lemma}[theorem]{Lemma}

\theoremstyle{definition}
\newtheorem{definition}[theorem]{Definition}

\newtheorem{notation}[theorem]{Notation}

\theoremstyle{remark}
\newtheorem{example}[theorem]{Example}
\newtheorem{remark}[theorem]{Remark}

\newcommand{\ncat}[1]{\mathbf{#1}} 
\newcommand{\nncat}[2]{\bm{\mathcal{#1}}\ncat{#2}} 
\newcommand{\TFS}{\ncat{TFS}}
\newcommand{\Set}{\ncat{Set}}
\newcommand{\Span}{\ncat{Span}}
\newcommand{\FinSet}{\ncat{FinSet}}
\newcommand{\Cat}{\ncat{Cat}}
\newcommand{\Fam}{\ncat{Fam}}
\newcommand{\Cob}{\ncat{Cob}}
\newcommand{\Grph}{\ncat{Grph}}

\newcommand{\WGrph}{\ncat{wGrph}}
\newcommand{\Man}{\ncat{Man}}
\newcommand{\Euc}{\ncat{Euc}}
\newcommand{\Alg}{\ncat{Alg}}

\newcommand{\WDalg}{\nncat{W}{D}\alg}
\newcommand{\WMonHom}{\ncat{WMonHom}}
\newcommand{\Int}{\ncat{Int}}

\newcommand{\FPCat}{\ncat{FPCat}}
\newcommand{\FCCat}{\ncat{FCCat}}
\newcommand{\RegCat}{\ncat{RegCat}}
\newcommand{\RegCate}{\RegCat_e}
\newcommand{\1}{\ncat{1}}
\newcommand{\DCF}{\ncat{DCF}}

\newcommand{\SOpd}{\nncat{S}{Opd}}
\newcommand{\SMulticat}{\nncat{S}{Multicat}}
\newcommand{\SMonCatl}{\nncat{S}{MC}_\mathrm{\ell}}
\newcommand{\SMonCats}{\nncat{S}{MC}} 

\newcommand{\nfun}[1]{\mathsf{#1}} 
\newcommand{\DDS}{\nfun{DDS}}
\newcommand{\CDS}{\nfun{CDS}}
\newcommand{\SpnS}{\nfun{Spn}} 
\newcommand{\ADS}[1][]{\nfun{Mch_{\mathit{#1}}}}
\newcommand{\dADS}[1][]{\nfun{Mch^d_{\mathit{#1}}}}
\newcommand{\tADS}[1][]{\nfun{Mch^t_{\mathit{#1}}}}

\newcommand{\tdADS}[1][]{\nfun{Mch^{td}_{\mathit{#1}}}}

\newcommand{\iADS}[1][]{\nfun{Mch^{in}_{\mathit{#1}}}}

\newcommand{\iADSe}[1][]{\nfun{Mch^{\epsilon\textrm{-}in}_{\mathit{#1}}}}
\newcommand{\iADSd}[1][]{\nfun{Mch^{\delta\textrm{-}in}_{\mathit{#1}}}}
\newcommand{\xADS}[1][]{\nfun{Mch^\mathbf{x}_{\mathit{#1}}}}
\newcommand{\Cntr}{\nfun{Cntr}}
\newcommand{\Fnc}{\nfun{Fnc}}
\newcommand{\Ext}[1]{\nfun{Ext}_{#1}}

\newcommand{\Tr}{\nfun{Tr}}
\newcommand{\im}{\nfun{im}}
\newcommand{\Img}{\nfun{Img}}
\newcommand{\CM}{\nfun{CM}}
\newcommand{\Fr}{\nfun{Fr}}

\newcommand{\Traj}[1][1]{C^{#1}}

\newcommand{\len}{\nfun{len}}
\newcommand{\Stk}{\nfun{Stk}_0}
\newcommand{\End}{\nfun{End}}

\newcommand{\tn}[1]{\textnormal{#1}}

\newcommand{\upd}[1]{#1^{\tn{upd}}}
\newcommand{\rdt}[1]{#1^{\tn{rdt}}}
\newcommand{\dyn}[1]{#1^{\tn{dyn}}}
\newcommand{\op}{^{\tn{op}}}
\newcommand{\Sync}{\tn{Sync}}
\newcommand{\Lan}{\tn{Lan}}

\newcommand{\sdom}[1]{\lambda_0(#1)}
\newcommand{\scod}[1]{\rho_0(#1)}
\newcommand{\Trns}[1]{\tn{Tr}_{#1}}
\newcommand{\src}{\tn{src}}
\newcommand{\tgt}{\tn{tgt}}

\newcommand{\Fact}{\tn{Fact}}

\newcommand{\tw}[1]{#1_{\tn{tw}}}
\newcommand{\cnst}{\tn{cnst}}
\newcommand{\Path}{\tn{Path}}

\newcommand{\wh}{\widehat}
\newcommand{\inp}[1]{{#1}^\mathrm{in}}
\newcommand{\out}[1]{{#1}^\mathrm{out}}
\newcommand{\ii}[1]{{#1}^\mathrm{i}}
\newcommand{\oo}[1]{{#1}^\mathrm{o}}
\newcommand{\dProd}[1]{\wh{#1}} 
\newcommand{\Xin}[1][]{\dProd{\inp{X}_{#1}}}
\newcommand{\Xout}[1][]{\dProd{\out{X}_{#1}}}
\newcommand{\Yin}{\dProd{\inp{Y}}}
\newcommand{\Yout}{\dProd{\out{Y}}}
\newcommand{\Zin}{\dProd{\inp{Z}}}
\newcommand{\Zout}{\dProd{\out{Z}}}
\newcommand{\fin}{\dProd{\inp{\phi}}}

\newcommand{\finl}[1]{\dProd{\inp{\phi}_{#1}}}
\newcommand{\fout}{\dProd{\out{\phi}}}

\newcommand{\wt}[1]{\widetilde{#1}}
\newcommand{\ul}[1]{\underline{#1}}
\newcommand{\ol}[1]{\overline{#1}}
\newcommand{\Shv}[1]{\wt{#1}}
\newcommand{\ascat}[1]{\ol{#1}} 
\newcommand{\cat}[1]{\mathcal{#1}}
\newcommand{\singleton}[1][\ast]{\{#1\}}
\newcommand{\RR}{{\mathbb{R}}}
\newcommand{\RRnn}{{\mathbb{R}_{\geq0}}} 

\newcommand{\NN}{{\mathbb{N}}} 
\newcommand{\ZZ}{\mathbb{Z}}
\newcommand{\WW}[1]{\cat{W}_{#1}}
\newcommand{\ti}{\ell}
\newcommand{\Psh}[1]{\mathsf{Psh}(#1)}
\newcommand{\SubPsh}[1]{\mathsf{SubPsh}(#1)}
\newcommand{\Sub}{\mathsf{Sub}}
\newcommand{\simrightarrow}{\xrightarrow{\raisebox{-4pt}[0pt][0pt]{\ensuremath{\sim}}}}
\newcommand{\undop}{\cat{O}}
\newcommand{\Inp}{\mathrm{-}}
\newcommand{\dash}{\textrm{-}}

\newcommand{\scirc}{{\hspace{1pt}\scriptstyle\circ\hspace{1pt}}}

\newcommand{\bbot}{\mathbin{\rotatebox[origin=c]{90}{$\bot$}}}
\newcommand{\To}[1]{\xrightarrow{#1}}
\newcommand{\alg}{\dash\Alg}

\newcommand{\Conduche}{Conduch\'{e}\ }
\renewcommand{\ss}{\subseteq}
\newcommand{\yoneda}[1]{\mathsf{Yon}_{#1}}
\newcommand{\fakeell}{{\color{white}{\ell}}}

\newcommand{\ShRst}[4][]{#2|_{#1[#3,#4#1]}}
\newcommand{\ullimit}{\ar[dr,phantom,very near start,"\lrcorner" description]}

\newcommand{\from}{\leftarrow}
\newcommand{\From}[1]{\xleftarrow{#1}}
\newcommand{\tto}{\rightrightarrows}

\renewcommand{\implies}{$\Rightarrow$}
\newcommand{\assShf}{\mathsf{asSh}}

\newcommand{\algmap}[6]
{\begin{tikzcd}[column sep=.7in,row sep=.2in,ampersand replacement=\&]\WW{#1}\ar[dr,"{#4}"]\ar[dd,"\WW{#2}"']\&\\
\phantom{A}\ar[r,phantom,pos=.35,"\Downarrow{\scriptstyle #6}"description]\&\Cat\\
\WW{#3}\ar[ur,"{#5}"']\&
\end{tikzcd}}

\DeclareMathOperator{\Ob}{ob}

\DeclareMathOperator{\id}{id}
\DeclareMathOperator*{\colim}{colim}
\DeclareMathOperator{\Hom}{Hom}
\DeclarePairedDelimiter{\ceil}{\lceil}{\rceil}
\DeclarePairedDelimiter{\floor}{\lfloor}{\rfloor}

\newcommand{\twocell}[3][]{\arrow[draw=none,to path={(dom#2.center)--(cod#2.center)\tikztonodes}]{}[anchor=center,#1]{\Downarrow #3}}

\newcommand{\twocellA}[2][]{\twocell[#1]{A}{#2}}

\newcommand{\SmallBox}[3]
{\begin{tikzpicture}[oriented WD, baseline=(Y.south), bb Small]
\node[inner sep=.1cm] [bb={1}{1}] (X) {$\scriptstyle #3$};
\draw[label] node[left=.1 of X_in1] (Y) {$#1$}
             node[right=.1 of X_out1] {$#2$};
\end{tikzpicture}}

\tikzcdset{
   arrow style=tikz,
   diagrams={>={Classical TikZ Rightarrow[angle=63:4pt, line width=.6pt]}},
   arrows={semithick}
}

\tikzset{
  tick/.style={postaction={
    decorate,
    decoration={markings, mark=at position 0.5 with {\draw[-] (0,.4ex) -- (0,-.4ex);}}}
  },
  tickx/.style={
    postaction={ decorate,
      decoration={markings,
        mark=at position 0.5 with {
          \fill circle [radius=.28ex];
        }
      }
    }
  }
}
\newcommand{\tickar}
{\begin{tikzcd}[baseline=-0.5ex,cramped,sep=small,ampersand replacement=\&]{}\ar[r,tick]\&{}\end{tikzcd}}
\tikzset{
   dom/.style={append after command={coordinate[alias=dom#1]}},
   domA/.style={dom=A}, domB/.style={dom=B},
   domC/.style={dom=C}, domD/.style={dom=D},
   domE/.style={dom=E}, domF/.style={dom=F},
   cod/.style={append after command={coordinate[alias=cod#1]}},
   codA/.style={cod=A}, codB/.style={cod=B},
   codC/.style={cod=C}, codD/.style={cod=D},
   codE/.style={cod=E}, codF/.style={cod=F}
}

\tikzset{
   oriented WD/.style={
      every to/.style={out=0,in=180,draw},
      label/.style={
         font=\everymath\expandafter{\the\everymath\scriptstyle},
         inner sep=0pt,
         node distance=2pt and -2pt},
      semithick,
      node distance=1 and 1,
      decoration={markings, mark=at position .5 with {\arrow{stealth};}},
      ar/.style={postaction={decorate}},
      execute at begin picture={\tikzset{
         x=\bbx, y=\bby,
         every fit/.style={inner xsep=\bbx, inner ysep=\bby}}}
      },
   bbx/.store in=\bbx,
   bbx = 1.5cm,
   bby/.store in=\bby,
   bby = 1.75ex,
   bb port sep/.store in=\bbportsep,
   bb port sep=2,
   bb port length/.store in=\bbportlen,
   bb port length=4pt,
   bb min width/.store in=\bbminwidth,
   bb min width=1cm,
   bb rounded corners/.store in=\bbcorners,
   bb rounded corners=2pt,
   bb small/.style={bb port sep=1, bb port length=2.5pt, bbx=.4cm, bb min width=.4cm, bby=.7ex},
   bb Small/.style={bb port sep=1, bb port length=2.5pt, bbx=.5cm, bb min width=.5cm, bby=1ex},
   bb/.code 2 args={
      \pgfmathsetlengthmacro{\bbheight}{\bbportsep * (max(#1,#2)+1) * \bby}
      \pgfkeysalso{draw,minimum height=\bbheight,minimum width=\bbminwidth,outer sep=0pt,
         rounded corners=\bbcorners,thick,
         prefix after command={\pgfextra{\let\fixname\tikzlastnode}},
         append after command={\pgfextra{\draw
            \ifnum #1=0{} \else foreach \i in {1,...,#1} {
               ($(\fixname.north west)!{\i/(#1+1)}!(\fixname.south west)$) +(-\bbportlen,0) coordinate
               (\fixname_in\i) -- +(\bbportlen,0) coordinate (\fixname_in\i')}\fi 
            \ifnum #2=0{} \else foreach \i in {1,...,#2} {
               ($(\fixname.north east)!{\i/(#2+1)}!(\fixname.south east)$) +(-\bbportlen,0) coordinate
               (\fixname_out\i') -- +(\bbportlen,0) coordinate (\fixname_out\i)}\fi;
         }}}
   },
   bb name/.style={append after command={\pgfextra{\node[anchor=north] at (\fixname.north) {#1};}}}
}

\allowdisplaybreaks[1]

\title{Dynamical Systems and Sheaves}
\author{Patrick Schultz}
\address[PS]{Broad Institute, 415 Main St, Cambridge, MA 02142, USA}
\email{pschultz@broadinstitute.org}
\author{David I. Spivak}
\address[DS]{Massachusetts Institute of Technology, 77 Massachusetts Ave. Cambridge, MA 02139}
\email{dspivak@mit.edu}
\author{Christina Vasilakopoulou}
\address[CV]{University of California, Riverside, 900 University Avenue, 92521, USA}
\email{vasilak@ucr.edu}

\begin{document}

\begin{abstract}
A categorical framework for modeling and analyzing systems
in a broad sense is proposed. These systems should be thought of as `machines' with inputs and outputs,
carrying some sort of signal that occurs through some notion of time.
Special cases include continuous and discrete dynamical systems (e.g. Moore machines). Additionally,
morphisms between the different types of systems allow their translation in a common
framework. A central goal is to understand the systems that result from arbitrary interconnection  of component subsystems, possibly of different types, as well as establish conditions that ensure totality and determinism compositionally.
The fundamental categorical tools used here include lax monoidal functors, which provide a language of compositionality,
as well as sheaf theory, which flexibly captures the crucial notion of time.
\end{abstract}

\thanks{Schultz, Spivak and Vasilakopoulou were supported by AFOSR grant FA9550--14--1--0031 and NASA grant NNH13ZEA001N}

\maketitle

\setcounter{tocdepth}{1}
\tableofcontents

\section{Introduction}

\begin{figure}[h]
\begin{center}
\begin{tikzpicture}[oriented WD, bbx = .5cm, bby =.8ex, bb min width=.5cm, bb port length=2pt, bb port sep=1]
  \node[bb={1}{3}] (X11) {};
  \node[bb={2}{1}, right=1.5 of X11] (X12) {};
  \node[bb={3}{2}, red, below right=of X12] (X13) {};
  \node[bb={2}{2}, fit={(X11) (X12) (X13) ($(X12.north)+(0,2)$) ($(X13.east)+(.5,0)$)}, dashed] (Y1) {};
  \node[bb={2}{1}, below left=6 and 0 of X13] (X21) {};
  \node[bb={0}{2},below left=of X21] (X22) {};
  \node[bb={1}{2}, fit=(X21) (X22), dashed] (Y2) {};
  \node[bb={2}{3}, fit={($(Y1.north)+(0,1)$) ($(Y2.south)-(0,1)$) ($(Y1.west)-(.25,0)$) ($(Y1.east)+(.25,0)$)}] (Z) {};
  \begin{scope}[gray]
  \draw[ar] (Z_in1') to (Y1_in1);
  \draw[ar] (Z_in2') to (Y2_in1);
  \draw (Y2_in1') to[in looseness=2] (X21_in1);
  \draw[ar] (X22_out1) to (X21_in2);
  \draw (X22_out2) to (Y2_out2');
  \draw[ar] (Y2_out2) to (Z_out3');
  \draw (Y1_in1') to (X11_in1);
  \draw (X13_out2) to (Y1_out2');
  \draw (Y1_in2') to[in looseness=2] (X13_in3);
  \draw (X11_out3) to[in looseness=2] (X13_in2);%
  \draw (X12_out1) to (X13_in1);
  \draw (X11_out2) to (X12_in2);
  \draw (Y1_out1) to (Z_out1');
  \draw[ar] (Y1_out2) to (Z_out2');
  \draw (X21_out1) to (Y2_out1');
  \draw[ar] let \p1=(Y2.north east), \p2=(Y1.south west), \n1={\y1+2*\bby}, \n2=\bbportlen in
  	(Y2_out1) to[in=0] (\x1+\n2,\n1) -- (\x2-\n2,\n1) to[out=180] (Y1_in2);
  \draw[ar] let \p1=(X13.north east), \p2=(X12.north west), \n1={\y2+\bby}, \n2=\bbportlen in
  	(X13_out1) to[in=0] (\x1+\n2,\n1) -- (\x2-\n2,\n1) to[out=180] (X12_in1);
  \draw let \p1=(X12.north west), \p2=(X13.north east), \n1={\y1+2*\bby}, \n2=\bbportlen in
  	(X11_out1) to (\x1-2*\n2,\n1) -- (\x2+2*\n2,\n1) to[out=0] (Y1_out1');
  \end{scope}
  \end{tikzpicture}
  \end{center}
  \caption{Compositional analyses facilitate the 
  rearrangement as well as the replacement of internal components}
  \label{fig:compositionall}
\end{figure}
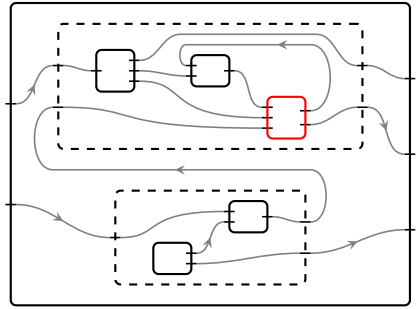

The broad goal of this paper is to suggest a solid categorical framework for understanding and simulating systems of systems.
While the current work is certainly theoretical, our ultimate interest is to
make the case that category theory can be useful for modeling real-world systems, certainly not a novel thesis as we discuss in the related work section. 
We consider broadly defined \emph{open} (dynamical) systems that take in, process, and send out material.
Their complexity level varies greatly; for example, they can be used
to model anything from an electric circuit, to a chemistry experiment, to a robot.
Moreover, they are designed to be interconnected: the material output
by one system is sent to, and received by, another. 
The central idea to which the current work adheres, building on earlier works like \cite{Rupel.Spivak:2013a,Spivak:2013b,Vagner.Spivak.Lerman:2015a},
is that a single system may arise by wiring together any number of component open subsystems; see \cref{fig:compositionall}.
Analyzing a composite system is often intractable, because its complexity is generally exponential in the number of subsystems.
Hence it is often crucial that the analysis be \emph{compositional} \cite{NIST_report} because such analyses can be applied to the subsystems independently, and the results can be composed in a specific sense. This also means that 
the analysis is robust to redesign: improvements can arise from reconfiguring any one of the numerous parts and subparts of large-scale systems (see \cref{fig:compositionall}) at any level of a hierarchy, which can itself be re-structured, and the analyses of unaffected systems remain valid. 

The theory of monoidal categories and operads provides an excellent formalism in which to study compositionality. In our approach,
an object in the symmetric monoidal category $\WW{\cat{C}}$ of $\cat{C}$-\emph{labeled boxes}
and \emph{wiring diagrams} (for some category $\cat{C}$) looks like a box in \cref{fig:compositionall}, thought of as an interface for an open dynamical system with input and output ports through which it can interact with its environment.
A morphism in $\WW{\cat{C}}$ describes how systems can be interconnected together to form new systems: for example, any of the two dotted composite systems in \cref{fig:compositionall} are morphisms in the underlying operad, whereas that picture constitutes an operad composition of a 3-ary (top) and a 2-ary (bottom) morphism producing a 5-ary one. This framework is a cornerstone for the present work and will be described in detail; wiring diagrams become the algebraic operations for combining dynamical systems.


The \emph{typing} category $\cat{C}$ indicates the sort of information or material that each component in a system may send or receive. For example, if $\cat{C}=\Set$, each port is associated to a set of possible inputs/outputs signals. To properly model the flow of such signals,
we must decide whether changes are happening continuously, at discrete time steps, or in some combination; that is how the fundamental issue of time enters our formalization. In particular, for a discrete dynamical system the flow occurs on the ticks of a global clock, whereas for continuous dynamical systems it occurs continuously throughout intervals of time. Furthermore, 
when different types of systems are composed, their notions of time must be faithfully translated into a common one,
keeping in mind that there is no `best' notion of time: the more expressive a language is for describing systems, the more difficult it is to answer questions about them.

The categories, or more precisely toposes, that capture the required variations of time in this paper are sheaves on intervals in $\RR$ or in $\NN$, structures also studied in \cite{Lawvere:1986a}. Such typing categories, herein denoted $\Shv{\Int}$ or $\Shv{\Int_N}$, endow input and output ports with sets of allowed trajectories of information over each (discrete or continuous) interval of time: these trajectories can restrict over smaller intervals of time, and compatible ones can glue to produce a trajectory over the union of intervals. 
There is also a (slice) topos $\Shv{\Int}/\Sync$ encompassing `synchronized continuous time', in which each continuous interval is assigned a phase $\theta\in [0,1)$. If one wishes to combine continuous-time systems with discrete ones that operate on the same clock, the result is in the topos of \emph{synchronous} sheaves; see \cref{rem:whysynchrony}.



Having established the interface and the time framework, our world of interacting open systems is expressed as a \emph{wiring diagram algebra}, namely a lax monoidal functor $F\colon\WW{\cat{C}}\to\Cat$ or equivalently an operad morphism from the underlying operad. For any such functor and $X$ a $\cat{C}$-labeled box, $FX$ is the category of $F$-systems (for example, Moore machines or continuous dynamical systems) with input and output ports determined by the shape of $X$. Given $F$-systems `inhabiting' the boxes of a picture like \cref{fig:compositionall}, the lax monoidal structure along with functoriality of $F$ coherently produce a composite $F$-system that inhabits the outer box. Moreover, any functor $H\colon\cat{C}\to\cat{D}$ induces a strong monoidal functor $\WW{H}\colon\WW{\cat{C}}\to\WW{\cat{D}}$ capable of changing the input/output types. A morphism, then, of wiring diagram algebras is a \emph{monoidal opindexed 1-cell} \cite{Moeller.Vasilakopoulou:2018a}, i.e.\ a monoidal natural transformation
\begin{displaymath}
\begin{tikzcd}[column sep=.7in,row sep=.03in]
\WW{\cat{C}}\ar[dr, "F"]\ar[dd, "\WW{H}"'] \\ \ar[r,phantom,"\Downarrow{\scriptstyle\alpha}"description] & \Cat \\
\WW{\cat{D}}\ar[ur,"G"']
\end{tikzcd}
\end{displaymath}
with components $\alpha_{X}\colon FX\to G(\WW{H}(X))$ that specifically translate $F$-systems to $G$-systems.

Our work is based on a class of very general, span-like wiring diagram algebras, which we simply call \emph{machines}.
A lax monoidal functor $\SpnS_{\cat{C}}\colon\WW{\cat{C}}\to\Cat$ (\cref{prop:C_Span_Algebra}) gives rise to certain open systems called \emph{discrete}, \emph{continuous}, or \emph{synchronous} machines, corresponding to the situation in which the typing category $\cat{C}$ is the topos of discrete, continuous, or synchronous interval sheaves respectively. 
Not only do ordinary discrete and continuous dynamical systems translate into discrete and continuous machines via an embedding of algebras, but also they can ultimately compose to one another since all machines embed into the algebra of synchronous ones, $\ADS_\Sync\colon\WW{\Shv{\Int}/\Sync}\to\Cat$. 

A very important feature of our machines is that---in all time versions---there exist subalgebras of \emph{inertial}, \emph{deterministic}, and \emph{total} machines. In more detail, certain lifting conditions examine whether a system can determine a piece of a future output trajectory, whether a state could evolve in more than one way, or whether a state could potentially `die', when the current state is subject to some input trajectory. The fact that such subcategories form algebras themselves means that the respective definitions are carefully chosen so as to be closed under arbitrary compositions, including feedback. For that, the formalism of $\cat{W}_\cat{C}$-algebras is essential; such a result would not hold if traced monoidal categories were used instead, roughly for the same reason that \cite{Abramsky.Blute.Panangaden:1999a} introduces trace ideals: ``there appears to be a tension between having identities and having compact closed structure'', or in our terminology, identity maps are not inertial. Finally, for all above machine variations, there exist \emph{contracted machines} defined as systems of the same type that comply with a given `contract', expressed as a sub-presheaf. A logic of behavior contracts in a slightly different setting was later given in \cite{Schultz.Spivak:2017a}.

\subsection*{Related work and future goals}

The present work is part of a far more general, ongoing endeavor by many researchers to categorically formalize general processes and exhibit the advantages of such compositional analyses. Our results aim to contribute towards a better understanding of systems, using the theory of monoidal categories to capture nesting and sheaf theory to capture time. In particular, as mentioned earlier, the topos $\Shv{\Int}$ of sheaves on real intervals was also studied by Lawvere in \cite{Lawvere:1986a} for similar purposes; this becomes evident especially in \cref{appendix:DSF}. Notably, all of our models of time assume a common reference time frame for all components in the system: real time progresses for everyone at the same rate. If one wants a more flexible model for non-interacting parallel subsystems (i.e.\ concurrency), or to allow for relativistic effects, $\Int$-sheaves would have to be replaced by a different formalism such as event structures \cite{Winskel:1987a} or action structures \cite{Milner:1996a}.

The parts of this work oriented to discrete dynamical systems connect to the well-established theory of composing finite state automata (Moore machines) or more generally cyber-physical systems, see e.g. \cite{Lee.Seshia:2017}. In our framework, the basic operations are feedback and parallel placement, expressed as morphisms and monoidal product in the wiring diagram category; for example, serial composition (or CASCADE) is derived from those. An advantage of this framework is that every wiring diagram arrangement of automata, like the one in \cref{fig:compositionall}, can be directly given a precise mathematical formulation as their composite. This constitutes a universal way of composing systems, agnostic to the complications of interconnections; we plan to further pursue this point---as well as compositions of Moore machines with different types of systems---according to current challenges in the respective fields. 

Our general system notion is essentially a span inside $\Shv{\Int}$, which is in fact equivalent to a \emph{discrete \Conduche fibration} \cite{Johnstone:1999a} as explained in \cref{sec:B3}. In \cite{Fiore:2000a,Bunge.Fiore:2000a}, the authors present a very similar strategy of expressing a general process via a discrete \Conduche fibration (therein called \emph{unique factorization lifting} functor) over monoids of intervals, also distinguishing between discrete and continuous variations. They also discuss \emph{hybrid} systems which can in fact be modeled as continuous machines, see \cite{Schultz.Spivak:2017a}. Thus there exist strong similarities between the core formalisms of the two works which should be further investigated, while also examining certain differences: first, the distinction between input and output for us is essential---we are interested in the conditions under which deterministic and total systems compose---whereas their systems communicate with the outer world via a single \emph{control} space. Moreover, the two works' scopes differ a great deal: we focus on nesting and composing systems of different types, whereas they focus on the designs of languages and their logic.

Such an investigation is also relevant to connections of the current framework to a more behavioral approach to systems, e.g. \cite{Willems.Polderman:2013a}, which abandons the input/output distinction: our span formalization of machines is already symmetric in that sense. This is also related to machines giving an immediate instance of a \emph{hypergraph category}, recently shown to be equivalent to algebras for an operad of cospans \cite{Fong.Spivak:2018}. On the other hand, total machines in our context seem to be strongly associated to the concept of \emph{open maps} in the theory of bisimulation, see e.g. \cite{Joyal.Nielsen.Winskel:1996}, where again the distinction between input and output is critical.

Our work is also connected to the theory of traced monoidal categories \cite{Joyal.Street.Verity:1996a} and PROPs in the following way. Recall that a \emph{prop} is a symmetric strict monoidal category whose monoid of objects is $(\NN,0,+)$, the free  monoid on one generator. A \emph{wheeled prop} (see e.g.\ \cite{Markl:2009a}) is a traced, symmetric, strict monoidal category whose underlying (cancellable) monoid of objects is again $(\NN,0,+)$. The diagrams in a wheeled prop look very similar to diagrams in this paper like \cref{picture_wiring_diagram}, so it is worth elaborating on their differences. First of all, it was shown in \cite{Spivak.Schultz.Rupel:2016a} that the category of wheeled props is equivalent to that of algebras on the operad $1\text{-}\Cat\mathbf{Cob}$ of oriented 1-cobordisms. The operad of wiring diagrams does indeed compare to that of cobordisms, but the two are not equivalent. 
In more detail, in a wheeled prop there is only one generating object, but we have a generating object for each sheaf; thus one should begin by using a colored wheeled prop. Moreover, wires in our diagrams are allowed to split or terminate; thus one should assume that each color in the prop is equipped with a comonoid structure. Finally, wires in our diagrams are not allowed to feed straight through; this is a sort of \emph{dialectica} condition in the sense of \cite{DePaiva:1990}, which renders it unreasonable to ask for a map from cobordisms to $\WW{}$ 
since the usual axiomatization of traced monoidal categories depends on the identity. Thus while it is true that every wheeled colored prop with a comonoid structure on every object induces a wiring diagram algebra in our sense, the converse does not hold.

Our work is also similar in spirit to that of \cite{Katis.Sabadini.Walters:2000a} on systems with boundary. This work builds on a notion of a \emph{category with feedback}, closely related to traced monoidal categories  but with a built-in notion of delay, used to define an \emph{ifo} (input-feedback-delay) system: examples of those include circuits, algorithms as well as Mealy machines. A significant example is that of spans of graphs, which are shown to model discrete time processes; since graphs are precisely sheaves on discrete intervals $\Shv{\Int}_N$ as explained by \cref{prop:ShIntN_grphs}, our interests fully align in this case. When it comes to continuity, the authors keep the same model by perceiving discrete spaces as infinitesimal motions, whereas our continuous machines live as spans inside a completely different topos---sheaves on real intervals. 

Finally, our notion of machines was motivated by a project with NASA and Honeywell, in which we were tasked with giving semantics to a system of interacting systems, such as the US national airspace system. This work later morphed into \cite{Schultz.Spivak:2017a} which emphasizes the logical aspects, whereas the current paper is geared toward an understanding of control-theoretic issues: when do systems of systems, with components of different types, still exhibit deterministic control based on inputs. The work is ongoing, e.g.\ we hope to understand the semantics of probabilistic structures (internal valuations), so that we can model probabilistic behaviors, such as stochastic processes.

\subsection*{Plan of the paper}

In \cref{ch:wiring_diagrams_and_algebras},
we briefly review monoidal categories and (colored) operads as well as the theory of their algebras, and we describe the category $\WW{\cat{C}}$ of labeled boxes and wiring diagrams for a typing category $\cat{C}$. 
We recall the $\WW{\cat{C}}$-algebras of discrete and continuous dynamical systems, for which wires carry only static information---the set of symbols or the space
of parameters that drive the system---and finally we introduce the general $\cat{C}$-span algebras.

We bring time into the picture in \cref{ch:Int}, where we define sites $\Int$ and $\Int_N$
of continuous and discrete-time intervals respectively, and consider their toposes of sheaves. 
We also provide morphisms from discrete and continuous sheaves to a slice topos $\Shv{\Int}/\Sync$, where $\Sync$
is a \emph{synchronizing} sheaf whose sections can be thought of as phase shifts.

We define abstract dynamical systems, which we call \emph{machines}, in \cref{ch:ADS}.
They form $\WW{\cat{C}}$-algebras, and so do various subclasses such as machines that are total and/or deterministic. For $\cat{C}=\Shv{\Int},\Shv{\Int_N}$ or $\Shv{\Int}/\Sync$ we obtain different time-typed machines, designed to encompass as general systems as possible,
and we also discuss a notion of \emph{safety contracts} for machines. 

In \cref{ch:maps_to_machines} we provide algebra morphisms that translate between almost every
combination of machines discussed previously. In particular, we convert our older, static definitions of discrete and continuous dynamical systems 
into our new language of time-based machines, and also discrete and continuous machines into synchronous machines.
We give sufficient conditions for these translations to preserve totality and determinism. 

Finally, in \cref{appendix:DSF} we give an alternative categorical
construction of the sheaf topos $\Shv{\Int}$, namely as the category of \emph{discrete \Conduche fibrations} over
the monoid of nonnegative real numbers. This equivalence is derived from results in \cite{Johnstone:1999a}
and provides another viewpoint for machines and their variations. 

\subsection*{Acknowledgments}

We greatly appreciate our collaboration with Alberto Speranzon and Srivatsan Varadarajan,
who have helped us to understand how the ideas presented here can be applied in practice
(specifically for modeling the National Airspace System) and who provided motivating examples with which to test and often augment the theory.
We also thank the anonymous reviewers for valuable suggestions; in particular, such a suggestion led to a more abstract formalism
of system algebras, explained in \cref{sec:Span_Systems}.

\section{Wiring diagrams and algebras}\label{ch:wiring_diagrams_and_algebras}

This paper is about interconnecting systems in order to build more complex systems.
The notion of building one object from many is nicely captured using operads and their
algebras: an operad describes the ways 'objects' can be combined, and an operad algebra gives semantics by describing the 'objects' themselves. 
The operads used in this article in fact all underlie monoidal categories. While we assume
familiarity with monoidal categories, we review some key facts and the relationship to operads in \cref{sec:background}.
Standard references for these topics include \cite{Joyal.Street:1993a} and \cite{Leinster:2004a}.
The notion of a multicategory (colored operad) was in fact introduced much earlier by Lambek \cite{Lambek:1969a},
along with the crucial observation that every monoidal category gives rise to such a structure.

In \cref{sec:operad_of_wiring_diagrams} we restrict our interest to a certain symmetric monoidal category
$\WW{}$ (and its associated operad) studied in previous works like \cite{Vagner.Spivak.Lerman:2015a}.
The objects and morphisms in $\WW{}$ classify a certain sort of string diagrams,
which for reasons explained elsewhere \cite{Spivak.Schultz.Rupel:2016a}
we call wiring diagrams. The category of wiring diagram algebras, i.e.\ monoidal functors $\WW{}\to\Cat$,
is the realm where various notions of systems live and can be composed according to their wiring arrangements.
In \cref{sec:DandCDS} we elaborate
on two important examples of those algebras, namely discrete and continuous dynamical systems, whereas
in \cref{sec:Span_Systems} we describe a very general class of span-like algebras, whose
semantics is fundamental for our description of abstract systems in \cref{ch:ADS}.

\subsection{Background on monoidal categories and operads}\label{sec:background}

We denote a monoidal category by $(\cat{V},\otimes,I)$. 
Recall that a \emph{lax monoidal} functor $F\colon(\cat{V},\otimes,I)\to(\cat{W},\otimes,I)$ comes equipped with natural structure morphisms
\begin{displaymath}
F_{c,d}\colon F(c)\otimes F(d)\to F(c\otimes d)
\quad\tn{and}\quad
F_I\colon  I \to F(I)
\end{displaymath}
satisfying well-known coherence axioms.
The functor $F$ is \emph{strong monoidal} if these structure morphisms are isomorphisms. If the
monoidal categories $\cat{V}$ and $\cat{W}$ are moreover symmetric, with $b_{c,d}\colon c\otimes d\simrightarrow
d\otimes c$, then $F$ is \emph{symmetric} if the structure maps appropriately
commute with the symmetry isomorphisms.
A \emph{monoidal} natural transformation $\alpha\colon F\Rightarrow G$
between monoidal functors is an ordinary natural transformation
whose components furthermore commute with the structure morphisms. For detailed descriptions, see e.g.\ \cite{Joyal.Street:1993a}.

We denote by $\SMonCatl$ the 2-category of symmetric monoidal categories,
symmetric lax monoidal functors and monoidal natural transformations, and by
$\SMonCats$ the corresponding 2-category with strong monoidal functors.

Two fundamental examples of cartesian monoidal categories are
$(\Set,\times,\singleton)$ and $(\Cat,\times,\1)$. In our context, we will
often refer to symmetric lax monoidal functors $\cat{V}\to\Cat$
as ($\Cat$-valued)
$\cat{V}$-\emph{algebras}, and to monoidal natural transformations between them as
$\cat{V}$-\emph{algebra maps}; hence we denote $\cat{V}\alg\coloneqq\SMonCatl(\cat{V},\Cat)$.
The terminology comes from the closely-connected world of operads, to which we next turn.

Recall that a \emph{colored operad}, or \emph{multicategory},
$\cat{P}$ consists of a set of objects (colors) $\Ob\cat{P}$, a hom-set of $n$-ary operations
$\cat{P}(c_1,\ldots,c_n;c)$ for each $(n+1)$-tuple of objects, an identity operation
$\id_c\in\cat{P}(c;c)$, and a composition formula 
\begin{displaymath}
\cat{P}(c_1,\ldots,c_n;c)\times\cat{P}(c_{11},\ldots,c_{1k_1};c_1)\times\cdots\times\cat{P}(c_{n1},\ldots,c_{nk_n};c_n)\to
\cat{P}(c_{11},\ldots,c_{nk_n};c)
\end{displaymath}
that can be visualized \cite[Fig.~2-B]{Leinster:2004a} as
\begin{equation}\label{eq:operadcomp}
 \begin{tikzpicture}[triangle/.style = {fill=gray!20, regular polygon, regular polygon sides=3}]
\path (0,0) node [triangle,draw,shape border rotate=-90,inner sep=0pt,label=178:$\vdots$] (a) {$\theta_n$}
(0,4) node [triangle,draw,shape border rotate=-90,inner sep=0pt,label=178:$\vdots$] (b) {$\theta_1$}
(3,2) node [triangle,draw,shape border rotate=-90,label=135:$c_1$,label=230:$c_n$,label=178:$\vdots$] (c) {$\theta$}
(9,2) node [triangle,draw,shape border rotate=-90,inner sep=-25pt,label=178:$\vdots$] (d) {$\qquad\theta\circ(\theta_1,\ldots,\theta_n)$};
\draw [-] (a) .. controls +(right:2cm) and +(left:1cm).. (c.220); 
\draw [-] (b) .. controls +(right:2cm) and +(left:1cm).. (c.140);
\draw [-] (d) to node [above] {$c$} (11.7,2);
\draw [-] (c) to node [above] {$c$} (4,2);
\draw [-] (7.2,3.5) to node [above] {$c_{11}$} (7.85,3.5);
\draw [-] (7.2,.5) to node [below] {$c_{nk_n}$} (7.85,.5);
\draw [-] (-.85,4.3) to node [above] {$c_{11}$} (-.3,4.3);
\draw [-] (-.85,3.6) to node [below] {$c_{1k_1}$} (-.3,3.6);
\draw [-] (-.85,0.3) to node [above] {$c_{n1}$} (-.3,0.3);
\draw [-] (-.85,-0.4) to node [below] {$c_{nk_n}$} (-.3,-0.4);
\node () at (5.5,2) {$\mapsto$};
\end{tikzpicture}
\end{equation}
This data is subject to associativity and unitality axioms,
and the operad is moreover \emph{symmetric} if there are compatible permutation actions on the hom-sets,
$\cat{P}(c_1,\ldots,c_n;c)\simrightarrow\cat{P}(c_{\sigma(1)},\ldots,c_{\sigma(n)};c)$.
An \emph{operad functor} $F\colon\cat{P}\to\cat{P'}$ consists of a mapping on objects,
and a mapping on hom-sets $\cat{P}(c_1,\ldots,c_n;c)\to\cat{P}'(Fc_1,\ldots,Fc_n;Fc)$
that preserves composition, symmetries and identities. An \emph{operad transformation}
$\alpha\colon F\Rightarrow G$ 
consists of unary operations $\alpha_c\in\cat{P'}(Fc;Gc)$
that are compatible with respect to the composition formula. Therefore we obtain
a 2-category $\SOpd$, often denoted $\SMulticat$ in the literature,
of symmetric colored operads; see e.g.\ \cite{Hermida:2000a,Leinster:2004a}.

There exists a 2-functor, called the \emph{underlying operad} functor as in \cite[Example 2.1.3]{Leinster:2004a},
\begin{equation}\label{eq:underlying_operad}
\undop\colon\SMonCatl\rightarrow\SOpd
\end{equation}
mapping each symmetric monoidal category $\cat{V}$ to the operad with
$\Ob(\undop\cat{V})\coloneqq\Ob\cat{V}$ and
$\undop\cat{V}(c_1,\ldots,c_n;c)\coloneqq\cat{V}(c_1\otimes\cdots\otimes c_n,c)$.
Note that for $c=c_1\otimes\cdots\otimes c_n$, there is a unique morphism in $\undop\cat{V}$
corresponding to $\id_c$, called the \emph{universal morphism} for
$(c_1,\dots,c_n)$. To each lax monoidal functor $F\colon\cat{V}\to\cat{W}$, we can assign an operad functor
$\undop F$ with the same mapping on objects. On morphisms, it sends some $f\colon c_1\otimes\cdots\otimes c_n\to c$ in $\cat{V}$ to the composite
\[
Fc_1\otimes\cdots\otimes Fc_n\xrightarrow{F_{c_1,\cdots,c_n}}F(c_1\otimes\cdots\otimes c_n)\xrightarrow{Ff}Fc
\]
in $\cat{W}$.
Lastly, a monoidal natural transformation $\alpha$ involves exactly the data needed to define the operad transformation $\undop\alpha$.

The functor $\undop$ is clearly faithful, but it is also full:
for any operad map $G\colon\undop\cat{V}\to\undop\cat{W}$, we can define a functor
$F\colon\cat{V}\to\cat{W}$ via the action of $G$ on objects and unary morphisms. Its lax monoidal structure is given by the action of $G$ on universal morphisms.
Moreover, there is also a bijection between the respective 2-cells,
which exhibits $\undop$ as a fully faithful 2-functor; see \cite[Remark 9.5]{Hermida:2000a}
for relevant results. Not all operads underlie monoidal categories---the functor $\undop$ is not essentially surjective---however, all operads
of interest in this paper do happen to be in the image of $\undop$. 

For typographical reasons, we denote
the operad $\undop\Cat$ underlying the cartesian monoidal category of categories
using the same symbol $\Cat$, when no confusion arises. If $\cat{P}$ is an operad,
a ($\Cat$-valued) $\cat{P}$-\emph{algebra} is an operad functor $\cat{P}\to\Cat$.
The category of $\cat{P}$-algebras is denoted $\cat{P}\alg\coloneqq\SOpd(\cat{P},\Cat)$,
see \cite[Def 2.1.12]{Leinster:2004a}. This nomenclature agrees with that above, since
for any monoidal category $\cat{V}$,
the fully faithful $\undop$ induces an
isomorphism between the corresponding categories of algebras
\begin{equation}\label{eq:algebras_iso}
\cat{V}\alg\simrightarrow(\undop\cat{V})\alg.
\end{equation}

\begin{remark}\label{rem:everything_is_pseudo}
In this article, we will in fact work with $\cat{W}$-\emph{pseudoalgebras}
for various monoidal categories $\WW{}$, rather than ordinary algebras. This amounts to viewing
$\cat{W}$ as a monoidal 2-category (with trivial 2-cells), and considering
\emph{weak monoidal pseudofunctors} $A\colon(\WW{},\otimes,I)\to(\Cat,\times,\1)$. More explicitly,
for any pair of maps $\psi,\phi$ in $\WW{}$, we have only a natural isomorphism of functors $A(\psi\circ\phi)\cong A(\psi)\circ A(\phi)$
rather than a strict equality; similarly for the identities, along with coherence axioms.
Moreover, the lax monoidal structure is pseudo, meaning that any previously commutative diagram of the axioms
now only commutes up to natural isomorphism (with certain
coherence conditions tying them together). 
In fact, $\cat{W}$-algebras are really \emph{monoidal indexed categories} in the sense of \cite{Moeller.Vasilakopoulou:2018a}.

More details of such structures can be found e.g.\ in \cite{Day.Street:1997a},
where $\WMonHom(\cat{A},\cat{B})$ denotes the category of weak monoidal pseudofunctors and monoidal pseudonatural transformations, for any two monoidal 2-categories $\cat{A},\cat{B}$.
Due to coherence theorems, there is an equivalence $\cat{W}\alg\simeq\WMonHom(\cat{W},\Cat)$
and so we safely omit such details in the presentation below,
to keep the already-technical material more readable. Thus we will speak of the pseudoalgebra
$A$ as above simply as a $\WW{}$-algebra and the pseudomaps between them simply as $\WW{}$-algebra morphisms.
\end{remark}

\subsection{The operad of labeled boxes and wiring diagrams}\label{sec:operad_of_wiring_diagrams}

Operads (or monoidal categories) of various `wiring diagram shapes' have been considered in
\cite{Spivak:2013b,Rupel.Spivak:2013a}. More recently, \cite{Spivak.Schultz.Rupel:2016a}
showed a strong relationship between the category of algebras on a certain operad (namely
$\Cob$, the operad of oriented 1-dimensional cobordisms) and the category of traced monoidal
categories; results along similar lines were proven in \cite{Fong:2016a}. The operads $\WW{}$ we use here
are designed to model what could be called `cartesian traced categories without identities'.%
\footnote{The reason to study traced categories without identities is that there is no
`identity' dynamical system. Attempting to define such a thing, one loses the trace
structure; similar problems arise throughout computer science and engineering applications:
see \cite{Ageron:1996a,Selinger:1997a}. A solution to this problem was described in
\cite{Abramsky.Blute.Panangaden:1999a}, using trace ideals. The present operadic approach is another solution,
which appears to be roughly equivalent.} 

For any category $\cat{C}$, a \emph{$\cat{C}$-typed finite set} is a finite set $X$
together with a function $\tau\colon X\to\Ob\cat{C}$ assigning to every
element $x\in X$ an object $\tau(x)\in\cat{C}$, called its \emph{type}.
We often elide the typing function in notation when it is implied from the context, writing $X$ rather than $(X,\tau)$.
The $\cat{C}$-typed finite sets form a category $\TFS_{\cat{C}}$,
where morphisms are functions $f\colon X\to X'$ that respect the types, i.e.\
$\tau(fx)=\tau(x)$; we call such an $f$ a \emph{$\cat{C}$-typed function}.%
\footnote{
The category $\TFS_{\cat{C}}$ of typed finite sets in $\cat{C}$ is
closely related to the \emph{finite family fibration}
$\Fam(\cat{C}\op)\to\FinSet$ (see \cite[Definition 1.2.1]{Jacobs:1999a});
namely $\TFS_{\cat{C}}$ is the category of cartesian arrows in $\Fam(\cat{C}\op)$, or equivalently cartesian arrows in $\Fam(\cat{C})$.
Replacing $\TFS_{\cat{C}}$ with $\Fam(\cat{C}\op)$ throughout the paper would yield similar
definitions and results; for example, taking product of types \cref{eq:dependent_product} is
well-defined in this context, $\prod\colon\Fam(\cat{C}\op)\op\to\cat{C}$.
}
Thus we have
\begin{displaymath}
\TFS_{\cat{C}}\cong\FinSet/{\Ob\cat{C}}
\end{displaymath}
and $\cat{C}$ is called the \emph{typing} category.
Notice that $(\TFS_\cat{C},+,\emptyset)$ is a
cocartesian monoidal category via
\[(X,\tau)+(X',\tau')=\big(X+X',(\tau,\tau')\big)\]
where $X+X'$ is the disjoint union of $X$, $X'$ and $(\tau,\tau')$ denotes the co-pairing of functions.
Moreover, any functor $F\colon\cat{C}\to\cat{D}$ induces
a (symmetric strong monoidal) functor $\TFS_F\colon\TFS_\cat{C}\to\TFS_\cat{D}$, sending $(X,\tau)$ to
$(X,F\circ\tau)$. We have thus constructed a functor $\TFS_{(\Inp)}\colon\Cat\to\SMonCats$.

If $\cat{C}$ has finite products, then we can assign to each $\cat{C}$-typed finite set $(X,\tau)$
the product of all types of its elements, denoted $\dProd{X}\in\cat{C}$:
\begin{equation}\label{eq:dependent_product}
\dProd{X}\coloneqq\prod_{x\in X}\tau(x)
\end{equation}
This induces a functor
$\dProd{(\Inp)}_\cat{C}\colon\TFS_{\cat{C}}\op\to\cat{C}$, which is strong monoidal because there are
isomorphisms $\dProd{X_1}\times\dProd{X_2}\cong\dProd{X_1+X_2}$ and
$\dProd{\emptyset}\cong1$. In fact, if $\FPCat$ is the category of
finite-product categories and functors which preserve them, and $i\colon\FPCat\to\SMonCats$
realizes a finite-product category
with its canonical cartesian monoidal structure,
then the functors $\dProd{(\Inp)}_\cat{C}$ are the components of a pseudonatural transformation
\begin{equation}\label{eq:producttransf}
\begin{tikzcd}[row sep=.1in]
& |[alias=domA]|\Cat\ar[dr,bend left=15,"\TFS_{(-)}\op"] & \\
\FPCat\ar[rr,bend right,"i"' codA]\ar[ur,bend left=15,>->] && \SMonCats
\twocellA{\dProd{(\Inp)}}
\end{tikzcd}
\end{equation}
where the ordinary category $\FPCat$ and functor $\TFS$ are viewed as 2-structures
with the trivial 2-cells and mappings respectively.

\begin{example}\label{ex:types_examples}
Typed finite sets will be heavily used in the context of wiring diagrams; here we consider a few typing categories
$\cat{C}$, whose objects will serve as types of elements later.
In general, if $(X,\tau)$ is a typed finite set,
we think of each element $x\in X$ as a \emph{port}, while its associated type $\tau(x)\in\Ob\cat{C}$ specifies
what sort of information passes through that port.
\begin{compactenum}
\item When $\cat{C}=\Set$, each port carries a \emph{set} of possible signals;
there is no notion of time. Ports of this type are relevant when studying
discrete dynamical systems.
\item When $\cat{C}=\Man$, (a small category equivalent to) the category of second-countable smooth manifolds
and smooth maps between them, each port carries a manifold of possible signals, again without a notion of time.
Ports of this type are relevant when specifying continuous dynamical systems:
machines that behave according to (ordinary) differential equations on manifolds.
\item\label[item]{TFS_Int} When $\cat{C}=\Shv{\Int}$,
the category of \emph{interval sheaves} defined in \cref{sec:interval_site_Int},
each port carries a very general kind of time-based signal.
The two types of dynamical systems described above in terms of
sets and manifolds can be translated into this language, at which point the time-based dynamics itself becomes evident.
The topos of interval sheaves will be central in our construction of machines.
\end{compactenum}
\end{example}

We are now in position to define a symmetric monoidal category
$(\WW{\cat{C}},\oplus,0)$ for any typing category
$\cat{C}$, see \cite[\S 3]{Vagner.Spivak.Lerman:2015a} or \cite[\S 3.3]{Spivak:2015b}.
Its objects will be called \emph{$\cat{C}$-labeled boxes}; they are pairs $X=(\inp{X},\out{X})\in\TFS_{\cat{C}}\times\TFS_{\cat{C}}$ of $\cat{C}$-typed finite sets. We can picture such an $X$ as
\begin{equation}\label{eq:box}
\begin{tikzpicture}[oriented WD, bbx=.1cm, bby =.1cm, bb port sep=.15cm]
	\node [bb={3}{3}] (X) {$X$};
	\draw[label]
		node[left=.1 of X_in1]  {$a_1$}
		node[left=.1 of X_in2]  {$\dotso$}
		node[left=.1 of X_in3]  {$a_m$}
		node[right=.1 of X_out1] {$b_1$}
		node[right=.1 of X_out2]  {$\dotso$}
		node[right=.1 of X_out3] {$b_n$};
\end{tikzpicture}
\end{equation}
Here, $\inp{X}=\{a_1,\ldots,a_m\}$ is the set of input ports, and $\out{X}=\{b_1,\ldots,b_n\}$
is the set of output ports. Each port $a_i$ or $b_j$ comes with its associated $\cat{C}$-type $\tau(a_i)$ or $\tau(b_j)\in\cat{C}$.

A morphism $\phi\colon X\to Y$ in $\WW{\cat{C}}$ is called a \emph{wiring diagram}. It consists
of a pair
\begin{equation}\label{eq:wiring_diag}
\left\{\begin{array}{l}
\inp{\phi}\colon\inp{X}\to\out{X}+\inp{Y} \\
\out{\phi}\colon\out{Y}\to\out{X}\end{array}\right.
\end{equation}
of $\cat{C}$-typed functions
which express `which port is fed information by which'.
Graphically, we can picture such a map $\phi$, going from the inside box $X$ to the outside box $Y$, as
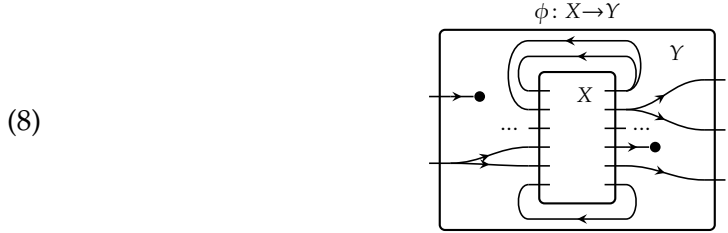
\begin{equation}\label{eq:wiringdiagpic}
\begin{tikzpicture}[oriented WD,baseline=(Y.center), bbx=2em, bby=1.2ex, bb port sep=1.2]
\node[bb={6}{6}] (X) {};
\node[bb={2}{3}, fit={($(X.north east)+(0.7,1.7)$) ($(X.south west)-(.7,.7)$)}] (Y) {};
\node [circle,minimum size=4pt, inner sep=0, fill] (dot1) at ($(Y_in1')+(.5,0)$) {};
\node [circle,minimum size=4pt, inner sep=0, fill] (dot2) at ($(X_out4)+(.5,0)$) {};
\draw[ar] (Y_in1') to (dot1);
\draw[ar] (X_out4) to (dot2);
\draw[ar] (Y_in2') to (X_in5);
\draw[ar] (Y_in2') to (X_in4);
\draw[ar] (X_out5) to (Y_out3');
\draw[ar] (X_out2) to (Y_out1');
\draw[ar] (X_out2) to (Y_out2');
\draw[ar] let \p1=(X.north west), \p2=(X.north east), \n1={\y1+\bby}, \n2=\bbportlen in
	(X_out1) to[in=0] (\x2+\n2,\n1) -- (\x1-\n2,\n1) to[out=180] (X_in1);
\draw[ar] let \p1=(X.north west), \p2=(X.north east), \n1={\y1+2*\bby}, \n2=\bbportlen in
	(X_out1) to[in=0] (\x2+\n2,\n1) -- (\x1-\n2,\n1) to[out=180] (X_in2);
\draw[ar] let \p1=(X.south west), \p2=(X.south east), \n1={\y1-\bby}, \n2=\bbportlen in
	(X_out6) to[in=0] (\x2+\n2,\n1) -- (\x1-\n2,\n1) to[out=180] (X_in6);	
\draw [label] node at ($(Y.north east)-(.5cm,.3cm)$) {$Y$}
              node at ($(X.north east)-(.4cm,.3cm)$) {$X$}
		      node[left=.1 of X_in3]  {$\dotso$}
		      node[right=.1 of X_out3] {$\dotso$}
		      node[above=of Y.north] {$\phi\colon X\to Y$}
	  ;
\end{tikzpicture}
\end{equation}
The identity wiring diagram $1_X=(\inp{1}_X\colon\inp{X}\to\out{X}+\inp{X},\out{1}_X\colon\out{X}\to\out{X})$
is given by the coproduct inclusion and identity respectively. Given another wiring diagram
$\psi=(\inp{\psi},\out{\psi})\colon Y\to Z$, their composite $\psi\circ\phi\eqqcolon\omega=
(\inp{\omega},\out{\omega})$ is given by 
\begin{gather*}
\inp{\omega}\colon\inp{X}\To{\inp{\phi}}\out{X}{+}\inp{Y}\To{1+\inp{\psi}}\out{X}{+}\out{Y}{+}\inp{Z}
\To{1+\out{\phi}+1}\out{X}{+}\out{X}{+}\inp{Z}\To{\nabla+1}\out{X}{+}\inp{Z}\\
\out{\omega}\colon\out{Z}\To{\out{\psi}}\out{Y}\To{\out{\phi}}\out{X}
\end{gather*}
Associativity and unitality can be verified; see \cite[Def 3.10]{Spivak:2015b} for details.

There is a monoidal structure on $\WW{\cat{C}}$ as follows. For two labeled boxes $X_1=(\inp{X_1},\out{X_1})$ and $X_2=(\inp{X_2},\out{X_2})$, their tensor product
$X_1\oplus X_2$ is defined to be $(\inp{X_1}+\inp{X_2},\out{X_1}+\out{X_2})$ and
amounts to parallel composition
\begin{displaymath}
\begin{tikzpicture}[oriented WD,baseline=(Y.center), bbx=1.3em, bby=1ex, bb port sep=.06cm]
 \node[bb={3}{3}] (X1) {};
 \node[bb={3}{3},below =.5 of X1] (X2) {};
 \node[fit=(X1)(X2),draw] {};
 \draw[label] 
              node at ($(X1.west)+(1,0)$) {$X_1$}
              node at ($(X2.west)+(1,0)$) {$X_2$}
              node[left=.1 of X1_in2]  {$\dotso$}
              node[right=.1 of X1_out2]  {$\dotso$}
              node[left=.1 of X2_in2]  {$\dotso$}
              node[right=.1 of X2_out2]  {$\dotso$};
 \draw (X1_in1) -- (-2.5,1.7);
 \draw (X1_out1) -- (2.5,1.7);
 \draw (X1_in3) -- (-2.5,-1.7);
 \draw (X1_out3) -- (2.5,-1.7);
  \draw (X2_in1) -- (-2.5,-5.6);
 \draw (X2_out1) -- (2.5,-5.6);
 \draw (X2_in3) -- (-2.5,-9.1);
 \draw (X2_out3) -- (2.5,-9.1);
 \end{tikzpicture}
\end{displaymath}
of boxes, viewed as a new box. The monoidal unit is $0=(\emptyset,\emptyset)$.
We refer to $(\WW{\cat{C}},\oplus,0)$ as the symmetric monoidal
category of \emph{$\cat{C}$-labeled boxes and wiring diagrams}

An arbitrary functor $F\colon\cat{C}\to\cat{D}$ induces a symmetric strong monoidal
functor $\WW{F}\colon\WW{\cat{C}}\to\WW{\cat{D}}$. It maps objects $X=((\inp{X},\inp{\tau}),(\out{X},\out{\tau}))$
to $\WW{F}(X)=((\inp{X},F\scirc\inp{\tau}),(\out{X},F\scirc\out{\tau}))$ which consist of the same finite set of input and output ports,
but with new typing functions defined by $F$. Together with the identity-like
action of $\WW{F}$ on morphisms, we obtain a functor
\begin{equation}\label{eq:W_}
\WW{(\Inp)}\colon\Cat\to\SMonCats.
\end{equation}
Even though this is defined on arbitrary categories and functors,
it turns out that in many of the examples that follow, using products \cref{eq:dependent_product}
and pullbacks is essential; see also \cref{rem:intermediate_wiring_category}.
Thus we often restrict to finitely complete categories
$\FCCat\rightarrowtail\Cat$ as the domain of $\WW{(\Inp)}$, e.g.\ see
\cref{prop:Mch_FCCat}.

Following \cref{ex:types_examples}, some wiring diagram categories we will use include $\WW{\Set},\WW{\Euc}$ and $\WW{\Shv{\Int}}$. Regarding induced functors between them, consider the example of 
the limit-preserving $U\colon\Euc\to\Set$ which maps a Euclidean space to its underlying
set. This naturally induces a symmetric strong monoidal functor
\begin{equation}\label{eq:W_U}
\WW{U}\colon\WW{\Euc}\to\WW{\Set}
\end{equation}
between the respective categories of wiring diagrams, important for \cref{DDS_to_CDS}.

If we apply the underlying operad functor \cref{eq:underlying_operad} to any monoidal category $\WW{\cat{C}}$,
we obtain the \emph{operad of wiring diagrams} $\undop\WW{\cat{C}}$.
An object, or color, is a $\cat{C}$-labeled box \cref{eq:box}, whereas for example a $5$-ary morphism $\phi\colon X_1,\ldots,X_5\to Y$
in $\undop\WW{\cat{C}}$ can be drawn like
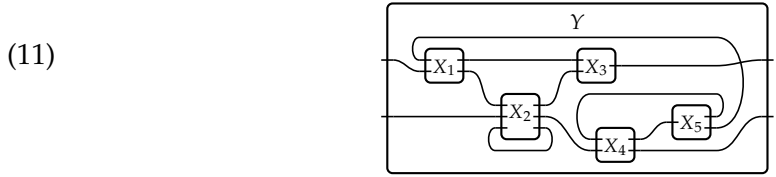
\begin{equation}\label{picture_wiring_diagram}
\begin{tikzpicture}[oriented WD, bb min width =.5cm, bbx=.5cm, bb port sep =1,bb port length=.08cm, bby=.15cm,baseline=.5]
\node[bb={2}{2},bb name = {\tiny$X_1$}] (X11) {};
\node[bb={3}{3},below right=of X11,bb name = {\tiny$X_2$}] (X12) {};
\node[bb={2}{1},above right=of X12,bb name = {\tiny$X_3$}] (X13) {};
\draw (X11_out1) to (X13_in1);
\draw (X11_out2) to (X12_in1);
\draw (X12_out1) to (X13_in2);
\node[bb={2}{2}, below right = -1 and 1.5 of X12, bb name = {\tiny$X_4$}] (X21) {};
\node[bb={1}{2}, above right=-1 and 1 of X21,bb name = {\tiny$X_5$}] (X22) {};
\draw (X21_out1) to (X22_in1);
\draw let \p1=(X22.north east), \p2=(X21.north west), \n1={\y1+\bby}, \n2=\bbportlen in
         (X22_out1) to[in=0] (\x1+\n2,\n1) -- (\x2-\n2,\n1) to[out=180] (X21_in1);
\node[bb={2}{2}, fit = {($(X11.north east)+(-1,3)$) (X12) (X13) ($(X21.south)$) ($(X22.east)+(.5,0)$)}, bb name ={\scriptsize $Y$}] (Z) {};
\draw (Z_in1') to (X11_in2);
\draw (Z_in2') to (X12_in2);
\draw (X12_out2) to (X21_in2);
\draw let \p1=(X22.south east),\n1={\y1-\bby}, \n2=\bbportlen in
  (X21_out2) to (\x1+\n2,\n1) to (Z_out2');
 \draw let \p1=(X12.south east), \p2=(X12.south west), \n1={\y1-\bby}, \n2=\bbportlen in
  (X12_out3) to[in=0] (\x1+\n2,\n1) -- (\x2-\n2,\n1) to[out=180] (X12_in3);
\draw let \p1=(X22.north east), \p2=(X11.north west), \n1={\y2+\bby}, \n2=\bbportlen in
  (X22_out2) to[in=0] (\x1+\n2,\n1) -- (\x2-\n2,\n1) to[out=180] (X11_in1);
\draw let \p1=(X13_out1), \p2=(X22.north east), \n2=\bbportlen in
 (X13_out1) to (\x1+\n2,\y1) -- (\x2+\n2,\y1) to (Z_out1');
\end{tikzpicture}
\end{equation}
Notice that the boxes are the objects and not the morphisms in the underlying operad $\cat{O}\WW{\cat{C}}$, as the earlier triangles representation for $n$-ary operad morphisms might suggest; therefore operadic composition as in \cref{eq:operadcomp} in this framework corresponds to a zoomed-in picture of a box three layers `deep', as in \cref{fig:compositionall}.

In what follows, our goal is to model various processes
as objects of $(\undop\WW{\cat{C}})\alg$, i.e.\ algebras for this operad.
Due to the isomorphism \cref{eq:algebras_iso} between algebras
for a monoidal category and for its underlying operad, we can identify such an operad
algebra with a lax monoidal functor from $(\WW{\cat{C}},+,0)$ to $(\Cat,\times,\1)$,
namely a $\WW{\cat{C}}$-algebra.
Given such an algebra $F\colon\WW{\cat{C}}\to\Cat$ and a $\cat{C}$-labeled box $X$ \cref{eq:box}, we refer to the objects of
the category $F(X)$ as \emph{inhabitants} of $X$. An algebra provides semantics to
the boxes, examples of which we will encounter at \cref{sec:DandCDS} as well as \cref{ch:ADS}.
Whereas the formal description of the algebras will be completely
in terms of lax monoidal functors, having an associated operadic
interpretation provides meaningful pictorial representations like \cref{picture_wiring_diagram}.

In more detail, given inhabitants of each inside box and their arrangement $\phi$,
an algebra $F$ constructs an inhabitant of the outside box,
thus associating a sort of \emph{composition formula} to $\phi$.
Indeed, for any algebra $F\colon\WW{\cat{C}}\to\Cat$, the composite
functor
\begin{equation}\label{eq:inhabitants_composition}
F(X_1)\times\cdots\times F(X_5)\xrightarrow{F_{X_1,\ldots,X_5}}F(X_1+\cdots+X_5)\xrightarrow{F(\phi)}F(Y)
\end{equation}
performs the following two steps. The $F$-inhabitants of five boxes $X_1,
\ldots,X_5$ are first combined---using the lax structure morphisms of $F$---to
form a single inhabitant of the parallel composite box $X\coloneqq X_1+\cdots+X_5$. The wiring diagram $\phi$ can then be considered 1-ary,
and the functor $F(\phi)$ converts the inhabitant of $X$ to an inhabitant of the outer box $Y$. 

Finally, there exists a category with objects symmetric lax monoidal functors $F\colon\cat{W}_\cat{C}\to\Cat$%
\footnote{In fact, $F$ should be a symmetric weak monoidal pseudofunctor and $\alpha\colon F\to G\WW{F}$ should be a monoidal pseudonatural transformation; we sweep these details under the rug as discussed in \cref{rem:everything_is_pseudo}.}
with domains categories of $\cat{C}$-labeled boxes and wiring diagrams, and morphisms
\begin{equation}\label{eq:morphWDalg}
\begin{tikzcd}[column sep=.7in,row sep=.1in]
\WW{\cat{C}}\ar[dr, "F"]\ar[dd, "\WW{F}"'] \\ \ar[r,phantom,"\Downarrow{\scriptstyle\alpha}"description] & \Cat \\
\WW{\cat{D}}\ar[ur,"G"']
\end{tikzcd}
\end{equation}
where $\WW{F}$ is the symmetric strong monoidal functor induced by some functor $F\colon\cat{C}\to\cat{D}$ as in \cref{eq:W_},
and $\alpha$ is a monoidal natural transformation.
Denote this category by $\WDalg$%
\footnote{This is a subcategory of $\mathbf{MonOpICat}$ of the tensor objects in (op)indexed categories, see \cite{Moeller.Vasilakopoulou:2018a}.} and notice there is an obvious forgetful functor $\WDalg\To{\mathrm{dom}}\SMonCats$.
All of our algebras examples like discrete and continuous dynamical systems from \cref{sec:DandCDS}, $\SpnS_\cat{C}$
from \cref{prop:C_Span_Algebra} and machines from \cref{ch:ADS} are objects of $\WDalg$,
whereas \cref{ch:maps_to_machines} describes various maps between them, which allow the translation of one
system kind to another.

\begin{remark}\label{rem:intermediate_wiring_category}
Suppose the typing category $\cat{C}$ has finite products and consider a $\cat{C}$-labeled box $X=(\inp{X},\out{X})$ as in \cref{eq:box}.
By forming the products $\Xin=\prod_{x\in\inp{X}}\tau(x)$ and $\Xout=\prod_{x\in\out{X}}\tau(x)$,
we can associate to the entire input side (all its ports) a single object $\Xin\in\cat{C}$,
and similarly associate $\Xout\in\cat{C}$ to the entire output side. 
Since the functor $\dProd{(\;\cdot\;)}_\cat{C}$ \cref{eq:producttransf}
is contravariant and strong monoidal,
a wiring diagram $\phi=(\inp{\phi},\out{\phi})$ \cref{eq:wiring_diag} induces a pair of morphisms in $\cat{C}$
\begin{equation}\label{eq:wiring_diag_Int}
\left\{\begin{array}{l}
\fin\colon\Yin\times\Xout\to\Xin \\
\fout\colon\Xout\to\Yout\end{array}\right.
\end{equation}
Intuitively, these now describe the direction of information flow in the wiring diagram:
in \cref{eq:wiringdiagpic} the box $X$ receives information only from the input of $Y$ as
well as from $X$'s output (feedback operation), whereas information that exits $Y$ only came from
the output of $X$. 

Through all examples of $\WW{\cat{C}}$-algebras $F\colon\WW{\cat{C}}\to\Cat$ found in this paper,
a pattern arises: $\cat{C}$ will be a finitely complete category, and we always begin by taking the product
of types. That is, all of our example algebras $F$ factor as
\begin{displaymath}
\begin{tikzcd}[row sep=.2in,column sep=.2in]
\WW{\cat{C}}\ar[rr,"F"]\ar[dr,"{X\mapsto\dProd{X}}"'] && \Cat \\
& \cat{\ol{W}}_\cat{C}\ar[ur,"\ol{F}"'] &
\end{tikzcd}
\end{displaymath}
where $\cat{\ol{W}}_\cat{C}$ is the monoidal category with objects $\Ob(\cat{C}\times\cat{C})$,
and morphisms $(X_1,X_2)\to(Y_1,Y_2)$ pairs of $\cat{C}$-morphisms
\begin{displaymath}
\left\{\begin{array}{l}
\phi_1\colon Y_1\times X_2\to X_1 \\
\phi_2\colon X_2\to Y_2\end{array}\right.
\end{displaymath}
with appropriate composition, identities and monoidal structure.
The map $\WW{\cat{C}}\to\cat{\ol{W}}_{\cat{C}}$ sends $X=(\inp{X},\out{X})\mapsto(\Xin,\Xout)$.
We will not mention $\cat{\ol{W}}_{\cat{C}}$ again, though we often slightly abuse notation by writing $F$ in place of $\ol{F}$.

The above description of $\cat{\ol{W}}_\cat{C}$ is reminiscent of compositional game theory \cite{Hedges:2018a}
and there also seem to be connections with bilenses and the Dialectica category \cite{DePaiva:1990}; such considerations are in the center
of future research goals.
\end{remark}

\subsection{Discrete and continuous dynamical systems}\label{sec:DandCDS}

As our primary examples, we consider two well-known classes of dynamical systems, namely \emph{discrete} and
\emph{continuous}, which have already been studied within the context of the operad of wiring diagrams
in previous works; see e.g.\ \cite{Rupel.Spivak:2013a,Spivak:2015b}.
A broad goal of this paper is to group these, as well as other
notions of systems, inside a generalized framework. This will be accomplished by constructing algebra maps
from these special case systems to the new, abstracted ones, described in \cref{ch:maps_to_machines}.

\begin{definition}\label{def:DDS}
A \emph{discrete dynamical system} with \emph{input set} $A$ and \emph{output set} $B$
consists of a set $S$, called the \emph{state set}, together with two functions
\begin{displaymath}
\begin{cases}
\upd{f}\colon A\times S\to S \\
\rdt{f}\colon S \to B
\end{cases}
\end{displaymath}
respectively called \emph{update} and \emph{readout functions},
which express the transition operation and the produced output of the machine. We refer to $(S,\upd{f},\rdt{f})$ as an \emph{$(A,B)$-discrete
dynamical system} or \emph{$(A,B)$-DDS}. If additionally an element $s_0\in S$ is chosen, we refer to $(S,s_0,\upd{f},\rdt{f})$
as an \emph{initialized $(A,B)$-DDS}.%
\footnote{
Note that if $A,B,S$ are finite, what we have called an initialized $(A,B)$-DDS is often called a \emph{Moore machine}.
}

A \emph{morphism} $(S_1,\upd{f_1},\rdt{f_1})\to(S_2,\upd{f_2},\rdt{f_2})$ of $(A,B)$-discrete dynamical systems
is a function $h\colon S_1\to S_2$ which commutes with the
update and readout functions elementwise, i.e.\ $h(\upd{f_1}(s,x))=\upd{f_2}(hs,x)$
and $\rdt{f_1}(s)=\rdt{f_2}(hs)$.
\end{definition}

The resulting category is denoted by $\DDS(A,B)$. Notice the similarity of the above description with a morphism inside $\ol{\cat{W}}$ of \cref{rem:intermediate_wiring_category}: a discrete dynamical system is a $\ol{\cat{W}}_\Set$-map $(S,S)\to(A,B)$ in that sense. 

We can now define a functor $\DDS\colon\WW{\Set}\to\Cat$, as in \cite[Def.~4.1,4.6,4.11]{Spivak:2015b}.
To an object $X=(\inp{X},\out{X})\in\WW{\Set}$ we assign the category $\DDS(X)\coloneqq\DDS(\Xin,\Xout)$ as defined above.
To a morphism (wiring diagram) $\phi\colon X\to Y$ \cref{eq:wiring_diag_Int}, we define $\DDS(\phi)\colon\DDS(X)\to\DDS(Y)$ to be the functor
which sends the system $(S,\upd{f},\rdt{f})\in\DDS(X)$ to the system $(S,\upd{g},\rdt{g})\in\DDS(Y)$ having the same state set
$S$, and new update and readout functions given by
\begin{gather*}
\upd{g}\colon\Yin\times S\xrightarrow{1\times\Delta}\Yin\times S\times S\xrightarrow{1\times\rdt{f}\times1}
\Yin\times\Xout\times S\xrightarrow{\fin\times1}\Xin\times S\xrightarrow{\upd{f}}S\\
\rdt{g}\colon S\xrightarrow{\rdt{f}}\Xout\xrightarrow{\fout}\Yout
\end{gather*}
also written, via their mappings on elements, as
\begin{equation}\label{eqn:DDS_phi}
\upd{g}(y,s)\coloneqq\upd{f}\Big(\fin\big(y,\rdt{f}(s)\big),s\Big)
\quad\tn{and}\quad
\rdt{g}(s)\coloneqq\fout\big(\rdt{f}(s)\big)
\end{equation}
It can be verified that $\DDS(\phi)$ preserves composition and identities in $\WW{\Set}$.
Moreover, $\DDS$ has a symmetric lax monoidal structure essentially given by cartesian product:
for systems $F=(S,\upd{f},\rdt{f})\in\DDS(X)$ and $G=(T,\upd{g},\rdt{g})\in\DDS(Y)$,
we define the functor $\DDS_{X,Y}\colon\DDS(X)\times\DDS(Y)\to\DDS(X+Y)$ to map $(F,G)$
to 
$$(S\times T, \Xin\times\Yin\times S\times T\To{\sim}\Xin\times S\times \Yin\times T\To{\upd{f}\times\upd{g}}S\times T,
S\times T\To{\rdt{f}\times\rdt{g}}\Xout\times\Yout).$$

\begin{proposition}\label{prop:algebra_DDS}
There exists a $\WW{\Set}$-algebra $\DDS$, for which the labeled boxes inhabitants are discrete dynamical systems.
\end{proposition}

In particular, if we have any interconnection arrangement of discrete dynamical systems like \cref{picture_wiring_diagram},
we can use this algebra structure to explicitly construct the state set, the update and the readout function
of the new discrete dynamical system which the subsystems form.

Notice that the minor abuse of notation $\DDS(X)\coloneqq\DDS(\Xin,\Xout)$ is explained by \cref{rem:intermediate_wiring_category}.
The algebra $\DDS$ is a prototype of many algebras throughout this paper, e.g.\ the abstract systems in \cref{ch:ADS}.
One will see strong similarities in what we describe next, the algebra of continuous dynamical systems
(details on which can be found in \cite[\S 4]{Vagner.Spivak.Lerman:2015a} or
\cite[\S 2.4]{Spivak:2015b}).

\begin{definition}\label{def:CDS}
Let $A,B$ be Euclidean spaces. An \emph{$(A,B)$-continuous dynamical system}
is a Euclidean space $S$, called the \emph{state space}, equipped with smooth functions
\begin{displaymath}
\begin{cases}
\dyn{f}\colon A\times S\to TS \\
\rdt{f}\colon S\to B
\end{cases}
\end{displaymath}
where $TS$ is the tangent bundle of $S$, such that $\dyn{f}$ commutes with the projections to $S$.
That is, $\dyn{f}(a,s)=(s,v)$ for some vector $v$; it is standard notation to write $(\dyn{f})_{a,s}\coloneqq v$. The maps $\dyn{f}$ and $\rdt{f}$
are respectively called the \emph{dynamics} and \emph{readout} functions. The first is an ordinary differential equation (with parameters in $A$),
and the second is an output function for the system.

A \emph{morphism} of $(A,B)$-continuous dynamical systems $(S_1,\dyn{f_1},\rdt{f_1})\to(S_2,\dyn{f_2},\rdt{f_2})$ is a smooth map $h\colon S_1\to S_2$
such that $\rdt{f_1}(s)=\rdt{f_2}(hs)$ and $Th(\dyn{f_1}(a,s))=\dyn{f_2}(a,hs)$, where $Th\colon TS_1\to TS_2$ is the derivative of $h$.
\end{definition}

The category of $(A,B)$-continuous dynamical systems is denoted by $\CDS(A,B)$.
We can now define a functor $\CDS\colon\WW{\Euc}\to\Cat$ as follows.
To an object $X=(\inp{X},\out{X})$ we assign the category $\CDS(X)\coloneqq\CDS(\Xin,\Xout)$.
To a morphism $\phi\colon X\to Y$, we define $\CDS(\phi)\colon\CDS(X)\to\CDS(Y)$ to be the functor which sends a continuous system 
$(S,\dyn{f},\rdt{f})$
to the continuous system $(S,\dyn{g},\rdt{g})$ where
\begin{equation}\label{CDS(phi)}
\dyn{g}(y,s)\coloneqq\dyn{f}\Big(\fin\big(y,\rdt{f}(s)\big),s\Big)
\quad\tn{and}\quad
\rdt{g}(s)\coloneqq\fout\big(\rdt{f}(s)\big)
\end{equation}
The functor's symmetric lax monoidal structure is again given by cartesian product.

\begin{proposition}
There exists a $\WW{\Euc}$-algebra $\CDS$, for which the labeled boxes inhabitants are continuous dynamical systems. 
\end{proposition}

Even though discrete and continuous dynamical systems have quite different notions of time and continuity, one can compare
\cref{eqn:DDS_phi} and \cref{CDS(phi)} to see that
their algebraic structure (e.g.\ the action on wiring diagrams) is very similar. One aspect of this similarity is summarized as follows.
\begin{proposition}\cite[Thm 4.26]{Spivak:2015b}\label{DDS_to_CDS}
For each $\epsilon>0$, we have a wiring diagram algebra map
\begin{displaymath}
\algmap{\Euc}{U}{\Set}{\CDS}{\DDS}{\alpha_\epsilon}
\end{displaymath}
where $\WW{U}$ is as in \cref{eq:W_U} and $\alpha_\epsilon$ is given by Euler's method of linear
approximation.
\end{proposition}

Thus any continuous $(A,B)$-dynamical system $(S,\dyn{f},\rdt{f})$ gives rise to a discrete $(UA,UB)$ dynamical system
$(US,\upd{f}_\epsilon,\rdt{f}_\epsilon)$ with readout function $\rdt{f_\epsilon}\coloneqq U\rdt{f}$ and with update function
given by the linear combination of vectors $\upd{f}_\epsilon(a,s)\coloneqq s+\epsilon\cdot(\dyn{f})_{a,s}$.

\subsection{Spans as \texorpdfstring{$\WW{}$}{W}-algebras}\label{sec:Span_Systems}

We conclude this section with an explicit construction of a $\WW{\cat{C}}$-algebra for any finitely complete
category $\cat{C}$, which eventually induces the most abstract
notion of a machine in \cref{ch:ADS}. Its special characteristics include the
span-like form of the systems inside the labeled boxes,
as well as naturally induced algebra maps between such systems.

First of all, recall that if $\cat{C}$ is a category with pullbacks, there is a bicategory $\Span_\cat{C}$ of
spans; its objects are the same as $\cat{C}$, and hom-categories $\Span_\cat{C}(X,Y)$ consist of spans $X\tickar Y\coloneqq X\from A\to Y$
in $\cat{C}$ as objects, and commutative diagrams
\begin{displaymath}
\begin{tikzcd}[sep=small]
 & A\ar[dl]\ar[dr]\ar[dd,"\alpha"] & \\
 X && Y \\
 & B\ar[ul]\ar[ur] & 
\end{tikzcd}
\end{displaymath}
as morphisms. Horizontal composition is given by pullbacks, so is associative only up to isomorphism.
If moreover $\cat{C}$ has finite products, therefore is finitely complete, then
spans $X\tickar Y$ can be equivalently viewed as maps $A\to X\times Y$, and so $\Span_\cat{C}(X,Y)\cong\cat{C}/(X\times Y)$.
The following result shows that for such $\cat{C}$, we can always define a symmetric lax monoidal functor from $\WW{C}$ to $\Cat$,
using the bicategory of $\cat{C}$-spans.

\begin{proposition}\label{prop:C_Span_Algebra}
For any finitely complete $\cat{C}$, there exists is a $\WW{\cat{C}}$-algebra 
\begin{equation}\label{eq:defSpnS}
\SpnS_\cat{C}\colon\cat{W}_\cat{C}\longrightarrow\Cat 
\end{equation}
for which labeled boxes inhabitants are spans in $\cat{C}$.
\end{proposition}

\begin{proof}
For each box $X=(\inp{X},\out{X})\in\Ob\WW{\cat{C}}$, we define the mapping on objects $\SpnS_\cat{C}(X)\coloneqq\Span_\cat{C}(\Xin,\Xout)$
to be the hom-category; it consists of $\cat{C}$-spans $S\to\Xin\times\Xout$.
For each wiring diagram $\phi\colon X\to Y$ \cref{eq:wiring_diag}, which by \cref{eq:wiring_diag_Int} determines morphisms in $\cat{C}$
\[
\fin\colon\Yin\times\Xout\to\Xin
\qquad\tn{and}\qquad
\fout\colon\Xout\to\Yout,\]
there is a functor $\SpnS_\cat{C}(\phi)\colon\Span_\cat{C}(\Xin,\Xout)\to\Span_\cat{C}(\Yin,\Yout)$
mapping some $\Xin\from S\to \Xout$ to the outside span below, formed as the composite (pullback)
\begin{equation}\label{eq:composite_span}
 \begin{tikzcd}
 &&& T\ar[drr]\ar[dll] &&& \\
 & \Yin\times\Xout\ar[dl,"\pi_1"']\ar[dr,"\fin"] && S\ar[dr]\ar[dl] && \Xout\ar[dl,equal]\ar[dr,"\fout"] & \\
 \Yin && \Xin && \Xout && \Yout\\
 \Yin\ar[rr,tick] && \Xin\ar[rr,tick] && \Xout\ar[rr,tick] && \Yout
 \end{tikzcd}
\end{equation}
This is clearly functorial, i.e. $\SpnS_\cat{C}(\phi)=(\id,\fout)\circ-\circ(\pi_1,\fin)$ precomposes
any span and span morphism with the right and left $\phi$-induced spans.
Finally, the functor $\SpnS_\cat{C}$ has a symmetric lax monoidal structure by using products symmetry $\sigma$ in $\cat{C}$:
\begin{equation}\label{eq:SpnS_monoidal}
\begin{tikzcd}[row sep=.02in,column sep=.25in]
\Span_\cat{C}(\Xin,\Xout)\times\Span_\cat{C}(\Zin,\Zout)\ar[r,"(\SpnS_\cat{C})_{X,Z}"] &
\Span_\cat{C}(\Xin\times\Zin,\Xout\times\Zout) \\
\left(S\To{p}\Xin\times\Xout,T\To{q}\Zin\times\Zout\right)
 \ar[r,mapsto] &
S\times T\To{\sigma\circ(p\times q)}\Xin\times\Zin\times\Xout\times\Zout
\end{tikzcd}
\end{equation}
\end{proof}

\begin{definition}\label{def:CSpanSystems}
For $\cat{C}$ a finitely complete category, the $\WW{C}$-algebra $\SpnS_\cat{C}$ described in
\cref{prop:C_Span_Algebra} is called the algebra of $\cat{C}$-\emph{span systems}.
\end{definition}

Explicitly, computing the composite \cref{eq:composite_span} for a $\cat{C}$-span
$S\to\Xin\times\Xout$ produces a span $T\to\Yin\times\Yout$
formed by taking the pullback along $\fin$ and composing with $\fout$:
\begin{equation}\label{eq:SpnS_on_arrows}
\begin{tikzcd}[row sep=.4in,column sep=.5in]
T\ar[r]\ar[d]\ullimit\ar[dd,bend right=5,out=300,in=240,dashed] &
S \ar[d] \\
\Yin\times\Xout\ar[r,"{(\fin,\pi_2)}"']\ar[d,"1\times\fout"] &
\Xin\times\Xout \\
\Yin\times\Yout &
\end{tikzcd}
\end{equation}
In particular, if $(f,g)\colon S\to\Xin\times\Xout$ is the span on the right and $k\colon T\to S$ is the top morphism, the left one is $(h,gk)$ inducing the composite $(h,gk\out{\phi})\colon T\to\Yin\times\Yout$.

In fact, since any map $f\colon A\to B$ in a category $\cat{C}$ with pullbacks induces an adjunction
$(f_!\dashv f^*)$ between the slice categories,
the functor $\SpnS_\cat{C}(\phi)\colon\Span_\cat{C}(\Xin,\Xout)\to\Span_\cat{C}(\Yin,\Yout)$
of the above proof can be equivalently expressed as
\begin{equation}\label{eq:slice_functors}
\cat{C}/(\Xin\times\Xout)\xrightarrow{(\fin,\pi_2)^*}\cat{C}/(\Yin\times\Xout)
\xrightarrow{(1\times\fout)_!}\cat{C}/(\Yin\times\Yout)
\end{equation}
between the slice categories.

\begin{remark}
It is the case that for any finitely complete category $\cat{C}$, the bicategory
$\Span_\cat{C}$ has the structure of a \emph{compact closed bicategory}, see \cite{Stay:2016a}.
Since on the 1-category level, every \emph{traced} monoidal category \cite{Joyal.Street.Verity:1996a}
gives rise to a compact closed one, there is evidence that the construction of \cref{prop:C_Span_Algebra}
could work if we replaced $\Span_\cat{C}$ with any traced symmetric monoidal \emph{bicategory}
in an appropriate sense.

Following an anonymous reviewer suggestion, for any such bicategory $\cat{K}$ 
where every object is equipped with a comonoid structure, we could define
a $\WW{\cat{K}}$-algebra as follows:
\begin{displaymath}
\begin{tikzcd}[row sep=.06in]
\cat{W}_K\ar[rr] && \Cat \\
(\inp{X},\out{X})\ar[rr,mapsto]\ar[dd,"\phi"'] && \cat{K}(\Xin,\Xout)\ar[dd,"\Phi"] \\
&&& \\
(\inp{Y},\out{Y})\ar[rr,mapsto] && \cat{K}(\Yin,\Yout)
\end{tikzcd}
\end{displaymath}
where $\Phi$ maps some $f\colon\Xin\to\Xout$ to $\Phi(f)=\mathrm{Tr}\left((\fout\times1)\circ\delta\circ f\circ\fin\right)$.
Explicitly, the map on which the trace acts is
\begin{displaymath}
 \Yin\times\Xout\xrightarrow{\fin}\Xin\xrightarrow{f}\Xout\xrightarrow{\delta}\Xout\times\Xout\xrightarrow{\fout\times1}\Yout\times\Xout.
\end{displaymath}
 
Formally defining the structure of a traced symmetric monoidal
bicategory, albeit possibly further clarifying the origin of such constructions,
is certainly beyond the scope of the current work, and not relevant to the main line.
Indeed, central cases of interest are what we will call deterministic and/or total machines, which do not form traced monoidal categories or
bicategories. Instead, they will form $\WW{}$-algebras enriched in $\Cat$. For relationships between these formalisms, the reader may
see \cite{Spivak.Schultz.Rupel:2016a}, which gives a tight relationship between enriched traced monoidal categories algebras on operads similar to $\WW{}$. 
\end{remark}

Finally, the proposition below shows that the mapping $\cat{C}\mapsto\SpnS_\cat{C}$ extends
to a functor with target category $\WDalg$, as described in \cref{sec:operad_of_wiring_diagrams}.

\begin{proposition}\label{prop:Mch_FCCat}
There exists a functor $\SpnS_{(\Inp)}\colon\FCCat\to\WDalg$ making the following diagram commute
\begin{equation}\label{eq:spantriangle}
  \begin{tikzcd}[row sep=.5in,column sep=.4in]
    &\WDalg\ar[d,"\mathrm{dom}"]\\
    \FCCat\ar[ur,"\SpnS_{(\Inp)}"]\ar[r,"{\WW{(\Inp)}}"',"\cref{eq:W_}"]&\SMonCats
  \end{tikzcd}
\end{equation}
\end{proposition}

\begin{proof}
By \cref{prop:C_Span_Algebra}, $\SpnS_\cat{C}\colon\WW{\cat{C}}\to\Cat$ is a symmetric lax monoidal functor.
Now given a functor $F\colon\cat{C}\to\cat{D}$ between typing categories which preserves finite limits, we
define a monoidal natural transformation $\SpnS_F$ as in the diagram below, where
$\WW{F}\colon\WW{\cat{C}}\to\WW{\cat{D}}$ is a (strong) monoidal functor as in \cref{eq:W_}: 
\begin{equation}\label{eq:defSpnSF}
\algmap{\cat{C}}{F}{\cat{D}}{\SpnS_\cat{C}}{\SpnS_\cat{D}}{\SpnS_F}
\end{equation}
Its components $(\SpnS_F)_X\colon\SpnS_\cat{C}(X)\to\SpnS_\cat{D}(\WW{F}(X))$ for each box $X\in\WW{\cat{C}}$
are functors that assign to any span $p\colon S\to\Xin\times\Xout$ the span $(\SpnS_F)_X(p)\coloneqq Fp$, where $Fp\colon FS
\to F(\Xin\times\Xout)\cong F\Xin\times F\Xout$.
Monoidality follows by $F(p\times q)\cong Fp\times Fq$
and (pseudo)naturality follows from the fact that $F$ preserves pullbacks: applying
$F$ on the feedback construction \cref{eq:SpnS_on_arrows} is the same as performing
that construction on $Fp$.
\end{proof}

Notice that by construction, if $F$ is a faithful functor, the algebra morphism $\SpnS_F$ is an embedding, namely has as components functors which are injective on objects and faithful.

This functor $\SpnS_{(\Inp)}$ gives rise to wiring diagram algebras of central importance for this work,
namely machines (\cref{ch:ADS}) as well as algebra maps between them (\cref{ch:maps_to_machines}). Moreover, certain subfunctors of it
end up capturing fundamental characteristics of systems, such as totality and determinism (\cref{thm:tot_det_variations}).
What comes next in order to reach these formalisms is the description of the appropriate typing categories,
which encompass notions of time.

\section{Interval sheaves}\label{ch:Int}

Guided by some seminal ideas by Lawvere \cite{Lawvere:1986a} and Johnstone
\cite{Johnstone:1999a}, in this chapter we describe a site $\Int$ whose objects can be considered as closed intervals of nonnegative length,
and whose morphisms are inclusions of subintervals. 
A family of subintervals covers an interval if they are jointly
surjective on points. We give several examples of sheaves on $\Int$, called \emph{interval sheaves}.
We then elaborate on variations of this site, giving rise to \emph{discrete} and \emph{synchronous} interval sheaves;
these ultimately correspond to different notions of time for our system models at \cref{ch:ADS}.

In \cref{appendix:DSF} we discuss an equivalence
between the topos of $\Int$-sheaves and the category of \emph{discrete \Conduche fibrations} over the additive monoid of
nonnegative real numbers, due to Peter Johnstone \cite{Johnstone:1999a}.

\subsection{The interval site \texorpdfstring{$\Int$}{Int}}\label{sec:interval_site_Int}

Let $\RRnn$ denote the linearly ordered poset of nonnegative real numbers,
and for any $p\in\RRnn$, let $\Trns{p}\colon\RRnn\to\RRnn$ denote the translation-by-$p$ function $\Trns{p}(\ti)\coloneqq p+\ti$. 

\begin{definition}\label{def:affine_intervals}
The category of \emph{continuous intervals}, denoted $\Int$, is defined to have:
\begin{compactitem}
\item objects $\Ob\Int=\{\ti\in \RRnn\}$,
\item morphisms $\Int(\ti',\ti)=\{\Trns{p}\mid p\in\RRnn\text{ and }p\leq \ti-\ti'\}$,
\item composition $\Trns{p}\circ\Trns{p'}\coloneqq\Trns{p+p'}$, and
\item $\id_\ti=\Tr_0$.
\end{compactitem}
We will sometimes denote $\Tr_p$ simply by $p\colon \ti'\to\ti$.
\end{definition}

The category $\Int$ above is the skeleton of (and hence equivalent to) the category
whose objects are closed, positively-oriented intervals $[a,b]\in\RR$,
and whose morphisms are orientation-preserving isometries, i.e.\ translations by which
one interval becomes a subinterval of another. Under the inclusion of the skeleton,
an object $\ti\in\Int$ is sent to the interval $[0,\ti]\ss\RR$, and we can consider
a map $\Trns{p}\colon\ti'\to\ti$ as a translation such that $[0,\ti']\cong[p,p+\ti']\ss[0,\ti]$:
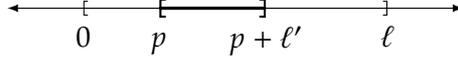
\begin{displaymath}
\begin{tikzpicture}
\draw[latex-latex] (-1,0) -- (5,0) ;
\draw[shift={(0,0)},color=black] node[below=3pt] {$0$};
\draw[shift={(1,0)},color=black] node[below=3pt] {$\fakeell p\fakeell$};
\draw[shift={(2.4,0)},color=black] node[below=3pt] {$p+\ti'$};
\draw[shift={(4,0)},color=black] node[below=3pt] {$\ti$};
\draw[{[-]},thick] (1,0) -- (2.4,0);
\draw[very thick] (1,0) -- (2.4,0);
\draw[{[-]}] (0,0) -- (4,0);
\end{tikzpicture}
\end{displaymath}

The following analogous definition is obtained by
replacing $\RRnn$ by the linear order $\NN$ of natural numbers.

\begin{definition}\label{def:integral_intervals}
The category of \emph{discrete intervals} $\Int_N$ has as objects the set of natural numbers,
$\Ob\Int_N=\{n\in\NN\}$, as morphisms $n'\to n$ the set of translations
$\Trns{p}$, where $p\in\NN$ and $p\leq n-n'$. The composition and identities are given by $+$ and $0$, similarly to \cref{def:affine_intervals}.
\end{definition}

In what follows, fix $\cat{R}$ (resp.\ $\cat{N}$) to be the additive monoid of
nonnegative real (resp.\ natural) numbers, viewed as a category with one object $*$.
\cref{prop:RRnn} gives a natural context for the interval categories $\Int$ and $\Int_N$ from
\cref{def:affine_intervals,def:integral_intervals}: they are the twisted arrow categories of $\cat{R}$ and $\cat{N}$ respectively.

\begin{definition}\cite[\S 2]{Johnstone:1999a}\label{def:twisted_arrow}
For any category $\cat{C}$, the \emph{twisted arrow category} of $\cat{C}$, denoted $\tw{\cat{C}}$, has morphisms of $\cat{C}$
as objects, and for $f\colon x\to y$ and $g\colon w\to z$ in $\cat{C}$, a morphism $f\to g$ in $\tw{\cat{C}}$ is a
pair $(u,v)$ making the diagram commute:
\[
\begin{tikzcd}
	x\ar[d,"f"']&w\ar[l,"u"']\ar[d,"g"]\\
	y\ar[r,"v"']&z
\end{tikzcd}
\]
\end{definition}

\begin{proposition}\label{prop:RRnn}
There is an isomorphism of categories
\begin{displaymath}
\Int\cong\tw{\cat{R}}\quad\mathrm{and}\quad\Int_N\cong\tw{\cat{N}}.
\end{displaymath}
\end{proposition}

\begin{proof}
The objects of $\tw{\cat{R}}$ are nonnegative real numbers $\ti\colon *\to *$,
and a morphism $\ti'\to\ti$ in $\tw{\cat{R}}$
is a pair of nonnegative real numbers $(p,q)$ such that $p+\ti'+q=\ti$. But $q$ is completely determined by $p$, $\ti$, and $\ti'$,
hence we can identify $(p,q)$ with the morphism $\Trns{p}\colon\ti'\to\ti$ of \cref{def:affine_intervals}.
Similarly for the discrete intervals.
\end{proof}

Finally, the one-object categories $\cat{R}$ and $\cat{N}$ themselves fall under the general description
of a factorization-linear category defined below. In particular, $\Fact(q)$ is called the \emph{interval category} in \cite{Lawvere:1986a} and denoted by
$[\![ q ]\!]$ in \cite{Johnstone:1999a}.

\begin{definition}\label{def:factorization_linear}
For any category $\cat{C}$ and morphism $q\colon a\to b$ in $\cat{C}$, define the \emph{$q$-factorization category},
denoted $\Fact(q)$, as follows. Its objects are triples $M=(m,f,g)$, where $f\colon a\to m$ and $g\colon m\to b$
are arrows in $\cat{C}$ composing to $q=g\circ f$.
If $M'=(m',f',g')$, a morphism $M\to M'$ in $\Fact(q)$ is an arrow
$k\colon m\to m'$ in $\cat{C}$ making the following triangles commute
\[
\begin{tikzcd}[row sep=.6in,column sep=.6in]
	a\ar[r,"f'"]\ar[d,"f"']\ar[dr,pos=.7,"q"]&m'\ar[d,"g'"]\\
	m\ar[r,"g"']\ar[ru,dashed,crossing over,pos=.7,"k"] & b
\end{tikzcd}
\]

We furthermore say that $\cat{C}$ is \emph{factorization-linear} if for every morphism $q$,
$\Fact(q)$ is a linear preorder, i.e.\ if for any pair of objects $M,M'\in\Fact(q)$,
\begin{compactenum}
\item there is at most one morphism $M\to M'$, and 
\item there exists either a morphism $M\to M'$ or a morphism $M'\to M$.
\end{compactenum}
\end{definition}

\begin{example}\label{ex:factorization_linear}\hfill

\begin{compactenum}
\item For any directed graph $G$, the free category $\cat{C}=\Fr(G)$ on $G$ is factorization-linear.
Since $\cat{N}$ is the free category on the terminal graph, it is factorization-linear.
\item The category $\cat{R}$ is also factorization-linear. Indeed, if $f,g,f',g'$ are nonnegative real numbers with $f+g=f'+g'$, then
there is exactly one real number $k$ such that $k+f=f'$ (iff $g'+k=g$), and $k\geq 0$ iff there is a map $(f,g)\to (f',g')$.
\end{compactenum}
\end{example}

We will later employ these characterizations in order to consider sheaves
on the categories of intervals. Below we fix some notation for presheaves on them.

\begin{notation}\label{presheaf_notation}
For any small category $\cat{C}$ (such as $\Int$ or $\Int_N$), we denote by $\Psh{\cat{C}}=[\cat{C}\op,\Set]$
the category of presheaves.
Consider an $\Int$-presheaf $A\colon\Int\op\to\Set$. For any continuous interval $\ti$, we refer to
elements $x\in A(\ti)$ as \emph{sections of $A$ on $\ti$}. We refer to sections of length $\ti=0$ as \emph{germs of $A$}.
For any section $x\in A(\ti)$ and any map
$\Trns{p}\colon\ti'\to\ti$, we write $\ShRst{x}{p}{p+\ti'}$ to denote the restriction
$A(\Trns{p})(x)\in A(\ti')$.
Similar notation applies to $\Int_N$.

We can graphically depict a section $x$ on a continuous interval $\ti\in\Int$, together with a restriction, and
a section $y$ on a discrete interval $n\in\Int_N$, as
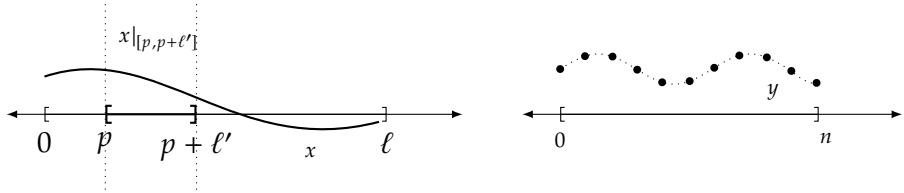
\begin{displaymath}
\begin{tikzpicture}[baseline=(bl)]
\node at (0,0) (bl) {};
\draw[latex-latex] (-1,0) -- (5,0);
\draw[shift={(-.5,0)},color=black] node[below=3pt] {$0$};
\draw[shift={(4,0)},color=black] node[below=3pt] {$\ti$};
\draw[shift={(.3,0)},color=black] node[below=3pt] {$p$};
\draw[shift={(1.5,0)},color=black] node[below=3pt] {$p+\ti'$};
\draw[{[-]}] (-.5,0) -- (4,0);
\draw[thick] (-.5,0.5) .. controls (1,1) and (2,-0.6) ..(3.9,-.1);
\draw[{[-]},thick] (0.3,0) -- (1.5,0);
\draw[dotted] (.3,-1) -- (.3,1.5);
\draw[dotted] (1.5,-1) -- (1.5,1.5);
\node[black] at (1,1) {$\scriptstyle\ShRst{x}{p}{p+\ti'}$};
\node[black] at (3,-.5) {$\scriptstyle x$};
\end{tikzpicture}\qquad
\begin{tikzpicture}[baseline=(bl)]
\node at (0,0) (bl) {};
\draw[latex-latex] (-.5,0) -- (4.5,0);
\draw[shift={(0,0)},color=black] node[below=3pt] {$\scriptstyle 0$};
\draw[shift={(3.5,0)},color=black] node[below=3pt] {$\scriptstyle n$};
\draw[{[-]}] (0,0) -- (3.4,0);
\draw[line width=3pt, line cap=round, dash pattern=on 0pt off 3.5\pgflinewidth] (0,.6) sin (.5,.8)
cos (1,.6) sin (1.5,.4) cos (2,.6) sin (2.5,.8) cos (3,.6) sin (3.5,.4);
\draw [dotted] (0,.6) sin (.5,.8) cos (1,.6) sin (1.5,.4) cos (2,.6) sin (2.5,.8) cos (3,.6) sin (3.4,.4);
\node[black] at (2.8,.3) {$\scriptstyle y$};
\end{tikzpicture}
\end{displaymath}

Certain classes of restrictions will come up often, so we create some special notation for them. 
For any $0\leq \ti'\leq\ti$, we write $\lambda_{\ti'}\colon A(\ti)\to A(\ti')$ and $\rho_{\ti'}\colon A(\ti)\to A(\ti')$ to denote the restrictions along the maps
$[0,\ti']\ss [0,\ti]$ and $[\ti-\ti',\ti]\ss[0,\ti]$ respectively. We refer to them as the \emph{left}
and \emph{right} restrictions of length $\ti'$. In particular if $\ti'=0$, we will often use the notation
\[\sdom{x}\coloneqq\ShRst{x}{0}{0}\in A(0)
\qquad\tn{ and }\qquad
\scod{x}\coloneqq\ShRst{x}{\ti}{\ti}\in A(0)\]
to denote the left and right endpoint germs of $x$.
\end{notation}

\subsection{\texorpdfstring{$\Int$}{Int}-sheaves}\label{sec:Int_sheaves}

We now wish to equip each of the categories $\Int$ and $\Int_N$ with a coverage. In \cite{Johnstone:1999a},
Peter Johnstone defined such a coverage for all twisted arrow categories $\tw{\cat{C}}$ with $\cat{C}$ factorization-linear.
For such $\cat{C}$, he moreover proved that the category of
sheaves on $\tw{\cat{C}}$ is equivalent to the category of discrete \Conduche fibrations over $\cat{C}$, a
brief discussion on which we postpone until \cref{appendix:DSF}; Johnstone's result is \cref{thm:DCFR_Int}.

\begin{definition}\label{def:Johnstone_coverage}
Let $\cat{C}$ be a factorization-linear category, \cref{def:factorization_linear}, and $\tw{\cat{C}}$ its twisted arrow category,
\cref{def:twisted_arrow}. The \emph{Johnstone coverage} on $\tw{\cat{C}}$ is defined as follows:
a covering family for an object $(h\colon x\to z)\in\tw{\cat{C}}$ is any pair of maps $(\id_x,g)\colon f\to h$
and $(f,\id_z)\colon g\to h$ in $\tw{\cat{C}}$, where $g\circ f=h$ in $\cat{C}$,
\[
\begin{tikzcd}
	x\ar[d,"f"']&x\ar[l,"\id_x"']\ar[d,"h"]&&y\ar[d,"g"']&x\ar[l,"f"']\ar[d,"h"]\\
	y\ar[r,"g"']&z&&z\ar[r,"\id_z"']&z
\end{tikzcd}
\]
It can be checked that
this is indeed a coverage \cite[\S 3.6]{Johnstone:1999a},
and we refer to the associated site as the \emph{Johnstone site} for $\tw{\cat{C}}$.%
\footnote{In fact, \cite[3.6]{Johnstone:1999a} explains that $\cat{C}$ being factorization
linear is stronger than necessary to define a coverage on $\tw{\cat{C}}$. However,
this stronger condition is necessary for \cref{cor:DCF_R=Int} to hold.}
\end{definition}

Applying the above definition to $\cat{R}$ and $\cat{N}$, which are
factorization-linear by \cref{ex:factorization_linear}, we obtain the Johnstone
sites $\Int$ and $\Int_N$ on their twisted-arrow categories, see \cref{prop:RRnn}.
Below we explicitly describe them, and fix our terminology.
  

\begin{definition}\label{def:Int_sheaves}
For any interval $\ti\in\Int$, and any real number $p\in[0,\ti]$, we say that the pair of subintervals $[0,p]$ and
$[p,\ti]$, or more precisely the morphisms $\Trns{0}\colon p\to\ti$ and
$\Trns{p}\colon(\ti-p)\to\ti$, form a cover, which we
call the \emph{$p$-covering family} for $\ti$. 
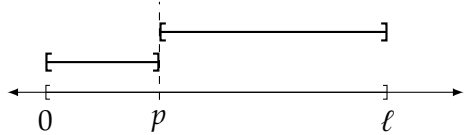
\[
\begin{tikzpicture}
\draw[latex-latex] (0,0) -- (6,0);
\draw[shift={(.5,0)},color=black] node[below=3pt] {$0$};
\draw[shift={(5,0)},color=black] node[below=3pt] {$\ti$};
\draw[shift={(2,0)},color=black] node[below=3pt] {$p$};
\draw[dashed] (2,-.1) -- (2,1.2);
\draw[{[-]}] (.5,0) -- (5,0);
\draw[{[-]},thick] (.5,.4) -- (2,.4);
\draw[{[-]},thick] (2,.8) -- (5,.8);
\end{tikzpicture}
\]
Similarly, the pairs $([0,p],[p,n])$ for any natural number $p\in\{0,\ldots,n\}$ form
the \emph{$p$-covering family} for $n\in\Int_N$. 

If $X$ is an $\Int$-presheaf,
we say that sections $x_1\in X(\ti_1)$ and $x_2\in X(\ti_2)$ are \emph{compatible} if the right endpoint of $x_1$ matches the left endpoint of $x_2$, i.e.\ if $\scod{x_1}=\sdom{x_2}$. Thus $X$ satisfies the sheaf axiom for the above coverage if, whenever $x_1$ and $x_2$ are compatible, there is a unique section $x_1*x_2\in X(\ti_1+\ti_2)$, called the \emph{gluing} of $x_1$ and $x_2$, such that $\lambda_{\ti_1}(x_1*x_2)=x_1$ and $\rho_{\ti_2}(x_1*x_2)=x_2$. 

We refer to an $\Int$-presheaf satisfying the sheaf axiom as a \emph{continuous interval sheaf}, or simply \emph{$\Int$-sheaf}. We denote
by $\Shv{\Int}$ the full subcategory of sheaves; it has the structure of a topos. Similarly, we define
\emph{discrete interval sheaves} or \emph{$\Int_N$-sheaves} and the topos $\Shv{\Int}_N$ thereof.

Each of the inclusions $U\colon\Shv{\Int}\hookrightarrow\Psh{\Int}$ and $\Shv{\Int}_N\hookrightarrow\Psh{\Int_N}$ has a left adjoint. In each case we denote it
$\assShf$ and call it the \emph{associated sheaf} or \emph{sheafification} functor. 
\end{definition}

The category of $\Int$-sheaves is certainly more complex than that of $\Int_N$-sheaves; indeed the latter is just the category of graphs.

\begin{proposition}\label{prop:ShIntN_grphs}
There is an equivalence of categories $\Shv{\Int}_N\simeq\Grph$,
between the category of discrete interval sheaves and the category of graphs.
\end{proposition}

\begin{proof}
To every sheaf $X\in\Shv{\Int}_N$, we associate the graph $G=(V,E,\src,\tgt)$
with $V\coloneqq X(0)$, $E\coloneqq X(1)$, $\src\coloneqq X(\Trns{0})$, and
$\tgt\coloneqq X(\Trns{1})$, where $\Trns{0},\Trns{1}\colon 0\to 1$ are the
two inclusions; see \cref{def:integral_intervals}. To every graph $G$ we associate
the $\Int_N$-sheaf $\Path(G)$ for which $\Path(G)(n)$ is the set of length-$n$ paths
in $G$ with the obvious restriction maps. This is
indeed a sheaf because two paths match if the ending vertex of one is the starting
vertex of the other, in which case the paths can be concatenated.

The above constructions are clearly functorial; we need to check that these functors
are mutually inverse. The roundtrip functor for a graph $G$ returns a graph with the
length-0 and length-1 paths in $G$, which is clearly isomorphic to $G$. For an
$\Int_N$-sheaf $X$, the roundtrip is again $X$ because the coverage on $\Int_N$
ensures that every section is completely determined by the length-0 and length-1 data.
\end{proof}

\begin{remark}
  \Cref{prop:ShIntN_grphs} is closely related to the classical nerve theorem for categories, as
  presented by \cite{Berger:2002a} and generalized by
  \cite{Leinster:2004b,Weber:2007a,Berger.Mellies.Weber:2012a}; see \cite[\S2.01]{Kock:2011a} for a
  clear and concise summary of this approach to nerve theorems.

  Our category $\Int_N$ is equivalent to the subcategory---written $\Delta_0$ in
  \cite{Kock:2011a}---of the simplicial category $\Delta$ containing all objects and all `free'
  morphisms (called \emph{immersions} in \cite{Berger:2002a}). A morphism $\phi\colon [m]\to[n]$ is
  free if $\phi(i+1) = \phi(i)+1$, i.e.\ if $\phi$ is distance preserving. In those sources, it is
  shown that the classical Segal condition characterizing which simplicial sets arise as the nerve
  of some category can be nicely phrased in terms of $\Delta_0$: a presheaf $f\colon\Delta\to\Set$
  is the nerve of a category if and only if the restriction of $f$ to $\Delta_0$ is a sheaf for a
  certain coverage, and moreover the category of sheaves for this coverage on $\Delta_0$ is
  equivalent to $\Grph$. This site $\Delta_0$ is equivalent to the Johnstone site $\Int_N$. We thank
  one of our reviewers for pointing out this connection.
\end{remark}

\begin{example}\label{ex:Int_N}
We do not include numerous examples of $\Int_N$-sheaves, assuming the reader
is familiar with graphs; the following construction will be used later. 

For any set $S$, consider the complete graph
$K(S)\coloneqq(S\times S\tto S)\in\Grph$, having a vertex for each $s\in S$
and an edge from $s\to s'$ for each $(s,s')\in S\times S$.
By \cref{prop:ShIntN_grphs}, there is a canonical way to view any
graph $G$ as an $\Int_N$-sheaf, under which $K(S)$ becomes the $\Int_N$-sheaf whose length-$n$ sections are $S$-lists $\langle s_0,\ldots,s_n\rangle$
of length $n+1$, namely $K(S)(n)=S^{n+1}$. This induces a functor $K\colon\Set\to\Shv{\Int}_N$.
\end{example}

\begin{example}\label{Int_sheaves_ex}
We collect some useful examples of $\Int$-sheaves. 

\begin{compactenum}
\item\label[item]{yoneda} For any $\ti\in\Int$, we have the representable presheaf $\yoneda{\ti}=\Hom_{\Int}(-,\ti)$. For each object $\ti'$ we have $\yoneda{\ti}(\ti')\cong\{p\in\RRnn\mid p\leq\ti-\ti'\}$, so in particular $\yoneda{\ti}(\ti')=\emptyset$ if $\ti'\geq\ti$.
For each morphism $\Trns{q}\colon\ti'\to\ti$, we have $\yoneda{\ti}(p)=p+q$. The coverage on $\Int$ is
sub-canonical, meaning that for any $\ell$, the representable presheaf is in fact a sheaf, i.e.\ $\yoneda{\ti}\in\Shv{\Int}$. One can think of $\yoneda{\ti}$ as the sheaf whose length $\ti'$ sections are all `placements of a length $\ti'$ interval inside of a length $\ti$ interval'.

\item\label[item]{sheaf_of_functions} For any set $S$, the sheaf of functions $\Fnc(S)$ assigns to each interval $\ti$ the set of functions
$\Fnc(S)(\ti)=\{f\colon[0,\ti]\to S\}$, with the obvious restriction maps. The functor $\Fnc\colon\Set\to\Shv{\Int}$ is a right adjoint, see \cref{prop:point}.

\item For any set $S$, the constant sheaf $\cnst(S)$ is defined by $\cnst(S)(\ti)\coloneqq S$ for any $\ti\in\Int$,
and identity restrictions maps. The functor $\cnst\colon\Set\to\Shv{\Int}$ has both a right adjoint $\Shv{\Int}(1,-)\colon\Shv{\Int}\to\Set$ and a left adjoint $\pi_0\colon\Shv{\Int}\to\Set$, given by $\pi_0(X)\coloneqq X(0)/\!\!\sim$, where $\sim$ is the equivalence relation generated by $\sdom{x}\sim\scod{x}$ for any $\ti$ and $x\in X(\ti)$.

\item \label[item]{ex:gluable_smooth_functions} 
If $X$ is a $C^n$-manifold, there is a sheaf of $C^n$ curves in $X$; we just need to be careful with the length-$0$ sections. Define a sheaf $C^n(X)$ of $C^n$ curves through $X$ by
\begin{equation}\label{eqn:Cn}
C^n(X)(\ti)\coloneqq \left\{(\epsilon,f)\mid\epsilon>0, f\in C^n\big((-\epsilon,\ti+\epsilon),X\big)\right\}/\sim
\end{equation}
where we set $(\epsilon,f)\sim(\epsilon',f')$ if one is the restriction of the other, say $f'|_{(-\epsilon,\ti+\epsilon)}=f$.

For $\ti>0$, the closed interval $[0,\ti]$ is a manifold with boundary, so in fact
the isomorphism $C^n(X)(\ti)\cong\{f\colon[0,\ti]\to X\mid f\in C^n\}$ is the usual definition of continuously
$n$-times differentiable functions out of $[0,\ti]$. However, \cref{eqn:Cn} makes sense when $\ti=0$ also: it defines the
set of $n$-jets in $X$. It should be emphasized that the length-$0$ sections of $C^n(X)$ are not the points in $X$,
so the notation $C^n(X)$ could be considered misleading. However, it is the right definition
to define a sheaf of $C^n$ curves, because gluing two $C^n$ curves whose endpoints have the same $n$-jet results in a $C^n$ curve. 

Also notice that if $m\leq n$, there is a sheaf morphism $C^n(X)\to C^m(X)$.


\item\label[item]{ex:synchronizing} Let $0,1\colon\yoneda{0}\to\yoneda{1}$ denote the image under $\yoneda{}\colon\Int\to\Shv{\Int}$
of the left and right endpoint inclusions $\Tr_0,\Tr_1\colon 0\to 1$ in $\Int$. Then we define the \emph{periodic synchronizing
sheaf} to be the quotient $\Sync\coloneqq\yoneda{1}/(0=1)$. For an interval $\ti$, we may identify
\begin{equation}\label{sync(l)}
\Sync(\ti)\cong\RR/\ZZ\cong\{0\leq[\theta]<1\} 
\end{equation}
where $[\theta]$ denotes the equivalence class modulo $\ZZ$ of $\theta\in\RR$ and may be called the \emph{phase}. One can imagine a length-$\ti$ section of $\Sync$ as a helix of height $\ti$; two different sections differ only by their phase (turning the helix by some angle $\theta$).
Given a subinterval $\Tr_p\colon\ti'\to\ti$, we have $\Sync(\Tr_p)(\theta)=[\theta+p]$. That is, we simply restrict the helix to the subinterval $[p,p+\ti']\ss [0,\ti]$. See \cref{rem:whysynchrony} for how this sheaf arises in our framework.
\end{compactenum}
\end{example}

Some of the above constructions are naturally connected in the following way.

\begin{proposition}\label{prop:point}
Consider the functor $\Stk\colon\Shv{\Int}\to\Set$, given by $\Stk(A)\coloneqq A(0)$.
It is the inverse image of an essential geometric morphism
\begin{displaymath}
\begin{tikzcd}[sep=.6in]
\Set\ar[r,shift right=2.5,"\Fnc"']\ar[r,shift left=2.5,"(\Inp)\cdot\yoneda{0}"] &
\Shv{\Int}\ar[l,"\Stk" description]
\end{tikzcd}
\end{displaymath}
i.e.\ $\Stk$ is an essential point of the topos $\Shv{\Int}$.
\end{proposition}

\begin{proof}
The right adjoint of $\Stk$ is $\Fnc(\Inp)$ defined in \cref{Int_sheaves_ex}(\ref{sheaf_of_functions}).
The left adjoint of $\Stk$ sends $S$ to the copower $S\cdot\yoneda{0}\in\Shv{\Int}$, where $\yoneda{0}$
is as in \cref{Int_sheaves_ex}(\ref{yoneda}).
\end{proof}

The following notions of \emph{extension}
for a section and a sheaf will be of central importance in \cref{ch:ADS}.

\begin{definition}\label{state_extension}\label{Ext}
Let $A\in\Shv{\Int}$ be a sheaf, and $\epsilon\geq 0$. For a section $a\in A(\ti)$, if
$a'\in A(\ti+\epsilon)$ is a section with $\ShRst{a'}{0}{\ti}=a$, we call $a'$ an \emph{$\epsilon$-extension of $a$}.
Moreover, we define the \emph{$\epsilon$-extension sheaf} $\Ext{\epsilon}(A)$ of $A$ by assigning to $\ti$ the set
\[\Ext{\epsilon}(A)(\ti)\coloneqq A(\ti+\epsilon)\]
induced by the functor $\Int\to\Int$ sending $\ti\mapsto\ti+\epsilon$ and $\Tr_p\mapsto\Tr_p$.
\end{definition}

Essentially, the extension sheaf of $A$ includes, for each interval of time, information (i.e. sections) for $\epsilon$-longer intervals;
this will be essential when discussing feedback of systems, later formalized e.g. by \cref{def:inertial_map}.

The above defines an endofunctor $\Ext{\epsilon}\colon\Shv{\Int}\to\Shv{\Int}$%
\footnote{It should be clear that this functor $\Ext{}$ has nothing to do with the $\Ext{}$ functor from homological algebra.}.
There exist two natural transformations $\lambda,\rho\colon\Ext{\epsilon}\Rightarrow\id_{\Shv{\Int}}$
whose components $\Ext{\epsilon}(A)\to A$ are given by applying either left or right restriction
as in \cref{presheaf_notation},
\begin{equation}\label{eq:Ext_maps}
\begin{tikzcd}
\Ext{\epsilon} A(\ti)=A(\ti+\epsilon)
\ar[r,shift left,"\lambda_\ti"]
\ar[r,shift right,"\rho_\ti"']
& A(\ti).
\end{tikzcd}
\end{equation}
That is, for $x\in A(\ti+\epsilon)$
we have $\lambda(x)\coloneqq\lambda_{\ti}(x)=\ShRst{x}{0}{\ti}$ and $\rho(x)\coloneqq\rho_{\ti}(x)=\ShRst{x}{\epsilon}{\ti+\epsilon}$.
Similarly for a discrete-interval sheaf $A\in\Shv{\Int}_N$ and $n\in\NN$, we define the extension sheaf $\Ext{n}(A)$ and the
natural transformations $\lambda,\rho$ with components $\Ext{n}(A)\to A$.

We finish this section by collecting a few useful lemmas.
\begin{lemma}\label{lemma:ext_pbs}
For any $\epsilon\geq 0$ the functor $\Ext{\epsilon}\colon\Shv{\Int}\to\Shv{\Int}$ commutes with all limits.
\end{lemma}
\begin{proof}
Limits in $\Shv{\Int}$ are taken pointwise.
\end{proof}

\begin{lemma}\label{lemma:ext}
Suppose given a commutative triangle of $\Shv{\Int}$-sheaves
\[
\begin{tikzcd}
	\Ext{\epsilon}S\ar[r,"h"]\ar[rd,"\lambda"']&S'\ar[d,"k"]\\
	&S
\end{tikzcd}
\]
Then $h$ is an epimorphism (monomorphism) if and only if all components of the underlying
presheaf morphism $Uh$ are surjective (injective).
\end{lemma}

\begin{proof}
First of all, the functor $U\colon\Shv{\Int}\subseteq\Psh{\Int}$ is faithful, so if $h$ is any sheaf map
and $Uh$ is an epimorphism (monomorphism) then so is $h$. Moroever, since $U$ is right adjoint to the sheafification $\assShf$, it preserves monomorphisms.

Now suppose $h$ is an epimorphism; we will show that $h_\ti$ is a surjection, for any $\ti\in\Int$, so choose $s'\in S'(\ti)$. Recall
that $h$ being an epimorphism means that there exists a cover  $0=\ti_0\leq\cdots\leq\ti_n=\ti$ such that each restriction
$s'_i\coloneqq\ShRst{s'}{\ti_i}{\ti_{i+1}}=h_{\ti_{i+1}-\ti_i}(s_i)$ for some $s_i\in\Ext{\epsilon}S(\ti_{i+1}-\ti_i)$. We may assume
$n=2$. Since the $s'_i$ are compatible in $S'$ and since $\lambda(s_i)=k(s'_i)$, we have that $\lambda(s_1)$ and
$\lambda(s_2)$ are compatible sections of $S$. Thus we can glue
$s\coloneqq\lambda(s_1)*s_2\in S(\ti_2+\epsilon)$, giving $s\in\Ext{\epsilon}S(\ti)$ with $h_{\ti}(s)=s'$.
\end{proof}

\subsection{Synchronization}\label{sec:synchronization}

For what follows, it is essential that we compare the toposes
$\Shv{\Int}_N$ and $\Shv{\Int}$ as the far ends of the time spectrum: from sections over continuous
intervals of time to those over specific, equally-spaced ticks of the clock. The following remark discusses why the notion of a \emph{synchronous} sheaf is required for a coherent common time-framework.

\begin{remark}\label{rem:whysynchrony}
Suppose $X$ is a discrete-interval sheaf and we want to translate it into a continuous-interval sheaf $X'$. What should the $X'$-sections of length $1.5$ be? Taking them to be $X(1)$ does not give rise to any reasonable sort of restriction map. Instead, the adjoint to the obvious forgetful functor from continuous to discrete sheaves, denoted $\Sigma_i$ in \cref{prop:int_sheaves_presheaves}, naturally adds a notion of \emph{phase}---a number between 0 and 1---so as to enable arbitrary real-number restrictions. In particular, $\Sigma_i(\1)$ of the terminal discrete sheaf is the sychronizing sheaf of \cref{Int_sheaves_ex},(\ref{ex:synchronizing}). Thus while there is a left adjoint functor from discrete- to continuous-time sheaves that preserves pullbacks (because sections of matching phases can be glued), it does not preserve the terminal object: this means it can be replaced by a geometric morphism into the slice topos $\Shv{\Int}/\Sigma_i(\1)$, which can later be used as a type-changing functor \cref{eq:W_} for wiring diagram algebras.
\end{remark}

Let us explicitly establish the above observations. If $i\colon\Int_N\to\Int$ is the evident inclusion, it induces an adjunction between the respective toposes of sheaves as follows.

\begin{proposition}\label{prop:i_pres_covers}
The functor $i\colon\Int_N\to\Int$ preserves covering families.
\end{proposition}
\begin{proof}
This follows directly from \cref{def:Int_sheaves}.
\end{proof}

Let $i^*\colon\Psh{\Int}\to\Psh{\Int_N}$ denote the functor given by
pre-composing a presheaf $X\colon\Int\op\to\Set$ with $i\op$, and let $\assShf$ denote
the sheafification as in \cref{def:Int_sheaves}; the following fact is well-known.

\begin{proposition}\label{prop:int_sheaves_presheaves}
The functor $i\colon\Int_N\to\Int$ induces an adjunction $\Sigma_i\dashv\Delta_i$ between sheaf toposes
\begin{equation}\label{eqn:int_sheaves_presheaves}
\begin{tikzcd}[sep=.6in]
\Shv{\Int}_N\ar[d,shift left=1.7,"U","{\bbot}"']\ar[r,shift left=1.7,"\Sigma_i","{\bot}"'] &
\Shv{\Int}\ar[d,shift left=1.7,"U","{\bbot}"']\ar[l,shift left=1.7,"\Delta_i"] \\
\Psh{\Int_N}\ar[u,shift left=1.7,"\assShf"]\ar[r,shift left=1.7,"\Lan_i","{\bot}"'] &
\Psh{\Int}\ar[u,shift left=1.7,"\assShf"]\ar[l,shift left=1.7,"i^*"]
\end{tikzcd}
\end{equation}
such that $U\circ\Delta_i=i^*\circ U$, hence also $\Sigma_i\circ\assShf\cong\assShf\circ\Lan_i$.
\end{proposition}
\begin{proof}
Because $i$ preserves covers (\cref{prop:i_pres_covers}), the composite $\Shv{\Int}\To{U}\Psh{\Int}\To{i^*}\Psh{\Int_N}$
factors through $\Shv{\Int}_N$, defining $\Delta_i$. Since $U$ is fully faithful and has a left adjoint $\assShf$, it can be seen that
the composite $\Sigma_i\coloneqq\assShf\circ\Lan_i\circ U$ is left adjoint to $\Delta_i$.
\end{proof}

We will see in \cref{cor:pb_and_sync_as_sigma} that
$\Sigma_i\dashv\Delta_i$ is not in general a geometric morphism 
because the left adjoint
$\Sigma_i$ does not preserve the terminal object (nor binary products).
However, what is important for our purposes is to see that it preserves pullbacks, thus
is a pre-geometric morphism in the sense of \cite[{}~A4.1.13]{Johnstone:2002a}.
This fact will be deduced from a direct formula calculation for $\Sigma_i$ as a coproduct (denoted using $\sqcup$). For any $x\in\RR$,
let $\floor{x},\ceil{x}\in\ZZ$ denote the floor (resp.\ ceiling) of $x$, i.e.\ the largest
integer $\floor{x}\leq x$ (resp.\ the smallest integer $\ceil{x}\geq x$).

\begin{proposition}\label{prop:Sigma_formula}
For any discrete-interval sheaf $X\in\Shv{\Int}_N$, there is an isomorphism
\begin{equation}\label{eqn:Sigma_formula}
\Sigma_iX(\ti)\cong\bigsqcup_{r\in[0,1)}{X\big(\ceil{r+\ti}\big)}
\end{equation}
Moreover, $\Sigma_i$ commutes with the forgetful functor $U$ in \cref{eqn:int_sheaves_presheaves},
i.e.\ $U\circ\Sigma_i\cong\Lan_i\circ U$.
\end{proposition}

\begin{proof}
The pointwise left Kan extension $C\coloneqq\Sigma_iX=\Lan_iX$ is computed as a colimit
\begin{equation}\label{eqn:colimit_Sigma}
C(\ti)=\colim_{\substack{n\in\Int_N\\(i(n)\to\ti)\in\Int\op}}X(n)=\int^{n\in\Int_N}\Int\op\big(i(n),\ti\big)\times X(n).
\end{equation}
The indexing category $(i\op\downarrow\ti)$ is a poset with the property that every element maps to a unique maximal element. Such maximal elements correspond to maps $r\colon\ti\to n$ in $\Int$,
for which $0\leq r<1$ and $0\leq n-(r+\ti)<1$, i.e.\ $n=\ceil{r+\ti}$.
Thus for each $\ti$, the set $C(\ti)$ is exactly as in the formula \cref{eqn:Sigma_formula}: an element of $C(\ti)$ is
a pair $(r,x)$, where $0\leq r<1$ and $x\in X(n)$.  We need to show that $C$ is a sheaf, but to do so we must
better understand the restriction maps; see \cref{ex.restriction} for intuition.

So suppose given $(r,x)$ as above. Given a map $p\colon\ti'\to\ti$, let $a\coloneqq \floor{p+r}$, $b\coloneqq\ceil{p+r+\ti'}$,
$n'\coloneqq b-a$, and $r'\coloneqq p+r-a$, so $0\leq r'<1$. We have $a\colon n'\to n$ in $\Int_N$ such that $r\circ p=a\circ r'$.
Given a section
$x\in X\big(\ceil{r+\ti}\big)$, the restriction map for $C$ is given by
\begin{equation}\label{eqn:rest_formula}
C(p)(r,x)=\left(r',\ShRst{x}{a}{b}\right).
\end{equation}

Suppose given sections $(r,x)\in C(\ti)$ and $(r',x')\in C(\ti')$, where $x\in X(n)$ and $x'\in X(n')$.
If they are compatible then, by the above restriction formula \cref{eqn:rest_formula}, we must have $r'=\ti+r-\floor{\ti+r}$ and
\[
\ShRst[\big]{x}{\floor{\ti+r}}{\ceil{\ti+r}}=\ShRst{(r,x)}{\ti}{\ti}=\ShRst{(r',x')}{0}{0}=\ShRst[\big]{x'}{\floor{r'}}{\ceil{r'}}.
\]
If we denote by $n_0$ the length of this section, then $n_0$ is either $0$ or $1$, depending on whether
$\ti+r=\floor{\ti+r}$ or not. Either way, $x$ and $x'$ are compatible sections of $X$ and can be glued to form $x*x'\in X(n-n_0+n')$.
Thus we have shown that $C=\Lan_iX$ is a sheaf. It follows that $U\circ\assShf\circ\Lan_i=\Lan_i$, so $U\circ\Sigma_i=\Lan_i\circ U.$
\end{proof}

\begin{example}\label{ex.restriction}
Given a discrete sheaf $X$, a length $\ti$ section of the continuous sheaf $C\coloneqq\Sigma_i(X)$ is given by choosing a phase $r\in[0,1)$ and an element $x\in X(n)$, where $n\coloneqq\ceil{r+\ti}$. Thus $x$ is a section of $X$ on the smallest discrete interval containing the continuous interval, when embedded to start at $r$. In black is a picture of a map $r\colon\ti\to n$ in $\Int$, where $n=6$, $r=\frac{2}{3}$, $\ell=\frac{29}{6}$, and $r+\ell=\frac{11}{2}$:
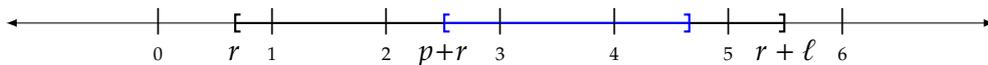
\begin{displaymath}
\begin{tikzpicture}[decoration={brace, amplitude=5pt}, text depth=1ex, text height=2ex]
  \draw[latex-latex] (-2,0) -- (11,0) ;
  \draw (1,-.4) node {$r$} (8.25,-.4) node {$r+\ti$};
  \draw[{[-]}, thick] (1,0) -- (8.25,0);
  \draw[blue, {[-]}, thick] (3.75,0) -- (7,0);
  \node at (3.75,-.4) {$p{+}r$};
  \foreach \x in {0,...,6} {\draw (1.5*\x,0) node {|} (1.5*\x,-.4) node {\tiny\x};}
\end{tikzpicture}
\end{displaymath}
therefore a section in $C\left(\frac{29}{6}\right)$ could be a pair $\left(\frac{2}{3},x\in X(6)\right)$.

Suppose we want to restrict $(\frac{2}{3},x)$ to the blue interval. In this case we visually we see that $p+r=2.5$, so we first take the new phase to be its fractional part, $r'\coloneqq 0.5$. We then restrict $x$ to the smallest discrete interval containing the blue interval when embedded to start at $0.5$, i.e.\ $\ShRst{x}{2}{5}$.
\end{example}

The periodic synchronizing sheaf $\Sync$ from \cref{Int_sheaves_ex}(\ref{ex:synchronizing}) arises as $\Sigma_i$ of the terminal $\Int_N$-sheaf.

\begin{corollary}\label{cor:pb_and_sync_as_sigma}
The functor $\Sigma_i$ preserves pullbacks but not the terminal object; indeed we have
\[\Sigma_i(\cnst\singleton)\cong\Sync.\]
\end{corollary}

\begin{proof}
The fact that $\Sigma_i$ preserves pullbacks can be checked by hand from \cref{eqn:Sigma_formula};
formally, it is due to the fact that the indexing category $(i\op\downarrow\ti)$ in
\cref{eqn:colimit_Sigma} is a coproduct of filtered categories,
see \cite[Ex.~1.3(vi)]{Adamek.Borceux.Lack.Rosicky:2002a}.

The fact that $\Sigma_i(\cnst\singleton)\cong\Sync$ can again be checked by hand, e.g.\ its set of length $\ell$ sections is given by the coproduct $\bigsqcup_{r\in[0,1)}1\cong[0,1)$. Another approach is to note that $\Sigma_i$ can be regarded as a functor $\Sigma_i\colon\Grph\to\Shv{\Int}$ by
\cref{prop:ShIntN_grphs}. Using \cref{eqn:colimit_Sigma} one can check that the image of the single-vertex graph \fbox{$\bullet$} is the
sheaf $\yoneda{0}$ and that the image of the walking edge graph $E=\fbox{$\bullet^0\to\bullet^1$}$
is $\yoneda{1}$, from \cref{Int_sheaves_ex}(\ref{yoneda}). The terminal $\Int_N$-sheaf $\1=\cnst\singleton$ corresponds to the quotient graph $E/(0=1)$.
Since $\Sigma_i$ preserves colimits, we indeed have $\Sigma_i(\cnst\singleton)\cong\Sync$.
\end{proof}

Our next goal is to compare the discrete-interval sheaf topos $\Shv{\Int}_N$ with the slice topos $\Shv{\Int}/\Sync$
over the synchronizing sheaf \cref{sync(l)}, whose objects $X\to\Sync$ we sometimes call \emph{synchronous sheaves}. An object $X\to\Sync$ of this should be visualized as an ordinary $\Int$-sheaf $X$, where every section of $X$ has been assigned a section of $\Sync$, i.e.\ a helix at some phase $\theta\in[0,1)$ or a watch-hand at some position. Two sections $x,y$ are `in sync' if they both have the same phase.

\begin{proposition}\label{prop:synchronized_geom_morph}
The adjunction $\Sigma_i\dashv\Delta_i$ between $\Shv{\Int}$ and $\Shv{\Int}_N$ from \cref{prop:int_sheaves_presheaves}
factors through a geometric morphism
\begin{equation}\label{eq:Sigma_i'}
\begin{tikzcd}
\Shv{\Int}_N\ar[r,shift left=1.7,"\Sigma'_i","{\bot}"']&\Shv{\Int}/\Sync\ar[l, shift left=1.7, "\Delta'_i"]
\end{tikzcd}
\end{equation}
in which the left adjoint $\Sigma_i'$ is fully faithful.
\end{proposition}

\begin{proof}
Define a functor $\Sigma'_i$ which sends $X\in\Shv{\Int}_N$ to the synchronous sheaf $\Sigma_i(X\To{!}\1)$.
Because $\Sigma_i$ preserves pullbacks, it can be verified that
$\Sigma'_i$ preserves all finite limits.

We now construct its right adjoint.
First note that $\Delta_i\colon\Shv{\Int}\to\Shv{\Int}_N$ is given by pre-composing with the inclusion $i\colon\Int_N\to\Int$. Applying it to $\Sync$ results in the constant sheaf $\Delta_i(\Sync)\cong\cnst(\RR/\ZZ)\in\Shv{\Int}_N$; indeed, it is constant because restricting $\Sync$ along maps in $\Int_N$ does not change the phase. 
Hence, given a synchronous sheaf $Y\to\Sync$, define $\Delta'_i(Y)$
to be the following pullback in $\Int_N$:
\[
\begin{tikzcd}[column sep=.5in]
	\Delta'_i(Y)\ar[r]\ar[d]\ar[dr,phantom,description,"\lrcorner" very near start]&\Delta_i(Y)\ar[d]\\
	\cnst\singleton\ar[r]&\Delta_i(\Sync)
\end{tikzcd}
\]
where the bottom map is the unit $\cnst\singleton\to\Delta_i\Sigma_i(\cnst\singleton)$; see \cref{cor:pb_and_sync_as_sigma}.
Since $\cnst\singleton\in\Shv{\Int}_N$ is a terminal object, it follows immediately that $\Delta_i'$
is right adjoint to $\Sigma_i'$.

The adjunction $\Sigma_i\dashv\Delta_i$ can be verified to be
the composite of $\Sigma'_i\dashv\Delta'_i$ with the adjunction
\begin{equation}\label{slice_adjunction}
\begin{tikzcd}
	\Shv{\Int}/\Sync\ar[r,shift left=1.7,"\Sigma_{\Sync!}","{\bot}"']&\Shv{\Int}\ar[l,shift left=1.7,"\Delta_{\Sync!}"]
\end{tikzcd}
\end{equation}
between slice categories as induced by the unique map $\Sync!\colon\Sync\to\singleton$ in $\Shv{\Int}$.
That is, $\Sigma_{\Sync!}(X\to\Sync)=X$ and $\Delta_{\Sync!}(Y)=(Y\times\Sync\to\Sync)$.
Finally, using \cref{prop:Sigma_formula}, one checks that
the unit $\eta\colon X\to\Delta'_i\Sigma'_iX$ is an isomorphism, so
$\Sigma'_i$ is fully faithful.
\end{proof}

Having introduced and related the toposes of interval sheaves in order to appropriately capture behaviors as sections over lengths of time, we now proceed to the formalization of abstract systems as algebras for wiring diagram operads, typed in those sheaves.

\section{Machines as generalizations of dynamical systems}\label{ch:ADS}

In this chapter, we turn to this work's basic goal: the explicit description of general sorts of processes, called \emph{machines},
in terms of interval sheaves. These machines take on many forms, from discrete and continuous dynamical
systems to more general deterministic or total systems,
to even more abstract objects.
What all these examples have in common is that a collection of machines can be interconnected to
form a new machine, and this operation is coherent with respect to nesting of wiring diagrams. In other words,
all of our machines will be algebras on an operad of
wiring diagrams, as defined in \cref{ch:wiring_diagrams_and_algebras}. Subsequently, in \cref{ch:maps_to_machines}
we produce algebra maps that relate all of these various examples.

Each machine will have an interface consisting of input and output ports, where each port carries a sheaf
of `possible signal behaviors', as in
\begin{displaymath}
\begin{tikzpicture}[oriented WD, bbx = 1cm, bby =.5cm, bb min width=1cm, bb port length=4pt, bb port sep=1]
	\node[bb={1}{1}] (X) {$S$};
	\draw[label]
		node [left=2pt of X_in1] {$A$}
		node [right=2pt of X_out1] {$B$}
		;
	\end{tikzpicture}
\end{displaymath}
The possible behaviors of the whole machine---including what occurs on its ports---is also represented by a sheaf $S$. 
We first give a very general definition: a machine is a span $A\from S\to B$. Then 
we restrict our attention to \emph{total} processes, which have the property that for any initial state and compatible input behavior,
there exists a compatible state evolution. 
We also consider \emph{deterministic} processes, which have the property
that the state evolution is unique given any input. Total and deterministic machines are discussed
in \cref{sec:total_det_inert}, and their compositionality is proven in \cref{sec:feedback}.
In the next chapter, we will see that while discrete dynamical systems as defined in \cref{def:DDS} are always both total and deterministic,
this is not the case for continuous dynamical systems of \cref{def:CDS} (basically because solutions to general
ODEs on unbounded domains may not exist and may not be unique; see \cref{rem:ODE_not_td}).

In \cref{safety_contracts} we discuss another related notion, namely that of 
\emph{contracts} on machines---which may or may not be satisfied by a given machine---that express behavioral guarantees in
terms of inputs and outputs. We will show that contracts also form a wiring diagram algebra.

\subsection{Continuous machines}\label{sec:ADS_only}

We begin with the most general notion of process that we will use. It primarily serves
the purpose of being all-inclusive, so that the various systems of interest can be seen as special cases.
The type of information handled is formalized by continuous interval sheaves, \cref{def:Int_sheaves},
whereas the algebraic structure is formalized by our preliminary work on span algebras in \cref{sec:Span_Systems}. 
While such machines ostensibly have a notion of input and output, the definition is in fact symmetric,
so it fits in with the work of \cite{Willems.Polderman:2013a}.

\begin{definition}\label{def:ADS}
Let $A,B\in\Shv{\Int}$ be continuous interval sheaves. A \emph{continuous $(A,B)$-machine} is a span
\begin{displaymath}
\begin{tikzcd}[row sep=15pt]
& S\ar[dl,"\ii{p}"']\ar[dr,"\oo{p}"] & \\
A && B
\end{tikzcd}
\end{displaymath}
in the topos $\Shv{\Int}$; equivalently, it is a sheaf $S$
together with a sheaf map $p\colon S\to A\times B$.
\end{definition}

We refer to $A$ as the \emph{input sheaf}, to $B$ as the \emph{output sheaf}, and to $S$ as the \emph{state sheaf}.
Their sections of arbitrary lengths $\ti$ should be thought of as all possible information inputted, worked out and outputted
by the machine, during a length of time $\ti$.
We call $\ii{p}$ the \emph{input sheaf map} and $\oo{p}$ the \emph{output sheaf map};
they are given by functions $\ii{p}_\ti\colon S(\ti)\to A(\ti)$, $\oo{p}_\ti\colon S(\ti)\to B(\ti)$
on sections of length $\ti$.
We denote by $\ADS(A,B)\coloneqq\Shv{\Int}/(A\times B)$, or sometimes by $\ADS[C](A,B)$,
the topos of all continuous $(A,B)$-machines.

\begin{example}\label{ex:first_ex_machines}
Here we collect some examples of continuous machines.
\begin{compactenum}
\item\label[item]{item:sheaf_maps}
If $f\colon A\to B$ is a morphism of $\Int$-sheaves, there is a continuous $(A,B)$-machine given by $(\id_A,f)\colon A\to A\times B$.
\item If $A$ and $B$ are $\Int$-sheaves, the identity span $A\times B\to A\times B$ corresponds to the machine for which input and output are completely uncoupled.
\item
For any sheaf $A$ and $\epsilon>0$, there is a machine (described in \cref{ex:delay_box}), which acts as an $\epsilon$-delay: the $A$ input at any time is output after a delay of $\epsilon$-seconds. 
\item In \cref{sec:DandCDS} we discussed discrete and continuous dynamical systems.
In \cref{DDS_to_ADS,CDS_to_ADS}, we will give algebra maps realizing them as machines.
\end{compactenum}
\end{example}

Let $\WW{\Shv{\Int}}$ denote the category of $\Shv{\Int}$-labeled boxes and wiring diagrams, as in \cref{sec:operad_of_wiring_diagrams}.
As in the case of discrete and continuous dynamical systems,
the important aspect of these newly-defined machines
is that they can be arbitrarily wired together to form new continuous
machines, i.e. they are algebras for the wiring diagram operad $\WW{\Shv{\Int}}$.
The following proposition is a corollary of \cref{prop:C_Span_Algebra}, for $\cat{C}=\Shv{\Int}$
the finitely complete category of continuous interval sheaves.

\begin{proposition}\label{prop:ADS}
Continuous machines form a $\WW{\Shv{\Int}}$-algebra $\ADS\colon\WW{\Shv{\Int}}\to\Cat$.
\end{proposition}

\begin{proof}
Recall from \cref{eq:dependent_product} that to each object (box) $X\in\WW{\Shv{\Int}}$ we can associate two sheaves, $\Xin,\Xout\in\Shv{\Int}$,
by taking the product of the types of the input and output $\Shv{\Int}$-typed sets.
The algebra is
\begin{equation}\label{eq:Mchphi}
\begin{tikzcd}[row sep=.06in]
\ADS\colon\WW{\Shv{\Int}}\ar[rr] && \Cat \\
(\inp{X},\out{X})\ar[rr,mapsto]\ar[dd,"\phi"'] && \Shv{\Int}/\Xin\times\Xout\ar[dd,"\ADS(\phi)"] \\
&&& \\
(\inp{Y},\out{Y})\ar[rr,mapsto] && \Shv{\Int}/\Yin\times\Yout
\end{tikzcd}
\end{equation}
where $\ADS(\phi)$ maps the span $(\ii{p},\oo{p})\colon S\to\Xin\times\Xout$ to
the span $(\ii{q},\oo{q})\colon T\to\Yin\times\Yout$ formed by first taking the pullback of 
$(\ii{p},\oo{p})$ along $\fin$ and then postcomposing with $\fout$, as showed in \cref{eq:SpnS_on_arrows}.
Since limits in $\Shv{\Int}$ are formed pointwise, we can describe the set of length-$\ti$ sections of the new state sheaf $T$ explicitly:
\begin{equation}\label{eq:newstatesheaf}
T(\ti)\cong\left\{(s,y)\in S(\ell)\times\Yin(\ti)\mid\ii{p}_\ti(s)=\finl{\ti}\big(y,\oo{p}_\ti(s)\big)\right\}.
\end{equation}
Roughly, $T$ represents state evolutions of the machine $(S,\ii{p},\oo{p})$ inhabiting the inside box $X$,
together with $\inp{\phi}$-compatible sections of $Y$-input.
With this representation, the input sheaf map of the formed machine is $\ii{q}_\ti(s,y)=y$,
and the output sheaf map is $\oo{q}=\fout\circ\oo{p}$. The symmetric monoidal structure is as in \cref{eq:SpnS_monoidal}.
\end{proof}

Therefore the above algebra of continuous machines is in fact the algebra of $\Shv{\Int}$-span
systems of \cref{def:CSpanSystems}.
Notice that even though for the topos $\Shv{\Int}$, the functor $\ADS(\phi)$ \cref{eq:slice_functors} has a right adjoint and preserves pullbacks,
it does not preserve the terminal object,
so it is again the inverse image part of a pre-geometric morphism but not of a geometric morphism.

\subsection{Total, deterministic, and inertial machines}\label{sec:total_det_inert}

Our goal in this section is to describe certain subclasses of continuous machines, \cref{def:ADS}.
Suppose a continuous machine is in a state germ $s_0\in S(0)$ corresponding to a certain input germ $a_0\coloneqq\ii{p}(s_0)\in A(0)$;
these can be thought as initial instantaneous values for the state and input of the machine.
If $a_0$ is extended to some longer input behavior $a\in A(\ti)$, e.g. information flows into the machine,
$s_0$ may or may not extend to some state behavior $s$ that accommodates that input $a$, i.e.\ with $\ii{p}(s)=a$, meaning that
the machine `runs'. There may be more than one such extension, or none at all.

\begin{equation}\label{eq:total_det}
\begin{tikzpicture}[baseline=(topdot)]
\draw[latex-latex] (-.5,0) -- (4.5,0);
\draw[shift={(0,0)},color=black] node[below=3pt] {0};
\draw[shift={(4,0)},color=black] node[below=3pt] {$\ti$};
\draw[{[-]}] (0,0) -- (4,0);
\node[black] (topdot) at (0,1.2) {\textbullet};
\node[black] at (-.3,1.2) {$\scriptstyle{s_0}$};
\node[black] at (0,0.2) {\textbullet};
\node[black] at (-.3,.2) {$\scriptstyle{a_0}$};
\draw[thick] (0,0.2) .. controls (.5,.5) and (2,-0.6) ..(3.9,-.1);
\draw[loosely dashed] (0,1.2) sin (.5,1.4) cos (1,1.2) sin (1.5,1) cos (2,1.2) sin (2.5,1.4) cos (3,1.2) sin (3.5,1) cos (3.9,1.1);
\draw[loosely dashed] (0,1.2) to[out=40,in=140] (3.9,1.5);
\node[black] at (3.5,.8) {$\scriptstyle s'$};
\node[black] at (3.5,2)  {$\scriptstyle s$};
\draw[|->] (0,1) -- (0,0.4);
\node[black] at (3.2,-.5) {$\scriptstyle a$};
\node[black](A) at (5.2,0) {$A(\ti)$};
\node[black](S) at (5.2,1.4) {$S(\ti)$};
\draw[->] (S) to node[right] {$\scriptstyle{\ii{p}_{\ti}}$} (A);
\end{tikzpicture}
\end{equation}

The idea is that a machine is \emph{total} if there is at least one state extension, and it is
\emph{deterministic} if there is at most one; their formal description (\cref{def:total_det_map})
will be accomplished by imposing conditions on the input and output sheaf maps. 
In particular, continuous machines are generally neither total
nor deterministic. 

Recall the \emph{$\epsilon$-extension} functor $\Ext{\epsilon}\colon\Shv{\Int}\to\Shv{\Int}$ introduced after
\cref{Ext}, where for a sheaf $A$ we have
$\Ext{\epsilon}A(\ti)=A(\ti+\epsilon)$, as well as the two natural transformations $\lambda$ and $\rho$ \cref{eq:Ext_maps}
of left and right restriction. For any continuous machine $(\ii{p},\oo{p})\colon S\to A\times B$, the pullback 
\begin{displaymath}
 \begin{tikzcd}
S' \ullimit\ar[r]\ar[d] & S\ar[d,"\ii{p}"] \\
\Ext{\epsilon}A\ar[r,"\lambda"'] & A
 \end{tikzcd}
\end{displaymath}
has as elements pairs $(a,s_0)$ of inputs over some length $\epsilon$ and state germs in the $\ii{p}$-fibre of the left restriction
$a|_{[0,0]}$. Therefore this already formalizes the initial state germ - input section situation of \cref{eq:total_det}.
We employ this to formally express the desired conditions on $\ii{p}$.

\begin{proposition}\label{prop:total_det_via_extensions}
Let $p\colon S\to A$ be a continuous sheaf map. For any $\epsilon\geq0$, consider the outer naturality square for $\lambda$ in $\Shv{\Int}$
\begin{equation}\label{dia:extension_lemma}
\begin{tikzcd}[sep=25pt]
\Ext{\epsilon}S \arrow[drr, bend left=25, "\lambda_S"] \arrow[ddr, bend right=30, "\Ext{\epsilon}p"']
\arrow[dr, "h^\epsilon"]&[-25pt]\\[-15pt]
& S' \ullimit\arrow[r]\arrow[d,"p'"']
& S\arrow[d, "p"] \\
& \Ext{\epsilon}A \arrow[r,"\lambda_A"'] & A
\end{tikzcd}
\end{equation}
where $h^\epsilon$ is the universal map to the pullback $S'$. The following are equivalent:
\begin{compactenum}
\item\label[item]{itm_sh_surj_all} for all $\epsilon>0$, the sheaf map $h^\epsilon\colon\Ext{\epsilon}S\to S'$ is an epimorphism (resp.\ monomorphism);
\item\label[item]{itm_fn_surj_all} for all $\epsilon>0$, the function $h^\epsilon_0\colon S(\epsilon)\to S'(0)$ is surjective (resp.\ injective);
\item\label[item]{itm_sh_surj_some} there exists $\delta>0$ such that for all $0\leq\delta'\leq\delta$ the sheaf map $h^{\delta'}\colon\Ext{\delta'}S\to S'$ is an epimorphism (resp.\ monomorphism);
\item\label[item]{itm_fn_surj_some} there exists $\delta>0$ such that  for all $0\leq\delta'\leq\delta$ the function $h^{\delta'}_0\colon S(\delta')\to S'(0)$ is surjective (resp.\ injective).
\end{compactenum}
\end{proposition}

\begin{proof}
When $\epsilon=0$, we have that $\lambda_A=\id$ and $\lambda_S=\id$, so $h^\epsilon=\id$. Thus it is clear that \cref{itm_sh_surj_all}\implies\cref{itm_sh_surj_some} and that
\cref{itm_fn_surj_all}\implies\cref{itm_fn_surj_some}, e.g.\ choose $\delta\coloneqq1$. By \cref{lemma:ext}, the sheaf map
$h^\epsilon$ is an epi (resp.\ mono) if and only if for all
$\ti\in\Int$ the function $h^\epsilon_\ti$ is epi (resp.\ mono). In particular, taking $\ti=0$ we have \cref{itm_sh_surj_all}\implies\cref{itm_fn_surj_all} and
\cref{itm_sh_surj_some}\implies\cref{itm_fn_surj_some}.

It thus suffices to prove
\cref{itm_fn_surj_some}\implies\cref{itm_sh_surj_all};
so assume \cref{itm_fn_surj_some} holds for some $\delta>0$, and choose $\epsilon>0$. Any $s'\in S'(\ti)$ can be identified with a pair
$(s,a)\in S(\ti)\times A(\ti+\epsilon)$
such that $p(s)=\ShRst{a}{0}{\ti}$. For the surjectivity claim, we have that $h^{\delta'}_0$ is surjective for all $0\leq\delta'\leq\delta$, and the goal is to show $h^\epsilon_\ti$ is surjective; it suffices to find an extension $\bar{s}\in S(t+\epsilon)$ of $s$ over $a$.

Since $\RR$ is a Euclidean domain, there exists a unique $N\in\NN$ and $0\leq\delta'\leq\delta$
such that $\epsilon=N\delta+\delta'$.
Then $(\ShRst{a}{\ti}{\ti+\delta'},\scod{s})\in S'(0)$
and since $h^{\delta'}_0$ is surjective, there exists an extension
$s_0'\in S(\delta')$ emanating from $\scod{s}$ above $\ShRst{a}{\ti}{\ti+\delta'}$.
Gluing this to $s$ we obtain
$s_0\coloneqq s*s_0'\in S(\ti+\delta')$ extending $s$ over $a_0\coloneqq\ShRst{a}{0}{\ti+\delta'}$.

If $N=0$ we are done, so we proceed by induction on $N$. Given $s_n\in S(\ell+\delta'+n\delta)$,
extending $s$ over $a_n\coloneqq\ShRst{a}{0}{\ti+\delta'+n\delta}$ we use the surjectivity of $h^\delta_0$
to extend it, exactly as above, to find $s_n'$ emanating from $\sdom{s_n}$ over $a_n$. At the end,
we have the desired lift $s_N\in S(\ti+\epsilon)$.

Assume now that $h_0^{\delta'}$ is injective for all $0\leq\delta'\leq\delta$. Then there is at most one extension
$s_0'\in S(\delta')$ emanating from $\scod{s}$ over $\ShRst{a}{\ti}{\ti+\delta'}$. Thus there is at most one $s_0\in S(\ti+\delta')$
extending $s$ over $\ShRst{a}{0}{\ti+\delta'}$. The proof concludes analogously to the surjective case above.
\end{proof}

A generalization and purely formal proof of the above result is given in \cref{lemma:small_extensions_big}.

\begin{definition}\label{def:total_det_map}
We will say that a sheaf morphism $p\colon S\to A$ is \emph{total}
(resp.\ \emph{deterministic}) if it satisfies the equivalent conditions of \cref{prop:total_det_via_extensions}.
\end{definition}

\begin{remark}
In \cref{appendix:DSF}, we discuss an equivalence of categories between $\Int$-sheaves
and discrete \Conduche fibrations over the monoid $\RRnn$.
Under that correspondence, a morphism $p\colon S\to A$  in $\Shv{\Int}$ is both total and deterministic
if and only if its associated functor $\ascat{p}\colon\ascat{S}\to\ascat{A}$
is a discrete opfibration; see \cref{sec:B3}.
\end{remark}

It would be reasonable to define a continuous machine $(\ii{p},\oo{p})\colon S\to A\times B$ to be total
(resp.\ deterministic) if and only if the input map $\ii{p}$ is, as this matches our intuitive notion \cref{eq:total_det}.
Indeed, condition \cref{itm_fn_surj_all} says that
for every pair $(s_0,a)\in S(0)\times A(\epsilon)$ with $a|_{[0,0]}=\ii{p}_0(s_0)$,
there exists some $s\in S(\epsilon)$ above $a$ which extends $s_0$, for every $\epsilon>0$.
However, such a definition turns out to not be closed under feedback composition.
In order for total or deterministic machines to form a $\cat{W}$-algebra, we must add an extra condition,
this time on their output map. We introduce a notion of \emph{inertia},
in which a machine's current state determines
not only its current output, but also a small amount of its future output.

\begin{definition}\label{def:inertial_map}
A sheaf map $p\colon S\to B$ is called \emph{$\epsilon$-inertial} when there exists a factorization
\begin{equation}\label{eq:inertial}
\begin{tikzcd}
	&
	\Ext{\epsilon}B\ar[d, "\lambda"]\\
	S\ar[r, "p"']\ar[ur, dashed, "\ol{p}"]&
	B
\end{tikzcd}
\end{equation}
through the $\epsilon$-extension of $B$ via a sheaf map $\ol{p}$; it is called \emph{inertial} when it is $\epsilon$-inertial for some $\epsilon>0$. Explicitly, $\ol{p}_\ti(s)|_{[0,\ti]}=p_\ti(s)\in B(\ti)$ for any section $s\in S(\ti)$.

We say a continuous machine $(\ii{p},\oo{p})\colon S\to A\times B$  is  \emph{inertial}
(resp.\ \emph{$\epsilon$-inertial}) if the output map $\oo{p}$ is. Inside $\ADS(A,B)$, the full category of
inertial (resp.\ $\epsilon$-inertial) machines is denoted $\iADS(A,B)$ (resp.\ $\iADSe(A,B)$).
\end{definition}

The following expression of inertiality via a lifting in a naturality square is equivalent, and may look more intuitive.

\begin{lemma}
A sheaf map $p$ is $\epsilon$-inertial when there exists a lift $\ol{p}$ as in
\begin{equation}\label{eqn.lifting}
\begin{tikzcd}[column sep=large]
	\Ext{\epsilon}S\ar[r, "\Ext{\epsilon}p"]\ar[d, "\lambda_B"']&
	\Ext{\epsilon}B\ar[d, "\lambda_S"]\\
	S\ar[r, "p"']\ar[ur, dashed, "\ol{p}"]&
	B
\end{tikzcd}
\end{equation}
\end{lemma}
\begin{proof}
By \cref{def:inertial_map}, it suffices to show that if $p$ is $\epsilon$-inertial then the top triangle commutes. Since $(\Ext{\epsilon}p)_\ell=p_{\ell+\epsilon}$ for $\ell\in\Int$, we need to show that the lower triangle in the following square commutes: 
\begin{displaymath}
\begin{tikzcd}[column sep=50pt]
	&
	B(\ell+2\epsilon)\ar[d, "\lambda_{\ell+\epsilon}"]\\
	S(\ell+\epsilon)\ar[ur, "\ol{p}_{\ell+\epsilon}"]\ar[r, "p_{\ell+\epsilon}"]\ar[d, "\lambda_\ell"']&
	B(\ell+\epsilon)\\
	S(\ell)\ar[ru, "\ol{p}_\ell"']
\end{tikzcd}
\end{displaymath}
The top triangle commutes because $p$ is $\epsilon$-inertial, and the outer square commutes because $\ol{p}$ is a sheaf map, therefore has natural components; thus the lower triangle does too.
\end{proof}

\begin{remark}
In most cases of interest, if a map $p\colon S\to B$ is inertial then the lift $\ol{p}$ in \cref{eq:inertial} is unique. Indeed, suppose $\lambda\colon \Ext{\epsilon}S\to S$ is surjective, meaning that any short section can be extended in some way; we refer to such sheaves $S$ as \emph{extensible}. In this case, by \cref{eqn.lifting}, there is at most one $\ol{p}$ with $\ol{p}\circ\lambda=\Ext{\epsilon}p$.
\end{remark}

Note that if $p$ is $\epsilon$-inertial and $\delta\leq\epsilon$,
then $p$ is $\delta$-inertial as well. Thus, for any pair of sheaves $A,B$ there is a fully faithful
functor $\iADSe(A,B)\to\iADSd(A,B)$ given by left restriction, and the colimit of the directed system is $\iADS$:
\begin{equation}\label{eq:directedcolimitinertial}
\iADS(A,B)\cong\colim_{\epsilon>0}\left(\iADSe(A,B)\right).
\end{equation}

We are now ready to define total and deterministic machines, in such a way that they are closed under all wiring diagram operations
(\cref{prop:tADS}).

\begin{definition}\label{def:total_det_machine}
A continuous machine $(\ii{p},\oo{p})\colon S\to A\times B$ is \emph{total}
(resp.\ \emph{deterministic}) if
\begin{compactitem}
\item the input map $\ii{p}$ is total (resp.\ deterministic) in the sense of \cref{def:total_det_map}; and
\item the output map $\oo{p}$ is inertial in the sense of \cref{def:inertial_map}.
\end{compactitem}
We denote the full subcategory of $\ADS(A,B)$ of total, resp.\ deterministic, $(A,B)$-machines by $\tADS(A,B)$, resp.\ $\dADS(A,B)$.
Their intersection $\tdADS(A,B)$ is the set of machines that
are both total and deterministic, i.e.\ machines for which the sheaf map $h^\epsilon$ in \cref{dia:extension_lemma} is an isomorphism.
\end{definition}

\begin{example}[Delay Box]\label{ex:delay_box}
For any sheaf $B$, we can define a total and deterministic $(B,B)$-machine $\SmallBox{B}{B}{D^\epsilon}$
that takes input of type $B$ and delays it for time $\epsilon$: it is the span
\[B\From{\rho}\Ext{\epsilon}B\To{\lambda}B\]
in $\Shv{\Int}$, where $\rho$ and $\lambda$ are the right and left restrictions as in \cref{eq:Ext_maps}.
It is clear that $\lambda$ is $\epsilon$-inertial. To see that $\rho$ is both total and
deterministic, one checks that the map $h^\epsilon\colon\Ext{2\epsilon}B\to\Ext{\epsilon}B\times_B\Ext{\epsilon} B$ in
\cref{dia:extension_lemma} is an isomorphism.

\end{example}

The following examples show that the conditions of totality and determinism on input sheaf maps of machines with non-inertial output maps are indeed not closed under arbitrary machines nesting, as implied earlier.

\begin{example}\label{ex:necc_and_suff}
Consider the wiring diagram $\phi\colon X\to 0$, see \cref{eq:wiringdiagpic}, shown below
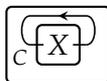
\[
\begin{tikzpicture}[oriented WD, bb small]
	\node[bb={1}{1}] (X) {$X$};
	\node[bb={0}{0},fit={($(X.north west)+(0,1)$) ($(X.south east)+(0,-1)$)}] (Y){};
	\draw[ar] let \p1=(X.north west), \p2=(X.north east), \n1=\bbportlen, \n2={\y1+\bby} in
		(X_out1) to[in=0] (\x2+\n1,\n2) -- (\x1-\n1,\n2) to[out=180] (X_in1);
	\draw[label]
		node [below left=1.5pt and .5pt of X_in1] {$C$};
\end{tikzpicture}
\]
which consists of $(\fin\colon\1\times C\To{\sim}C,\fout\colon C\To{!}\1)$ in $\Shv{\Int}$ as in \cref{eq:wiring_diag_Int}, where $\1$ is the terminal sheaf. Any inhabitant continuous machine $P\coloneqq (\ii{p},\oo{p})\colon S\to C\times C$ becomes, via \cref{eq:SpnS_on_arrows}, $\ADS(\phi)(P)=T\to\1\times\1$ with new sheaf of states given by $T(\ti)=\{s\in S(\ti)\;|\;\ii{p}_\ti(s)=\oo{p}_\ti(s)\}$ as in \cref{eq:newstatesheaf}, namely the collection of all state sections that have the same input and output.

Take $C=\Fnc(\ul{2})$ to be the sheaf of functions on the set $\ul{2}\coloneqq\{1,2\}$ as described in
\cref{Int_sheaves_ex}(\ref{sheaf_of_functions}), given by $C(\ti)=\{f\colon[0,\ti]\to\{1,2\}\}$.
The machine $(\id_C,\id_C)\colon C\to C\times C$ is such that its input map is deterministic (and total); indeed, given any input function $f\colon[0,\ti]\to\{1,2\}$ and an initial state germ above it, namely $f(0)$, then there exists a unique state extension mapping to $f$ via $\ii{p}$, namely $f$ itself. However, its composite $\ADS(\phi)(P)$ has the same state sheaf $C$ since all functions have the same input and output via the identities, whereas $C\to\1\times\1$ does not have a deterministic input anymore: there are multiple functions that only agree on $0$. 


Now take $C=\yoneda{1}$ to be the representable sheaf on the interval $1=[0,1]$ as described in \cref{Int_sheaves_ex}(\ref{yoneda}), given by $\yoneda{1}(\ti)=\Int(\ti,1)$. Notice that for lengths $\ti>1$, these sets are empty. Again, the machine $(\id_C,\id_C)\colon C\to C\times C$ has total (and deterministic), however its composite $C\to\1\times\1$ does not have a total input map anymore: the diagram \cref{dia:extension_lemma} would require for example $\lambda\colon\Ext{\epsilon}(C)(1)=C(1+\epsilon)\to C(1)$ to be surjective, when $C(1+\epsilon)=\emptyset$.

Notice that in both cases above, the output map $\id_C$ was not inertial. For example, suppose that there exists a factorization
$\begin{tikzcd} C\ar[r,"\ol{p}"description] & \Ext{\epsilon}(C)\ar[r,"\lambda"description] & C \end{tikzcd}$ of the identity sheaf morphism of the first machine. Considering the components $\ol{p}_1$ and $\ol{p}_2$ for $\epsilon=1$, naturality implies the commutativity of
\begin{displaymath}
 \begin{tikzcd}
  C(2)\ar[r,"\lambda_1"]\ar[d,"\ol{p}_2"'] & C(1)\ar[d,"\ol{p}_1"] \\
  C(3)\ar[r,"\lambda_2"] & C(2)
 \end{tikzcd}
\end{displaymath}
where the lower composite is the identity on $C(2)$; that would imply that the left restriction $\lambda$ is injective, which is clearly false.
\end{example}

\subsection{Closure under feedback}\label{sec:feedback}

For the remainder of this section, we write $\cat{W}$ to denote the symmetric monoidal category
$\WW{\Shv{\Int}}$ of $\Shv{\Int}$-labeled boxes and wiring diagrams. We saw in
\cref{prop:ADS} that continuous machines form a $\WW{}$-algebra $\ADS\colon\WW{}\to\Cat$;
now we show that inertial (\cref{def:inertial_map}),
as well as total machines and deterministic machines (\cref{def:total_det_machine}) do too.

\begin{proposition}\label{prop:iADS}
Inertial continuous machines (resp.\ $\epsilon$-inertial machines, for any $\epsilon>0$) 
form a subalgebra $\iADS\colon\cat{W}\to\Cat$ (resp.\ $\iADSe$). Moreover, we have a colimit of algebras
\[\iADS\cong\colim_{\epsilon>0}(\iADSe).\] 
\end{proposition}

\begin{proof}
Given a wiring diagram morphism $\phi\colon X\to Y$, we restrict \cref{eq:Mchphi} to obtain a functor
\[\iADSe(X)\ss\ADS(X)\To{\ADS(\phi)}\ADS(Y);\]
we will show that it factors through the full subcategory $\iADSe(Y)$, for the same $\epsilon>0$.
Suppose that $p=(\ii{p},\oo{p})\colon S\to\Xin\times\Xout$ is a machine and that $\oo{p}$
factors as $\lambda\circ\ol{p}$ \cref{eq:inertial}.
This is shown in the upper-middle square below
\begin{equation}\label{eq:compositeinertial}
\begin{tikzcd}[column sep=.4in,row sep=.2in]
T\ar[r,"k"']\ar[dd,"{(\ii{q},\oo{p}k)}"']\ullimit\ar[rrr,dashed,bend left=20,"\ol{q}"description] &
S\ar[r,dashed,"\ol{p}"']\ar[d,"{(\ii{p},\oo{p})}"']
& \Ext{\epsilon}(\Xout)\ar[r,"\Ext{\epsilon}(\fout)"']\ar[d,"\lambda"]
&[10pt] \Ext{\epsilon}(\Yout)\ar[dd,"\lambda"] \\
& \Xin\times\Xout\ar[r,"\pi_2"'] & \Xout\ar[dr,"\fout"] & \\
\Yin\times\Xout\ar[rr,"1\times\fout"']\ar[ur,"{(\fin,\pi_2)}"] && \Yin\times\Yout\ar[r,"\pi_2"'] & \Yout
\end{tikzcd}
\end{equation}
The construction for $\ADS(\phi)(S,p)=(T,q)$ as in \cref{eq:SpnS_on_arrows}
is the left-hand pullback, along with the bottom-left horizontal morphism, thus defining $q=(\ii{q},\oo{q})\colon T\to\Yin\times\Yout$.
The bottom (trapezoid-shaped) sub-diagram trivially commutes, and the
right-hand diagram commutes because $\lambda$ is a natural transformation. Thus the
outer square commutes, exhibiting a factorization of $\oo{q}$ through $\Ext{\epsilon}(\Yout)$.
Hence we do have a subfunctor $\iADSe\colon\WW{}\to\Cat$.

If $0<\delta\leq\epsilon$, the following diagram of categories commutes
\[
\begin{tikzcd}[column sep=.6in]
	\iADSe(X)\ar[r,"\iADSe(\phi)"]\ar[d]&\iADSe(Y)\ar[d]\\
	\iADSd(X)\ar[r,"\iADSd(\phi)"]&\iADSd(Y)
\end{tikzcd}
\]
so we have a natural transformation $\iADSe\Rightarrow\iADSd$, and the colimit of the directed system is $\iADS$ by \cref{eq:directedcolimitinertial}.
It remains to show that all these functors $\WW{}\to\Cat$, and maps between them, are monoidal. 

It is easy to see that the product \cref{eq:SpnS_monoidal} of $\epsilon$-inertial machines is $\epsilon$-inertial,
because $\Ext{\epsilon}$ preserves products of sheaves, see \cref{lemma:ext_pbs}.
Thus $\iADSe$ is indeed a $\WW{}$-algebra, and moreover
the natural transformation $\iADSe\Rightarrow\iADSd$ described above is clearly monoidal. We still
need to define a monoidal structure on the colimit $\iADS$: given output sheaf maps that factorize through $S\to\Ext{\epsilon_1}(\Xout)$ and
$P\to\Ext{\epsilon_2}(\Zout)$ for $\epsilon_1,\epsilon_2>0$, let $\epsilon\coloneqq\min(\epsilon_1,\epsilon_2)$. Then by left restriction (if necessary), the product output sheaf map factorizes through
$S\times P\to\Ext{\epsilon}(\Xout\times\Zout)$ therefore is also inertial for this new $\epsilon>0$.
\end{proof}

\begin{proposition}\label{prop:tADS}
Total machines 
form a subalgebra $\tADS\colon\cat{W}\to\Cat$. Similarly,
deterministic machines form a subalgebra $\dADS\colon\cat{W}\to\Cat$. Their intersection is a subalgebra $\tdADS\coloneqq\tADS\cap\dADS$.
\end{proposition}

\begin{proof}
Given an object $X=(\inp{X},\out{X})$, define the full subcategories $\tADS(X),\dADS(X)$
of $\ADS(\Xin,\Xout)$ as in \cref{def:total_det_machine}. 
Given a morphism $\phi\colon X\to Y$,
we restrict to obtain functors
\[
\tADS(X)\ss\ADS(X)\To{\ADS(\phi)}\ADS(Y) \qquad \mathrm{and} \qquad
\dADS(X)\ss\ADS(X)\To{\ADS(\phi)}\ADS(Y)
\]
described in \cref{prop:ADS}; we need to show that they factor through the respective full subcategories $\tADS(Y)$ and $\dADS(Y)$, and
the result for $\tdADS$ follows trivially. 

A total (resp.\ deterministic) machine $p=(\ii{p},\oo{p})\colon S\to\Xin\times\Xout$ is
$\epsilon$-inertial for some $\epsilon>0$, meaning that $\oo{p}$ is a composite
$\lambda\circ\ol{p}\colon S\to\Ext{\epsilon}\Xout\to \Xout$ \cref{eq:inertial}. By \cref{prop:iADS}, its image $q=(\ii{q},\oo{q})\colon T\to\Yin\times\Yout$ under $\ADS(\phi)$ is also $\epsilon$-inertial for the same $\epsilon$, so $\oo{q}=\lambda\circ \ol{q}$ for some $\ol{q}$ as exhibited by \cref{eq:compositeinertial}. In particular, from the top-right part of that diagram we deduce that this factorization coincides with one through $\Xout$, namely
\begin{equation}\label{eq:outqfactor}
\oo{q}=T\To{k}S\To{\ol{p}}\Ext{\epsilon}\Xout\To{\lambda}\Xout\To{\out{\phi}}\Yout.
\end{equation}
It suffices to check that the resulting input sheaf map $\ii{q}$ is total (resp.\ deterministic). 
We use the formulation in \cref{prop:total_det_via_extensions}(\cref{itm_sh_surj_some}); that is, we will show that 
the sheaf map $h^{\epsilon'}$ in \cref{dia:extension_lemma}
is an epimorphism (resp.\ monomorphism) for arbitrary $\epsilon'\leq\epsilon$.
But note that if $\oo{q}$ is $\epsilon$-inertial then it is also
$\epsilon'$-inertial, so we may take $\epsilon'=\epsilon$.
First of all, the bottom squares below are pullbacks in a straightforward way; we define $S'$ and $T'$ to be the top pullbacks
\begin{equation}\label{eqn:two_pbs}
\begin{tikzcd}[row sep=.2in]
	|[alias=a]|S'\ar[r]\ar[dd]
	\ar[ddd,dashed,pos=.3,"p'",to path={.. controls +(left:1.5) and +(left:1.5).. (\tikztotarget.west) \tikztonodes}]\ar[dr,phantom,pos=.05,bend right=20,"\lrcorner"']&S\ar[dd,"{(\ii{p},\ol{p})}"]\\
	& \phantom{A} \\
	\Ext{\epsilon}\Xin{\times}\Ext{\epsilon}\Xout\ar[dr,phantom,pos=.05,bend right=20,"\lrcorner" description]\ar[r,"\lambda\times1"']\ar[d,"\pi_1"']&\Xin{\times}\Ext{\epsilon}\Xout\ar[d,"\pi_1"]\\
	|[alias=b]|\Ext{\epsilon}\Xin\ar[r,"\lambda"']&\Xin
\end{tikzcd}
\;\;
\begin{tikzcd}[row sep=.2in]
	T'\ar[ddd,dashed,pos=.3,"q'",to path={.. controls +(left:1.5) and +(left:1.5).. (\tikztotarget.west) \tikztonodes}]
	\ar[r]\ar[dd]\ar[dr,phantom,pos=.05,bend right=20,"\lrcorner" description]&T\ar[d,"{(\ii{q},k)}"]\\
	& \Yin\times S\ar[d,"1\times\ol{p}"] \\
	\Ext{\epsilon}\Yin{\times}\Ext{\epsilon}\Xout\ar[dr,phantom,pos=.05,bend right=20,"\lrcorner" description]\ar[r,"\lambda\times1"']\ar[d,"\pi_1"']&\Yin{\times}\Ext{\epsilon}\Xout\ar[d,"\pi_1"]\\
	\Ext{\epsilon}\Yin\ar[r,"\lambda"']&\Yin
\end{tikzcd}
\end{equation}
Note that the map $p'$ on the left is the one defined in \cref{prop:total_det_via_extensions}, and similarly is $q'$ if we observe that the top pullbacks are along identities with respect to the second components. Therefore the outer pullbacks coincide with the ones in \cref{dia:extension_lemma} for the machines $S$ and $T$. Now consider the diagram
\[
\begin{tikzcd}[row sep=.17in,column sep=.2in]
	T'\ar[rr]\ar[dr]\ar[ddd]&&S'\ar[dr]\ar[ddd]\\
	&T\ar[rr,crossing over,pos=.4,"k"]\ar[d,crossing over,"{(\ii{q},k)}"]&&S\ar[ddd,"{(\ii{p},\ol{p})}" {near start,description}]\\
	&\Yin\times S &&&\\
	\Ext{\epsilon}\Yin\times\Ext{\epsilon}\Xout\ar[rr]\ar[dr]&&\Ext{\epsilon}\Xin\times\Ext{\epsilon}\Xout\ar[dr]\\
	&\Yin\times\Ext{\epsilon}\Xout\ar[rr,crossing over]\ar[from=uu,crossing over,pos=.7,"1\times\ol{p}"]
	&&\Xin\times\Ext{\epsilon}\Xout\ar[d,"1\times\lambda"]\\
	&\Yin\times\Xout\ar[rr,"{(\inp{\phi},\pi_2)}"]\ar[from=u,crossing over,"1\times\lambda"]\ar[d,"1\times\out{\phi}"]&&\Xin\times\Xout \\
	&\Yin\times\Yout
\end{tikzcd}
\]
where the two sides are precisely the pullbacks of \cref{eqn:two_pbs}. Notice how the two front vertical composites describe the $\epsilon$-inertial machine $S\to\Xin\times\Xout$ on the right, and its resulting (also $\epsilon$-inertial) machine $T\to\Yin\times\Yout$ on the left by \cref{eq:outqfactor}. Now the front rectangle is a pullback by definition \cref{eq:SpnS_on_arrows} of machine composition, and one can check that the lower bottom square is also a pullback. Therefore the back face is a pullback, which we re-write as the bottom square below: 
\begin{equation}\label{eq:finalpb}
\begin{tikzcd}
	\Ext{\epsilon}T\ar[dr,phantom,pos=.05,bend right=20,"\lrcorner" description]\ar[r]\ar[d,"i^\epsilon"']&\Ext{\epsilon}S\ar[d,"h^\epsilon"]\\
	T'\ar[r]\ar[dr,phantom,pos=.05,bend right=20,"\lrcorner" description]\ar[d]&S'\ar[d]\\
	\Ext{\epsilon}\Yin\times\Ext{\epsilon}\Xout\ar[r]&\Ext{\epsilon}\Xin\times\Ext{\epsilon}\Xout
\end{tikzcd}
\end{equation}
Since $\Ext{\epsilon}$ preserves limits by \cref{lemma:ext_pbs}, the big rectangle is a pullback, so the top is too.
This diagram takes place in a topos---in particular, a regular category---so if $h^{\epsilon}$ is epi (resp.\ mono) then $i^{\epsilon}$ is too. Therefore if $p$ is a total (resp.\ deterministic) machine, so is $\ADS(\phi)(p)$.

Finally, it is easy to check that if $\ii{p}\colon P\to\Xin$
and $\ii{q}\colon T\to\Zin$ are total (resp.\ deterministic),
their product $\ii{p}\times\ii{q}$ is too, so lax monoidality \cref{eq:SpnS_monoidal} follows.
\end{proof}

Notably, since the above conditions are only sufficient and not necessary --- the pullback \cref{eq:finalpb} may produce an epimorphism or monomorphism $i^\epsilon$ on the left even if the initial $h^\epsilon$ on the right is not --- it could be the case that a composite machine is total and deterministic without its sub-components being so. Even more, the two cases of \cref{ex:necc_and_suff} coincidentally maintain their totality and determinism respectively, despite not being inertial. In this work, our goal was to establish sufficient conditions for totality and determinism to be carried over from subsystems to their composite, and this is what was accomplished in \cref{prop:tADS}.

\subsection{Discrete and synchronous variations}\label{sec:discrete_variations}

One of the central purposes of the current work is to allow
discrete-time and continuous-time systems to be incorporated within the same
framework. In this section, we define discrete analogues of the machines described so far; all the results go through with nearly identical proofs, so we only provide the definitions and statements. In the following \cref{ch:maps_to_machines}, we compare the discrete and continuous machines via wiring diagram algebra maps.

Discrete machines have $\Int_N$-sheaves, rather than $\Int$-sheaves, labeling their input and output ports; see \cref{sec:Int_sheaves} were continuous and discrete interval sheaves were introduced.
All of the ideas of \cref{sec:total_det_inert,sec:feedback} can be modified to this setting,
with the main difference being that the positive real $\epsilon\in\RRnn$ used in the various extension sheaves (\cref{Ext}) is here replaced by the positive integer $1\in\NN$. For example, compare the following to \cref{prop:total_det_via_extensions,def:inertial_map,def:total_det_machine}.

\begin{definition}\label{def:ADS[N]}
Let $p\colon S\to A$ be a discrete sheaf map. Consider the outer naturality square for $\lambda$ in $\Shv{\Int_N}$

\[
\begin{tikzcd}[sep=25pt]
\Ext{1}S \arrow[drr, bend left=25, "\lambda"] \arrow[ddr, bend right=30, "\Ext{1}p"'] \arrow[dr, "h^1"]&[-25pt]\\[-15pt]
& S' \ullimit\arrow[r]\arrow[d,"p'"']
& S\arrow[d, "p"] \\
& \Ext{1}A \arrow[r,"\lambda"'] & A
\end{tikzcd}
\]
where $h_1\colon\Ext{1}S\to S'$ is the universal map to the pullback $S'$.
Then $p$ is \emph{total} (resp.\ \emph{deterministic}) if $h_1$ is an epimorphism (resp.\ monomorphism); $p$ is \emph{inertial} if it
factors through $\lambda\colon\Ext{1}A\to A$.

A \emph{discrete machine} is a span $p=(\ii{p},\oo{p})\colon S\to A\times B$ in
$\Shv{\Int}_N$. The machine is \emph{inertial} if the output map $\oo{p}$ is.
The machine \emph{total} (resp.\ \emph{deterministic}) if both the output map $\oo{p}$ is inertial and the input map
$\ii{p}$ is \emph{total} (resp.\ \emph{deterministic}).
\end{definition}

\begin{remark}\label{rem:simplifications_for_graphs}
By an argument similar to that in \cref{prop:total_det_via_extensions}, we can simplify \cref{def:ADS[N]} a great deal.
Recall from \cref{prop:ShIntN_grphs} that $p\colon S\to A$ can be identified with a graph homomorphism;
write $S=(\begin{tikzcd}[sep=18pt]S_1\ar[r,shift left=2pt,"\src"]\ar[r,shift right=2pt,"\tgt"']&S_0\end{tikzcd})$
where we let $S_n\coloneqq S(n)$ denote the set of length-$n$ sections of $S$,
and similarly for $A$. Consider the diagram of sets
\begin{equation}\label{eqn:tot_det_graphs}
\begin{tikzcd}[sep=25pt]
S_1 \arrow[drr, bend left=25, "\src"] \arrow[ddr, bend right=30, "p_1"'] \arrow[dr, "h"]&[-20pt]\\[-10pt]
& S' \ullimit\arrow[r]\arrow[d]
& S_0\arrow[d, "p_0"] \\
& C_1 \arrow[r,"\src"'] & C_0
\end{tikzcd}
\end{equation}
where $S'$ is the pullback and $h$ is the induced function.
Then $p$ is total (resp.\ deterministic) if and only if $h$ is surjective (resp.\ injective). Notice that $h=h^1_0$ from above.
\end{remark}

The proof of \cref{prop:ADS_N} follows that of \cref{prop:ADS,prop:iADS,prop:tADS}.
In particular, $\ADS[N]\colon\cat{W}_{\Shv{\Int}_N}\to\Cat$ is the algebra of $\Shv{\Int}_N$-span systems of
\cref{def:CSpanSystems}.

\begin{proposition}\label{prop:ADS_N}
Discrete machines form a $\cat{W}_{\Shv{\Int}_N}$-algebra $\ADS[N]\colon\cat{W}_{\Shv{\Int}_N}\to\Cat$. Inertial, total, deterministic discrete machines form subalgebras
$\iADS[N]$, $\tADS[N]$, $\dADS[N]$ (and $\tdADS[N]$).
\end{proposition}

To compare discrete and continuous systems in \cref{ch:maps_to_machines}, it will be useful to have a mediating construct that
handles `synchronization' in continuous systems.
Recall the synchronizing sheaf $\Sync\in\Shv{\Int}$ \cref{sync(l)} and the synchronous sheaves $X\to\Sync$
from \cref{sec:synchronization}. Given an object (resp.\ map, span) in
$\Shv{\Int}/\Sync$, we refer to its image under the forgetful
$\Shv{\Int}/\Sync\to\Shv{\Int}$ as the \emph{underlying} object
(resp.\ map, span) in $\Shv{\Int}$.

\begin{definition}\label{def:ADS_Sync}
A morphism $p\colon S\to A$ in $\Shv{\Int}/\Sync$ is \emph{inertial}
(resp.\ \emph{total, deterministic}) if the underlying map in $\Shv{\Int}$ is.

A \emph{synchronous machine} is a span $p=(\ii{p},\oo{p})\colon S\to A\times B$ in
$\Shv{\Int}/\Sync$. It is \emph{inertial} (resp.\ \emph{total, deterministic}) if the underlying
continuous machine $p$ in $\Shv{\Int}$ is.
\end{definition}

Again, the proof of \cref{prop:ADS_Sync} is similar to the analogous results proven above,
and in particular $\ADS[\Sync]\colon\cat{W}_{\Shv{\Int}/\Sync}\to\Cat$ is the algebra of $\Shv{\Int}/\Sync$-span systems of
\cref{def:CSpanSystems}.

\begin{proposition}\label{prop:ADS_Sync}
Synchronous machines form a $\cat{W}_{\Shv{\Int}/\Sync}$-algebra
$\ADS[\Sync]\colon\cat{W}_{\Shv{\Int}/\Sync}\to\Cat.$
Inertial, total, deterministic synchronous machines form subalgebras
$\iADS[\Sync]$, $\tADS[\Sync]$, $\dADS[\Sync]$ (and $\tdADS[\Sync]$).
\end{proposition}

Evidently, continuous machines receive input information continuously through time, discrete
machines receive information in discrete `ticks' $0,1,2,..$ of a global clock, whereas synchronous machines
also have a continuous input flow, this time together with an assigned phase $\theta\in[0,1)$.
In \cref{cor:1} we will show that there are algebra morphisms realizing any discrete or continuous machine as a synchronous machine.

\subsection{Safety contracts}\label{safety_contracts}

We conclude this chapter by a short discussion on
\emph{contracts} that we can impose on machines.
Intuitively, a contract for a machine should be a set of logical formulas that
dictate what sort of input/output behavior is \emph{valid} for it.
If the machine inhabits some labeled box $X$, its contract should be expressed relatively to the sections
of the input and output sheaves $\Xin,\Xout$.

For example, suppose we have a box $X$ and consider the following natural language contract for a discrete
machine inhabiting $X$: ``if I ever receive two True's in a row, I will output a False within 5
seconds''.
For any input/output pair $(i,j)\in\Xin(n)\times\Xout(n)$, we can evaluate
whether the machine satisfies the contract, at least on the interval $[0,6]$: if the first two inputs are True,
$\ShRst{i}{0}{1}=\langle \tn{True,\ True}\rangle$, then the last five outputs $\ShRst{j}{2}{6}$ should include a False.
For any section $s$ of longer length $n\geq 6$, we say the machine satisfies the contract if it does on every subinterval $[a,a+6]\ss[0,n]$.

But on shorter sections, it is unclear what to do: should one say the contract is satisfied on short intervals or not?
One choice is what is sometimes called a \emph{safety contract}. These have the property that if a section validates the contract,
then so does any restriction of it. Thus we formalize safety contracts on $X$ as sub-presheaves
$C\ss U(\Xin\times\Xout)\cong U\Xin\times U\Xout$ where $U$ is the forgetful functor into presheaves; namely
if a section $x=(i,j)$ is in $C$ then it is valid.%
\footnote{Using sheaves, rather than presheaves, for safety contracts does not have the intended semantics, because the sheaf condition
would imply that
the concatenation of any two valid behaviors is valid. It thus effectively disallows historical context
from being a consideration in validity. In the example above, any section of length 4 is valid, but some sections of length 8 are not
so the gluing of two valid sections is not always valid.}
Thus we arrive at the following definition.

\begin{definition}\label{def:safety_contract}
Let $A,B\in\Shv{\Int}$ be continuous interval sheaves, and $U\colon\Shv{\Int}\to\Psh{\Int}$
the forgetful functor. A \emph{safety contract} is a sub-presheaf $C\ss UA\times UB$.
We say a section $(a,b)\in A\times B$ \emph{validates} the contract if $(a,b)\in C$.
\end{definition}

Recall that images exist in any topos, in particular in $\Psh{\Int}$. 
If $p=(\ii{p},\oo{p})\colon S\to A\times B$ is a continuous machine as in \cref{def:ADS}, we say that it
\emph{validates} the contract, denoted $p\models C$, if the image $\im(Up)$
is contained in $C$, or equivalently $Up$ factors through it:
\begin{displaymath}
\begin{tikzcd}[row sep=.5in,column sep=.5in]
 & C \ar[d,>->] \\ US\ar[r,"{Up}"']\ar[ur,dashed] & UA\times UB
\end{tikzcd}
\end{displaymath}
In other words, $p$ validates $C$ if, for
every section $s\in S$, the associated input and output validates the contract as above, $\big(\ii{p}(s),\oo{p}(s)\big)\in C$.

The collection of all $(A,B)$-contracts is a poset---in fact a Heyting algebra; it is denoted
\begin{equation}\label{eq:Cntr}
\Cntr(A,B)\coloneqq\SubPsh{A\times B}\equiv\Sub_{\Psh{\Int}}(UA\times UB).
\end{equation}

In what follows, we often elide the forgetful functor $U$
so we may write $C\ss A\times B$ to denote $C\ss UA\times UB$;
in diagrams involving both sheaves and presheaves,
everything should be regarded as a presheaf.

\begin{remark}
As mentioned at the beginning of this section, contracts should really be written in a logical formalism. 
What we call safety contracts in \cref{def:safety_contract} are a reasonable semantics for this formalism;
a related logical formalization was later worked out in \cite{Schultz.Spivak:2017a}. 
\end{remark}

\begin{proposition}\label{prop:Cntr}
Safety contracts form a $\WW{\Shv{\Int}}$-algebra.
\end{proposition}

\begin{proof}
The proof proceeds like the ones for machines or span-like systems earlier, e.g. \cref{prop:ADS}.
To any box $X\in\WW{\Shv{\Int}}$, we associate the posetal category $\Cntr(\Xin,\Xout)$
of \cref{eq:Cntr}.
Given a wiring diagram $\phi\colon X\to Y$ and a contract $C\ss\Xin\times\Xout$, we form $D\coloneqq\Cntr(\phi)(C)$
\begin{equation}\label{eq:Cntr_on_arrows}
\begin{tikzcd}[sep=large]
&[-15pt] C'\ar[r]\ar[d,>->]\ullimit\ar[dl,two heads,bend right] & C \ar[d,>->] \\
D\ar[dr,>->, bend right] &\Yin\times\Xout\ar[r,"{(\fin,\pi_2)}"']\ar[d,"1\times \fout"] & \Xin\times \Xout \\
& \Yin\times\Yout &
\end{tikzcd}
\end{equation}
as the subpresheaf of $\Yin\times\Yout$  which the image of the pullback $C'$.
Analogously to \cref{eq:slice_functors}, if $(\exists_f\dashv f^*)\colon\Sub(B)\to\Sub(A)$
are the induced adjunctions for any $f\colon A\to B$ in the topos $\Psh{\Int}$,
the functor $\Cntr(\phi)$ is the composite
\begin{displaymath}
\SubPsh{\Xin\times\Xout}\xrightarrow{(\fin,\pi_2)^*}
\SubPsh{\Yin\times\Xout}\xrightarrow{\exists_{(1\times\fout)}}\SubPsh{\Yin\times\Yout}.
\end{displaymath}
Functoriality and unitality follow. 
The lax monoidal structure is again given by the cartesian product like \cref{eq:SpnS_monoidal},
i.e.\ natural maps $\SubPsh{X}\times\SubPsh{Z}\to\SubPsh{X\times Z}$.
\end{proof}

Finally, we may consider the inhabitant of a wiring diagram box to be a machine with an associated contract,
namely any contract validated by it. These also form an algebra,
essentially combining \cref{prop:ADS,prop:Cntr}; the composed inhabitant is the composite
machine, validating the associated composite contract.

\begin{definition}\label{def:contracted_machine}
Let $A,B\in\Shv{\Int}$ be sheaves. An \emph{$(A,B)$-contracted machine} is a pair $(S,C)\in\ADS(A,B)\times\Cntr(A,B)$ such that $S\models C$.

Contracted machines form a full subcategory $\CM(A,B)$ of $(\ADS\times\Cntr)(A,B)$.
\end{definition}

\begin{proposition}\label{prop:contracted}
Contracted machines form a $\WW{\Shv{\Int}}$-subalgebra of $\ADS\times\Cntr$.
\end{proposition}
\begin{proof}
For any $X=(\inp{X},\out{X})$, define
$$\CM(X)\coloneqq\CM(\Xin,\Xout)\subseteq\ADS(\Xin,\Xout)\times\Cntr(\Xin,\Xout).$$
Suppose $\phi\colon X\to Y$ is a wiring diagram. Given a contracted machine $(p\colon S\to\Xin\times\Xout,
m\colon C\rightarrowtail U\Xin\times U\Xout)\in\CM(\Xin,\Xout)$, applying $\ADS(\phi)\times\Cntr(\phi)$ produces
a pair $(q\colon T\to\Yin\times\Yout,n\colon D\to U\Yin\times U\Yout)$ given by
\cref{eq:SpnS_on_arrows} and \cref{eq:Cntr_on_arrows}. The following diagram
verifies that $q$ validates $n$:

\begin{displaymath}
\begin{tikzcd}
&[-10pt]T\ullimit\ar[dl,dashed,bend right=10pt,"\exists!"']\ar[dd]\ar[rr] &&
[-10pt]S\ar[dd,"p"]\ar[dl,dashed,bend right=10pt]\\
[-10pt]C'\ullimit\ar[d,two heads,bend right=10pt]\ar[dr,>->,bend right=25pt]\ar[rr]&& C &\\
D\ar[dr,>->,bend right=10pt,"n"']&\Yin\times\Xout\ar[rr]\ar[d]&&
\Xin\times\Xout\ar[from=ul,>->,bend right=25pt,crossing over,"m"]\\
&\Yin\times\Yout
\end{tikzcd}
\end{displaymath}

The left hand side composite is precisely the new machine $q$, due to the upper-side
pullback formed; clearly it factors through $n$. Monoidality is inherited from $\ADS\times\Cntr$,
essentially given by the cartesian product.
\end{proof}

We could of course replace $\ADS$ with any other version of a machine (total, deterministic
and discrete, synchronous)
as described in \cref{sec:total_det_inert,sec:discrete_variations},
and express safety contracts on those. The significance of \cref{prop:contracted} can be summarized as follows: if we
arbitrarily interconnect machines which validate specific contracts, the composite machine
they form is forced to satisfy their composite contract, which can be specified using the wiring diagram algebraic operations.
In other words, we could reason about the expected valid behavior of the total system only by looking at valid
behaviors, and not even state specifications, of the component subsystems.

\section{Maps between dynamical systems and machines}\label{ch:maps_to_machines}

In this final chapter of the main text, we describe wiring diagram algebra maps, i.e.\
monoidal natural transformations \cref{eq:morphWDalg}, between the various algebras we have considered so far.
More specifically, we are interested in maps from discrete and continuous dynamical systems of \cref{sec:DandCDS}
to the machines considered in \cref{ch:ADS},
as well as maps between the different kinds of machines (continuous, discrete, etc.) and contracts.
These maps, in fact embeddings, allow us to translate one sort of algebra to another, while ensuring consistency of serial, parallel,
and feedback composition.

Let us briefly elaborate on the previous sentence for two 
$\WW{}$-algebras, i.e. lax monoidal functors $F,G\colon\WW{}\to\Cat$. Recall that, given
$F$-inhabitants of boxes $f_i\in F(X_i)$ and a wiring diagram $\phi\colon X_1,\ldots,X_n\to Y$
like \cref{picture_wiring_diagram}, we obtain an inhabitant
$F(\phi)(f_1,\ldots,f_n)\in F(Y)$ of the outer box via \cref{eq:inhabitants_composition}.
Given a monoidal natural transformation $\alpha\colon F\Rightarrow G$,
its components $\alpha_{X_i}\colon F(X_i)\to G(X_i)$ map $F$-inhabitants to $G$-inhabitants. The axioms dictate the commutativity of
\begin{displaymath}
\begin{tikzcd}[column sep=.6in,row sep=.5in]
F(X_1)\times...\times F(X_n)\ar[r,"{F_{X_1,...,X_n}}"]\ar[d,"{\alpha_{X_1}\times...\times\alpha_{X_n}}"']
& F(X_1+...+X_n)\ar[r,"{F(\phi)}"]\ar[d,dashed,"{\alpha_{X_1+...+X_n}}"'] & F(Y) \ar[d,"\alpha_Y"] \\
G(X_1)\times...\times G(X_n)\ar[r,"{G_{X_1,...,X_n}}"'] & G(X_1+...+X_n)\ar[r,"{G(\phi)}"'] & G(Y)
\end{tikzcd}
\end{displaymath}
which explicitly means that whether we first compose the interior systems and then translate via $\alpha$ (upper composite),
or first translate and then compose (lower composite), the resulting systems are the same.

We explain in \cref{sec:maps_from_dynamical_systems} how to translate from various sorts of dynamical
systems to machines, and in \cref{sec:maps_between} how to translate between various sorts of machines.
Specifically, \cref{cor:1} establishes synchronous machines $\ADS[\Sync]$ as the common framework where discrete and continuous machines can both be mapped.

\subsection{Realizing dynamical systems as machines}\label{sec:maps_from_dynamical_systems}

Understanding how the motivating examples of discrete and continuous
dynamical systems of \cref{sec:DandCDS} fit into our generalized framework of machines,
makes the latter much more concrete and accomplishes one of the main goals of this work.

To begin with, we realize discrete dynamical systems (\cref{def:DDS}) as discrete machines (\cref{def:ADS[N]}), in fact
total and deterministic ones.
The former are $\Set$-based---the input and output signals are elements
of a set, as are the states---so the first step is to compare the typing categories of the respective wiring diagram algebras, namely sets with
$\Int_N$-sheaves. In \cref{ex:Int_N},
for any set $S$ we defined the $\Shv{\Int}_N$-sheaf $K(S)$,
with length-$n$ sections $K(S)(n)=S^{n+1}=\Hom_\Set(\{0,\ldots,n\},S)$, i.e.\ lists of length $n+1$ in $S$. This mapping
extends to a (finite-limit preserving) functor
\begin{equation}\label{eq:defofK}
K\colon\Set\to\Shv{\Int}_N
\end{equation}
therefore by \cref{eq:W_}, there is an induced strong monoidal functor $\WW{K}\colon\WW{\Set}\to\WW{\Shv{\Int}_N}$
which sends a pair of $\Set$-typed finite sets $(\tau\colon\inp{X}\to\Set,\tau'\colon\out{X}\to\Set)$
to the $\Shv{\Int}_N$-typed finite sets $(K\scirc\tau\colon\inp{X}\to\Shv{\Int}_N,K\scirc\tau'\colon\out{X}\to\Shv{\Int}_N)$.

Recall from \cref{prop:ADS_N} the algebra $\tdADS[N]\colon\WW{\Shv{\Int}_N}\to\Cat$ of total and deterministic discrete machines,
and from \cref{prop:algebra_DDS} the algebra $\DDS\colon\WW{\Set}\to\Cat$ of discrete dynamical systems.

\begin{proposition}\label{DDS_to_ADS}
There exists a morphism of wiring diagram algebras, i.e.\ monoidal natural transformation
\begin{displaymath}
\algmap{\Set}{K}{\Shv{\Int}_N}{\DDS}{\tdADS[N]}{\beta}
\end{displaymath}
\end{proposition}

\begin{proof}
We first define the component functor $\beta_X\colon\DDS(X)\to\tdADS[N](\WW{K}(X))$
for each box $X=(\inp{X},\out{X})\in\WW{\Set}$, mapping a discrete dynamical system on $X$ to a total and deterministic discrete machine
on $\WW{K}(X)$. We then we check that
the $\beta_X$'s satisfy the necessary naturality and monoidality conditions.

Recall that a discrete dynamical system on a $\Set$-labeled box $X$ is a tuple $F=(S,\upd{f},\rdt{f})$, where $S$ is a set,
and $\upd{f}\colon\Xin\times S\to S$ and $\rdt{f}\colon S\to\Xout$ are functions. 
For simplicity, denote by $\inp{\cat{X}}\coloneqq K(\Xin)$ and $\out{\cat{X}}\coloneqq K(\Xout)$
the image of the input and output sets under $K$ \cref{eq:defofK}, namely the discrete interval sheaves
with $n+1$-element lists as sections of length $n$. 

For each discrete dynamical system $F$, we can construct a discrete interval sheaf $\mathcal{S}$ by
\begin{equation}\label{eq:cat_S}
\mathcal{S}(n)\coloneqq
\left\{
(x,s)\colon \{0,\ldots,n\}\to\Xin\times S\;\middle|
\;s_{i+1}=\upd{f}\big(x_{i},s_{i}\big),\quad 0\leq i<n
\right\}
\end{equation}
i.e. finite lists of all pairs
of input and state elements, in the order processed by the dynamical system: the $(i+1)^{\text{th}}$ state is determined by
the $i^{\text{th}}$ state and input.
We can now define $\beta_X(F)$ to be the span $(\ii{p},\oo{p})\colon\mathcal{S}\to\inp{\mathcal{X}}\times\out{\mathcal{X}}$ in $\Shv{\Int}_N$,
where $\ii{p}(x,s)=x$ and $\oo{p}(x,s)=\rdt{f}(s)$; this is the discrete machine that corresponds to $F$. Notice that if $\beta_X(F)=\beta_X(G)$ for two $(\Xin,\Xout)$-discrete dynamical systems, then their state sets and update, readout functions are the same to begin with.

The machine $\beta_X(F)$ is inertial \cref{eq:inertial}, because we can factor its output sheaf map $\oo{p}$ through
$\mathcal{S}\to\Ext{1}\out{\mathcal{X}}$, sending $(x,s)\in\mathcal{S}(n)$ to the sequence
$$\langle \rdt{f}s_0,\ldots,\rdt{f}s_n,\rdt{f}(\upd{f}(x_n,s_n)) \rangle\in\out{\mathcal{X}}(n+1).$$
Since $s_{n+1}=\upd{f}(x_n,s_n)$, this factorization is compatible with restrictions. This is due to the
fact that every current state and input decide the subsequent output.

The machine is also total and deterministic, because the square below as in \cref{eqn:tot_det_graphs} is already a pullback:
\[
\begin{tikzcd}[column sep=large]
	\Xin\times\Xin\times S\ar[r,"{(\pi_1,\pi_3)}"]\ar[d,"{(\pi_1,\pi_2)}"']\ullimit&\Xin\times S\ar[d,"\pi_1"]
	\\
	\Xin\times\Xin\ar[r,"\pi_1"']&\Xin
\end{tikzcd}
\]
This is the case since, for each current state and input elements, the next input produces a uniquely determined state
of the dynamical system, essentially because $\upd{f}$ is a function.

Having described $\beta_F$'s mapping on objects, for any map of $(\Xin,\Xout)$-discrete dynamical systems $F_1\to F_2$ described in \cref{def:DDS}, the function $h\colon S_1\to S_2$
induces a machine morphism given by $(1\times h)^{n+1}\colon\cat{S}_1(n)\to\cat{S}_2(n)$, so $\beta_X$ is a (faithful) functor.
These are natural in $X$: if $\phi$ is a wiring diagram from $X$ to $Y$, then the  following commutes (up to iso, see \cref{rem:everything_is_pseudo})
\begin{equation}\label{eq:DDS_naturality}
\begin{tikzcd}[column sep=.75in,row sep=.3in]
\DDS(\Xin,\Xout)\ar[r,"\DDS(\phi)"]\ar[d,"\beta_X"'] & \DDS(\Yin,\Yout)\ar[d,"\beta_Y"] \\
\tdADS[N](\inp{\cat{X}},\out{\cat{X}})\ar[r,"{\tdADS[N](K\phi)}"'] & \tdADS[N](\inp{\cat{Y}},\out{\cat{Y}})
\end{tikzcd}
\end{equation}
Indeed, given a wiring diagram $\phi\colon X\to Y$ as in \cref{eq:wiring_diag_Int}, the
formulas \cref{eqn:DDS_phi} for $\DDS(\phi)$ can be applied to a discrete dynamical system
$(S,\upd{f},\rdt{f})$ to give rise to a new discrete system $(S,\upd{g},\rdt{g})$. This is sent by $\beta_Y$ to the machine
\begin{displaymath}
\cat{T}(n)=\left\{(y,s)\colon\{0,\ldots,n\}\to\Yin\times S\mid s_{i+1}=\upd{f}(\fin(y_i,\rdt{f}(s_i)),s_i)\right\}
\end{displaymath}
and appropriate input and output maps.
On the other hand, one can apply the formula $\ADS_N(\phi)$ from \cref{eq:SpnS_on_arrows} to the machine
$\cat{S}\to\inp{\cat{X}}\times\out{\cat{X}}$ from \cref{eq:cat_S}, i.e.\ take
the pullback along $(\fin,\pi_2)$ and post-compose with $1\times\fout$.
The resulting machines $\cat{T}\to\inp{\cat{Y}}\times\out{\cat{Y}}$ are isomorphic.

Finally, we can verify that $\beta_{X+Y}\circ\DDS_{X,Y}=(\tdADS[N])_{\cat{X},\cat{Y}}\circ(\beta_{X}\times\beta_{Y})$ so that the natural transformation $\beta$ is monoidal; this uses $\DDS$-monoidality from \cref{sec:DandCDS} and \cref{eq:SpnS_monoidal}.
\end{proof}

Next, we describe an analogous map that transforms continuous dynamical systems
(\cref{def:CDS}) to continuous machines (\cref{def:ADS}). Because we will be comparing
different sorts of machines, we will denote continuous machines by $\ADS[C]$ rather than simply by $\ADS$.

Let $\Traj[\infty]\colon\Euc\to\Shv{\Int}$ be the functor mapping a space $A$ to the
sheaf $\Traj[\infty](A)$ of $C^\infty$-trajectories in $X$, as described in \cref{Int_sheaves_ex}(\ref{ex:gluable_smooth_functions}),
which has as length-$\ti$ sections all smooth
functions $f\colon[0,\ti]\to A$. This functor preserves all finite limits,
and it induces a strong monoidal functor $\WW{\Traj[\infty]}\colon\WW{\Euc}\to\WW{\Shv{\Int}}$
again by \cref{eq:W_}.

\begin{proposition}\label{CDS_to_ADS}
There exists a wiring diagram algebra map
\begin{displaymath}
\algmap{\Euc}{\Traj[\infty]}{\Shv{\Int}}{\CDS}{\ADS[C]}{\delta}
\end{displaymath}
\end{proposition}

\begin{proof}
We first define the components $\delta_X\colon\CDS(X)\to\ADS[C](\WW{\Traj[\infty]}(X))$
for each box $X=(\inp{X},\out{X})\in\WW{\Euc}$, mapping a continuous dynamical system on $X$ to a continuous machine
on $\WW{\Traj[\infty]}(X)$.
We then check that the $\delta_X$ satisfy the necessary naturality and monoidality conditions.
Denote by $\inp{\cat{X}}$, $\out{\cat{X}}$ the sheaves $\Traj[\infty](\Xin),\Traj[\infty](\Xout)$ associated
to the box $\WW{\Traj[\infty]}(X)$ from \cref{Int_sheaves_ex}(\ref{ex:gluable_smooth_functions}).

By \cref{def:CDS} a continuous dynamical system on $X$ is a tuple $F=(S,\dyn{f},\rdt{f})$, where $S$ is a smooth manifold,
$\dyn{f}\colon\Xin\times S\to TS$ is the dynamics, and $\rdt{f}\colon S\to\Xout$
is a smooth map. Define a sheaf $\mathcal{S}\in\Shv{\Int}$ by
\begin{displaymath}
\mathcal{S}(\ti)\coloneqq
\left\{
(x,s)\colon [0,\ti]\to\Xin\times S\;\middle|\; x,s\tn{ are smooth and }\frac{ds}{dt}=\dyn{f}(x,s)
\right\}
\end{displaymath}
where $\frac{ds}{dt}$ is the derivative of $s\colon [0,\ti]\to S$;
notice that the trajectory $s$ is a solution to the differential
equation defining the dynamical system. We can now define $\delta_X(F)$
to be the span $\inp{\mathcal{X}}\from\mathcal{S}\to\out{\mathcal{X}}$ where the maps send $(x,s)$
to $x$ and to $\rdt{f}(s)$ respectively. Again, if $\delta_X(F)=\delta_X(G)$ for two continuous dynamical systems, from the above description we can deduce that $F=G$.

We next define $\delta_X$ on morphisms, using the chain rule: consider a map of $(\Xin,\Xout)$-continuous dynamical systems $F_1\to F_2$ as defined in \cref{def:CDS}, namely a smooth function $h\colon S_1\to S_2$ such that $\rdt{f_1}\circ h=\rdt{f_2}$ and the right-hand square commutes in
\[
\begin{tikzcd}
	{[0,\ti]}\ar[r,"{(x,s)}"]\ar[dr,bend right=15pt,"\frac{ds}{dt}"']
	&\Xin\times S_1\ar[r,"1\times h"]\ar[d,"\dyn{f_1}"]&[+10pt]\Xin\times S_2\ar[d,"\dyn{f_2}"]\\
	&TS_1\ar[r,"\frac{dh}{ds}"']&TS_2
\end{tikzcd}
\]
We define a machine morphism $\delta_X(F_1)\to\delta_X(F_2)$, i.e.\ a sheaf map $\mathcal{S}_1\to\mathcal{S}_2$ over $\inp{\mathcal{X}}\times\out{\mathcal{X}}$ as follows. For any $(x,s)\in\mathcal{S}_1(\ti)$, the left-hand triangle commutes,
i.e.\ $\frac{ds}{dt}=\dyn{f_1}(x,s)$. Thus the outer shape does too, so we find that $(x,h(s))\in \mathcal{S}_2(\ti)$ by the chain rule:
\[\dyn{f_2}\big(x,h(s)\big)=\frac{dh}{ds}\ \frac{ds}{dt}=\frac{d(h\circ s)}{dt}.\]
This map commutes with the projections to $\inp{\mathcal{X}}\times\out{\mathcal{X}}$, therefore $\delta_X$ is a (faithful) functor.

We will now show naturality and monoidality of $\delta$, similarly to \cref{eq:DDS_naturality}.
Suppose $\phi\colon X\to Y$ like \cref{eq:wiring_diag_Int}.
The new machine $T\coloneqq\ADS[C](\phi)(\delta_X(F))$ is defined by
\begin{displaymath}
\cat{T}(\ti)=\left\{(y,s)\colon[0,\ti]\to\Yin\times S\;\middle|\;\frac{ds}{dt}=\dyn{f}\big(s,\fin(y,\rdt{f}(s))\big)\right\}
\end{displaymath}
by applying construction \cref{eq:SpnS_on_arrows}:
first take $\cat{S}$'s pullback along $(\fin,\pi_2)$ and the post-compose with $(1\times\fout)$.
%
This composite $\cat{T}\to\inp{\cat{Y}}\times\out{\cat{Y}}$ is equal to
$\delta_X(\CDS(\phi)(F))$ as can be seen through \cref{CDS(phi)}. It can also be verified that $\delta$ is a monoidal transformation,
by checking that
$\delta_{X+Y}\circ\CDS_{X,Y}=\ADS_{\cat{X},\cat{Y}}\circ(\delta_{X}\times\delta_{Y})$.
\end{proof}

\begin{remark}\label{rem:ODE_not_td}
Unlike \cref{DDS_to_ADS}, the algebra map $\delta$ from \cref{CDS_to_ADS} does not factor through $\tdADS[C]$,
meaning that continuous dynamical systems do not generally correspond to total and deterministic continuous machines.
For example, the dynamical system $\dot{y}=y^2$ (inhabiting the closed box $\square$ with no inputs and no outputs) is not total. Indeed, the initial value $y(0)=1$ extends to a solution $y=\frac{1}{1-t}$ that exists only for $t<1$. Similarly, the dynamical system $\dot{y}=2\sqrt{|y|}$ is not deterministic. Indeed, the initial value $y(0)=0$ has the following solution for any $a$:
\[
\begin{cases}
	y=0 		&\tn{if } t\leq a\\
	y=(t-a)^2	&\tn{if } t\geq a
\end{cases}
\]
\end{remark}

\subsection{Maps between machines}\label{sec:maps_between}

We begin by providing maps between machines over the toposes of
continuous, discrete, and synchronous sheaves discussed in \cref{ch:ADS}.
These algebra maps naturally group together according to whether they refer to
general machines, or to total or deterministic variations; in fact
they are morphisms of $\WDalg$ in the more general setting described in \cref{sec:Span_Systems}.

As mentioned earlier, for $\cat{C}=\Shv{\Int}$, $\Shv{\Int}_N$, and $\Shv{\Int}_\Sync$,
the construction $\SpnS_{\cat{C}}$ of \cref{prop:C_Span_Algebra}
gives the wiring diagram algebras $\ADS[C]$, ${\ADS[N]}$, and ${\ADS[\Sync]}$ described at \cref{prop:ADS,prop:ADS_N,prop:ADS_Sync}.
Moreover, since the functors $\Sigma_i'\colon\Shv{\Int}_N\to\Shv{\Int}/\Sync$ from \cref{eq:Sigma_i'}
and $\Delta_{\Sync!}\colon\Shv{\Int}\to\Shv{\Int}/\Sync$ from \cref{slice_adjunction}
between the respective typing categories
preserve finite limits, they directly induce monoidal natural transformations $\SpnS_{\Sigma'_i}$
and $\SpnS_{\Delta_{\Sync!}}$ by \cref{prop:Mch_FCCat}. 

\begin{corollary}\label{cor:1}
There exist algebra maps realizing any continuous or discrete machine as a synchronous one:
\begin{displaymath}
\algmap{\Shv{\Int}_N}{\Sigma'_i}{\Shv{\Int}/_\Sync}{\ADS[N]}{\ADS[\Sync]}{\SpnS_{\Sigma'_i}}\qquad
\algmap{\Shv{\Int}}{\Delta_{\Sync!}}{\Shv{\Int}/_\Sync}{\ADS[C]}{\ADS[\Sync]}{\SpnS_{\Delta_{\Sync!}}}
\end{displaymath}
\end{corollary}

Essentially, the first algebra morphism maps a discrete machine $p\colon S\to\Xin\times\Xout$ to $\Sigma'_ip$,
a synchronous machine whose state, input and output sheaves are obtained from the old ones by applying $\Sigma_i((\Inp)\To{!}\1)$,
see \cref{prop:Sigma_formula}. Similarly, the second one gives the mapped span under $\Delta_{\Sync!}$,
where $\Delta_{\Sync!}(X)=(X\times\Sync\to\Sync)$ for the synchronizing sheaf \cref{sync(l)}. Notice that since $\Sigma_i'$ and $\Delta_{\Sync!}$ are both faithful functors (see proof of \cref{prop:synchronized_geom_morph}), their respective algebra maps have embeddings as components
as explained after \cref{prop:Mch_FCCat}.

This significant result accomplishes one of the main goals of this work, by bringing discrete and continuous systems in
a common environment, that of synchronous systems. For example, both discrete and continuous dynamical systems,
which can be realized as discrete and continuous machines respectively by \cref{DDS_to_ADS,CDS_to_ADS},
can now be further translated under the above algebra maps to synchronous machines. As a result, their arbitrary interconnections
can be studied in this common framework.

\begin{remark}
Using \cref{prop:Mch_FCCat} we can also recover the fact from \cite{Spivak:2015b} that extracting the steady states
of a dynamical system, and organizing them in terms of matrices, amounts to an algebra homomorphism. Indeed,
this follows from the obvious fact that for any finitely complete category $\cat{C}$ and object $c\in\cat{C}$, the $\Hom$-functor
$\cat{C}(c,-)\colon\cat{C}\to\Set$ is finitely complete, so we get a map of wiring diagram algebras
\[
\SpnS_{\cat{C}(c,-)}\colon\SpnS_{\cat{C}}\to\SpnS_{\Set}.
\]
For example if $\cat{C}=\Shv{\Int}$ and $c=\{*\}$ is the terminal sheaf (constant on one generator)
then $\cat{C}(\{*\},-)$ extracts the set of constant sections from any sheaf, and $\SpnS_{\cat{C}(\{*\},-)}$ is the steady state
extraction from \cite{Spivak:2015b}. If instead $c=\Sync$ is the synchronizing sheaf, we extract periodic cycles of length 1.
\end{remark}

Turning to maps preserving totality and determinism, recall that total and deterministic machines require a notion of extension:
for any extension of input, there exists (or is at most one) an extension of state to match.
We now abstract all the necessary structure to make sense of this at a higher level of generality.
Doing so will allow us to express
conditions of inertiality and totality/determinism from \cref{prop:total_det_via_extensions} and \cref{def:inertial_map},
for general span algebras like \cref{prop:C_Span_Algebra}.
We begin with a preliminary definition.

\begin{definition}
We define a \emph{Euclidean poset} to be a poset $(E,\leq)$ equipped with a symmetric monoidal structure $(+,0)$ such that
\begin{compactitem}
  \item the unit is minimal: $0\leq e$ for all $e\in E$
  \item for all $a,b\in E$, if $b\neq 0$ there exists $N\in\NN$ and $r\in E$ such that 
    \[
      a=N\cdot b+r
      \quad\tn{and}\quad
      r\leq b
    \]  
  where $N\cdot b\coloneqq b+\cdots+b$ ($N$ times).
\end{compactitem}
A \emph{morphism} of Euclidean posets is a strong monoidal functor.
\end{definition}

\begin{lemma}
Suppose $\nu\colon E\to E'$ is a morphism of Euclidean posets. If there exists $b\neq 0$ such that $\nu(b)=0$ then $\nu(a)=0$ for all $a$; such a morphism is called \emph{trivial}, otherwise it is \emph{nontrivial}.
\end{lemma}

To cast totality at this level of generality, we need our categories and functors to be regular,
rather than just finitely-complete; i.e. we restrict our general $\cat{C}$-span system functor
$\SpnS_{(\Inp)}\colon\FCCat\to\WDalg$ \cref{eq:spantriangle} to $\RegCat$,
the 2-category of regular categories and functors, with all natural transformations between them.
Recall that every topos is a regular category and every epimorphism in a topos is regular.
See \cite{barr:1971a} or \cite{Johnstone:2002a} for more on regular categories and functors.

For any regular category $\cat{C}$, $(\End(\cat{C}),\circ,1_\cat{C})$ is the monoidal category of regular endofunctors on $\cat{C}$. The following generalize the respective constructions from \cref{sec:total_det_inert}.

\begin{definition}\label{def:general_extensions}
For $\cat{C}\in\RegCat$, an \emph{extension structure} on $\cat{C}$
is a pair $(E,\Ext{})$, where $(E,\leq,+,0)$ is a Euclidean poset and $\Ext{}\colon E\op\to\End(\cat{C})$ is a strong monoidal functor.
\end{definition}

Spelling out what it means for $\Ext{}\colon E\op\to\End(\cat{C})$ to be strong monoidal,
\begin{compactitem}
  \item for all $e\in E$, the functor $\Ext{e}\colon\cat{C}\to\cat{C}$ preserves finite limits and regular epimorphisms;
  \item $\Ext{0}=\id_{\cat{C}}$, and for all $a,b\in E$, we have $\Ext{a+b}=\Ext{a}\circ\Ext{b}$;
  \item for any $e'\leq e$ there is a chosen map which we will denote $\lambda_{e,e'}\colon\Ext{e}\Rightarrow\Ext{e'}$. When $e'=0$ we abbreviate $\lambda_{e,0}$ by $\lambda_e\colon\Ext{e}\Rightarrow\id_{\cat{C}}$; and
  \item for all $a,b\in E$, there exist $N\in\NN$ and $r\in E$ such that $\Ext{a}=\Ext{N\cdot b}\circ\Ext{r}$.
\end{compactitem}

Using these extension structures we can formalize inertiality for the $\cat{C}$-span systems of \cref{def:CSpanSystems}, generalizing \cref{def:inertial_map}. Notice that we denote the components ${(\lambda_e)}_A\colon\Ext{e}A\to A$ of the natural transformation $\lambda_e$ by the same name, in order to avoid double subscripts; the (co)domain should clarify each component.

\begin{definition}\label{def:general_inertiality}
Given an extension structure on $\cat{C}$, we say that a map $p\colon S\to A$ in $\cat{C}$ is \emph{$\Ext{}$-inertial} if there
exists $0\neq e\in E$ such that $p$ factors through $\lambda_e\colon\Ext{e}(A)\to\Ext{0}(A)=A$.
For any box $X\in\WW{\cat{C}}$, we say that a $\cat{C}$-span $(\ii{p},\oo{p})\colon S\to\Xin\times\Xout$
is \emph{$\Ext{}$-inertial} if $\oo{p}$ is.
Denote $\SpnS^{\mathrm{in}}_\cat{C}(X)\ss\SpnS_\cat{C}(X)$ the full subcategory of all $\Ext{}$-inertial $\cat{C}$-span systems on $X$. 
\end{definition}

We can also generalize totality and determinism, using the following lemma (compare to \cref{prop:total_det_via_extensions}).
\begin{lemma}\label{lemma:small_extensions_big}
Let $\cat{C}$ be a regular category with extension structure $\Ext{}\colon E\op\to\End(\cat{C})$. Given $e\in E$ and $p\colon S\to A$,
define $h_e\colon\Ext{e}S\to S_e$ to be the universal map to the pullback denoted $S_e$ in
\begin{equation}\label{eqn:gen_total_det}
\begin{tikzcd}[sep=25pt]
\Ext{e}S \arrow[drr, bend left=25, "\lambda_e"] \arrow[ddr, bend right=30, "\Ext{e}p"']
\arrow[dr, "h_e"]&[-25pt]\\[-15pt]
& S_e \ullimit\arrow[r,"\lambda_e'"]\arrow[d,"p_e"']
& S\arrow[d, "p"] \\
& \Ext{e}A \arrow[r,"\lambda_e"'] & A
\end{tikzcd}
\end{equation}
Then the following are equivalent:
\begin{compactenum}
  \item\label[item]{item:biggy} for all $e\in E$, the map $h_e\colon\Ext{e}S\to S_e$ is a regular epimorphism (resp.\ a monomorphism, an isomorphism);
  \item\label[item]{item:smally} there exists $0\neq d\in E$ such that for all $d'\leq d$, the map $h_{d'}\colon\Ext{d'}S\to S_e$ is a regular epimorphism (resp.\ a monomorphism, an isomorphism).
\end{compactenum}
\end{lemma}
\begin{proof}
Clearly \cref{item:biggy}\implies\cref{item:smally}, so assume the latter holds for $d\neq0$. By the assumption that $E$ 
is a Euclidean poset, it suffices to show that for any $a,b\in E$ and $p\colon S\to A$, if both $h_{a}$ and $h_{b}$ have property
$\mathbf{P}$, where $\mathbf{P}$ is the property of being a regular epimorphism (resp.\ a monomorphism, an isomorphism) then so does $h_{a+b}$.

Consider the two diagrams below, where the left-hand diagram is as in \cref{eqn:gen_total_det}:
\[
  \begin{tikzcd}[baseline=(bot)]
    \Ext{a}S\ar[d,"h_a"']\\
    S_a\ar[r,"\lambda'_a"]\ar[d,"p_a"']\ullimit&S\ar[d,"p"]\\
    |[alias=bot]|\Ext{a}A\ar[r,"\lambda_a"']&A
  \end{tikzcd}
\hspace{.5in}
  \begin{tikzcd}[column sep=large,baseline=(bot)]
    \Ext{a+b}S\ar[d,"{\Ext{b}h_{a}}"']\\
    \Ext{b}S_a\ar[r,"{\Ext{b}\lambda_a}"]\ar[d,dashed,"h"']\ar[dd,bend right=40,end anchor=north west,"{\Ext{b}p_a}"']
    \ullimit&\Ext{b}S\ar[d,"h_b"']\\
    S_{a+b}\ar[r,"\lambda'_a"]\ar[d,"p_{a+b}"']\ullimit&S_b\ar[d,"p_b"']\ar[r,"\lambda_b'"]\ullimit&S\ar[d,"p"]\\
    |[alias=bot]|\Ext{a+b}A\ar[r,"{(\lambda_a)_{\Ext{b}A}}"']&\Ext{b}A\ar[r,"{(\lambda_b)_A}"']&A
  \end{tikzcd}
\]
For the right-hand diagram, begin by forming the two bottom pullback squares. By definition of $h_b$ \cref{eqn:gen_total_det}, we have $\Ext{b}p=p_b\circ h_b$, so the vertical rectangle is an application of $\Ext{b}$ to the left-hand diagram, and is thus a pullback because $\Ext{b}$ is regular. There is an induced map $h\colon\Ext{b}S_a\to S_{a+b}$, and we now have that each indicated square in the diagram is a pullback. 

Since $\Ext{a+b}p=\Ext{b}\Ext{a}p=\Ext{b}p_a\circ\Ext{b}h_a$ is the long vertical map, the universal property of pullback implies that $h\circ\Ext{b}h_a=h_{a+b}$. If $h_a$ and $h_b$ have property $\mathbf{P}$ then so does $h_{a+b}$ because its factors $h$ and $\Ext{b}h_a$ do: $\mathbf{P}$ is stable under pullbacks and $\Ext{b}$ preserves $\mathbf{P}$.
\end{proof}

\begin{definition}\label{def:general_tot_det}
Let $\cat{C}$ be a regular category with extension structure $\Ext{}\colon E\op\to\End(\cat{C})$,
and let $p\colon S\to B$ be a morphism. We say that $p$ is \emph{$\Ext{}$-total} (resp.\ \emph{$\Ext{}$-deterministic},
\emph{$\Ext{}$-total-determinstic}) if the map $h_{e}$ in \cref{eqn:gen_total_det} is a regular epimorphism
(resp.\ a monomorphism, an isomorphism) for all $e\in E$.

For any box $X\in\WW{\cat{C}}$, we call an $\Ext{}$-inertial span system $(\ii{p},\oo{p})\in\SpnS^{\mathrm{in}}_\cat{C}(X)$
(as in \cref{def:general_inertiality}) \emph{$\Ext{}$-total} (resp.\ \emph{$\Ext{}$-deterministic}, \emph{$\Ext{}$-total-deterministic})
if $\ii{p}$ is. Denote $\SpnS^{\mathrm{t}}_\cat{C}(X),\SpnS^{\mathrm{d}}_\cat{C}(X),\SpnS^{\mathrm{td}}_\cat{C}(X)
\ss\SpnS^{\mathrm{in}}_\cat{C}(X)$ the full subcategories of all $\Ext{}$-total (resp. $\Ext{}$-deterministic,
$\Ext{}$-total-deterministic) span systems on $X$.
\end{definition}

\begin{example}
For any $\epsilon\geq0$ the functor $\Ext{\epsilon}\colon\Shv{\Int}\to\Shv{\Int}$, as in \cref{Ext}, is a
endomorphism of $\Shv{\Int}$; we showed it is finitely complete in \cref{lemma:ext_pbs}, and it is not hard to check that it preserves epimorphisms (all of which are regular). The poset
$E=(\RR_{>0},\geq)$ of positive reals is a Euclidean poset, and the left restriction maps
$\Ext{\epsilon}\to\Ext{\epsilon'}$ for any $\epsilon\geq\epsilon'$ constitute an extension
structure $\Ext{}\colon E\op\to\End(\Shv{\Int})$. This extension structure can be lifted to one on $\Shv{\Int}/\Sync$. Namely, for every
$e\in\RR_{\geq0}$, we define $\Ext{e}(X\to\Sync)$ to be the composite $\Ext{e}X\to\Ext{e}(\Sync)\To{\lambda_e}\Sync$.
There is also an extension structure on $\Shv{\Int}_N$, such that \cref{def:general_extensions} generalizes the
definitions in \cref{sec:discrete_variations}, by taking $E=\NN$. 

In all the above cases, the notions of inertial, total, deterministic, and total-deterministic morphisms and span systems
generalize the older notions, as aimed. 
\end{example}

\begin{proposition}\label{prop:tot_det_variations1}
Let $\cat{C}$ be a regular category with an extension structure $\Ext{}\colon E\op\to\End(\cat{C})$.
Then there are symmetric lax monoidal functors
\[
\SpnS^{\mathrm{in}}_\cat{C},\;\SpnS^{\mathrm{t}}_\cat{C},\;\SpnS^{\mathrm{d}}_\cat{C},\;\SpnS^{\mathrm{td}}_\cat{C}
\colon\WW{\cat{C}}\longrightarrow\Cat
\]
defined as in \cref{def:general_inertiality,def:general_tot_det} on any $X\in\WW{\cat{C}}$,
which are subalgebras of $\SpnS_\cat{C}$ \cref{eq:defSpnS}.
\end{proposition}

\begin{proof}
The proof for inertiality closely follows that of \cref{prop:iADS}, and
totality and determinism follow that of \cref{prop:tADS}.
\end{proof}

These subfunctors constitute, on their own right, the mapping of a more general functor
on objects, namely regular categories. Towards that end, we define the following category.
\begin{definition}\label{def:mor_regcate}
Let $(\cat{C},E,\Ext{})$ and $(\cat{C'},E',\Ext{}')$ be regular categories with extension structures. A \emph{morphism}
between them is a pair $(F,\nu)$ where $F\colon\cat{C}\to\cat{C'}$ is a regular functor, $\nu\colon E\to E'$ is a nontrivial morphism of
Euclidean posets, and such that for all $e\in E$ the following diagram commutes (up to natural isomorphism)
\[
  \begin{tikzcd}
    \cat{C}\ar[r,"{\Ext{e}}"]\ar[d,"F"']&\cat{C}\ar[d,"F"]\\
    \cat{C'}\ar[r,"{\Ext{\nu(e)}'}"']&\cat{C'}
  \end{tikzcd}
\]
These form the category of \emph{regular categories with extensions}, denoted $\RegCate$, which naturally maps to
$\RegCate\to\RegCat\ss\FCCat$.
\end{definition}

\begin{theorem}\label{thm:tot_det_variations}
The restricted functor $\SpnS_{(\Inp)}\colon\RegCate\to\RegCat\hookrightarrow\FCCat\to\WDalg$ as in \cref{prop:Mch_FCCat}
has subfunctors
\begin{displaymath}
\SpnS^{\mathrm{in}}_{(\Inp)},\SpnS^{\mathrm{t}}_{(\Inp)},
\SpnS^{\mathrm{d}}_{(\Inp)},\SpnS^{\mathrm{td}}_{(\Inp)}\colon\RegCate\to\WDalg 
\end{displaymath}
which map any object $\Ext{}\colon E\op\to\End(\cat{C})$ in $\RegCate$ to the respective algebras of
\cref{prop:tot_det_variations1}.
\end{theorem}

\begin{proof}
For their mapping on morphisms, given some $(F,\nu)\colon (\cat{C},E,\Ext{})\to(\cat{C'},E',\Ext{}')$
as in \cref{def:mor_regcate}, we need to show that $\SpnS_F$ of \cref{eq:defSpnSF} appropriately restricts,
through its component functors, to a 2-cell of the form
\[
\algmap{\cat{C}}{F}{\cat{C'}}{\SpnS^\mathbf{x}_\cat{C}}{\SpnS^\mathbf{x}_\cat{C'}}{\SpnS^\mathbf{x}_F}
\]
where $\mathbf{x}\in\{\nfun{i},\;\nfun{t},\;\nfun{d},\;\nfun{td}\}$
stands for are all respective subalgebras of $\SpnS_\cat{C}$.
Therefore it suffices to show that in all three cases, a dashed arrow exists
for any $X\in\WW{\cat{C}}$, making it commute in $\Cat$:
\begin{equation}\label{eq:commuteswithinclusions}
  \begin{tikzcd}[column sep=.7in]
    \SpnS^\mathbf{x}_\cat{C}(X)\ar[d,hook]\ar[r,dashed]&\SpnS^\mathbf{x}_\cat{C'}(\cat{W}_FX)\ar[d,hook]\\
    \SpnS_\cat{C}(X)\ar[r,"(\SpnS_F)_X"']&\SpnS_\cat{C'}(\cat{W}_FX)
  \end{tikzcd}
\end{equation}
For inertial machines, the existence of a dashed arrow is implied by the fact that $\nu$ is nontrivial and $F\Ext{e}\cong\Ext{\nu(e)}F$
for all $e\in E$. For total (resp.\ deterministic) machines, the existence of a dashed arrow is implied by the
fact that $\Ext{e}$ preserves regular epimorphisms (resp.\ monomorphisms) for all $e$. Moreover, notice that if $(\SpnS_F)_X$ is an embedding, then clearly so is $(\SpnS^\mathbf{x}_F)_X$.
\end{proof}

We can now apply \cref{thm:tot_det_variations} to the regular categories $\Shv{\Int},\Shv{\Int}_N$ and $\Shv{\Int}/\Sync$ from
\cref{sec:Int_sheaves,sec:synchronization}, and to the functors $\Sigma_i'\colon\Shv{\Int}_N\to\Shv{\Int}/\Sync$
\cref{eq:Sigma_i'} and $\Delta_{\Sync!}\colon\Shv{\Int}\to\Shv{\Int}/\Sync$ \cref{slice_adjunction}, in order to obtain algebra maps
between the total and deterministic variations of continuous, discrete and synchronous machines,
as anticipated. Recall that as defined in \cref{ch:ADS}, the category of e.g. continuous machines is precisely
the category of $\Shv{\Int}$-span systems $\SpnS_{\Shv{\Int}}$, denoted $\ADS[C]$; similarly for the rest of terminology.

\begin{corollary}\label{cor:2}
Let
$\mathbf{x}\in\{\nfun{i},\;\nfun{t},\;\nfun{d},\;\nfun{td}\}$ stand for inertial, total, deterministic, or total-deterministic. There are 
algebra embeddings
\begin{displaymath}
\algmap{\Shv{\Int}_N}{\Sigma_i'}{\Shv{\Int}/_\Sync}{\xADS[N]}{\xADS[\Sync]}{\SpnS^\mathbf{x}_{\Sigma_i'}}\qquad
\algmap{\Shv{\Int}}{\Delta_{\Sync!}}{\Shv{\Int}/_\Sync}{\xADS[C]}{\xADS[\Sync]}{\SpnS^\mathbf{x}_{\Delta_{\Sync!}}}
\end{displaymath}
which translate the specific classes of discrete machines into synchronous machines of the same class,
and similarly for continuous into synchronous.
\end{corollary}

\begin{proof}
First of all, notice that both $\Sigma_i'\colon\Shv{\Int}_N\to\Shv{\Int}/\Sync$ and $\Delta_{\Sync!}\colon\Shv{\Int}\to\Shv{\Int}/\Sync$
are inverse image parts of geometric morphisms, so they are regular functors.
Clearly $(\Delta_{\Sync!},\id_{\RR_{\geq0}})$ is a map in $\RegCate$, so the second algebra maps follow
directly from \cref{thm:tot_det_variations}.

Now the inclusion $\nu\colon\NN\to\RR_{\geq 0}$ is a nontrivial morphism of Euclidean posets,
so to complete the proof, we need to show that $\Sigma'_i\circ\Ext{n}\cong\Ext{\nu(n)}\circ\Sigma'_i$ for any $n\in\NN$. It suffices to show
this for $n=1$, whence $\nu(n)=1$. 
Recall from \cref{prop:Sigma_formula} the formula
\[\Sigma_iX(\ti)\cong\bigsqcup_{r\in[0,1)}{X\big(\ceil{r+\ti}\big)}\]
where the extra data of $\Sigma_i'X$ is just the map $\Sigma_X\to\Sync$ given by $r$. Thus the result follows from the equation $\ceil{\ti+r+1}=\ceil{\ti+r}+1$.
\end{proof}

Notice that by construction of the above algebra morphisms, their components commute with the embeddings of each subclass of machines into the general discrete or continuous ones, as in \cref{eq:commuteswithinclusions}.

\begin{remark}
It should be noted that the only reason we need to work with regular categories and functors rather than
finitely complete ones, is for totalness. We can
replace $\RegCat$ by $\FCCat$ and \cref{thm:tot_det_variations} will still hold for inertial,
deterministic, and total-deterministic machines. That is, in \cref{def:general_extensions}, we could
ask only that $\cat{C}$ be finitely complete, and that endofunctors
$(\cat{C}\to\cat{C})\in\End(\cat{C})$ preserve finite limits. Similarly we can drop condition that
our finitely-complete functors $F\colon\cat{C}\to\cat{C'}$ are regular;
going through with appropriate changes to all constructions above would only exclude totalness results.
\end{remark}

The above \cref{cor:2} successfully restricts the more general \cref{cor:1} to the total and deterministic
variations of continuous, discrete and synchronous machines. Once again, (total, deterministic) synchronous
machines end up being the common framework where their discrete and continuous counterparts can be
studied together. The final proposition below ensures that contracted machines
described in \cref{safety_contracts} fit in the same picture. 

\begin{proposition}
There exists a functor $\Cntr_{(\Inp)}\colon\RegCat\to\WDalg$ making the diagram
\[
  \begin{tikzcd}
    &\WDalg\ar[d,"U"]\\
    \RegCat\ar[ur,"{\Cntr_{(\Inp)}}"]\ar[r,"{\WW{(\Inp)}}"']&\SMonCats
  \end{tikzcd}
\]
commute. Moreover, there exists a natural transformation
\begin{displaymath}
\begin{tikzcd}[row sep=.1in]
& \FCCat\ar[dr,"\SpnS_{(\Inp)}"] & \\
\RegCat\ar[ur,hook]\ar[rr,bend right=15,"\Cntr_{(\Inp)}"']\ar[rr,phantom,"{\scriptstyle\Downarrow\Img}"description] && \WDalg
\end{tikzcd}
\end{displaymath}
whose components translate each machine (therefore any total/deterministic subclass, \cref{thm:tot_det_variations})
into its 'maximal' validated safety contract.
\end{proposition}

\begin{proof}
Similarly to \cref{prop:Mch_FCCat}, this functor maps any regular category $\cat{C}$ to its algebra $\Cntr_\cat{C}\colon\cat{W}_\cat{C}\to\Cat$
of safety contracts essentially described in \cref{prop:Cntr}, i.e. $\Cntr_{\cat{C}}(X)\coloneqq\Sub_{\cat{C}}(\Xin\times\Xout)$
for any box $X\in\WW{\cat{C}}$. For any $\phi\colon X\to Y$ in $\cat{C}$, the functor $\Cntr_{\cat{C}}(\phi)$ is given by
the construction \cref{eq:Cntr_on_arrows} (seen inside an arbitrary regular category), and the symmetric lax monoidal structure
follows since products of inclusions are inclusions. Now any regular functor $F\colon\cat{C}\to\cat{D}$
preserves epi-mono factorization and pullbacks, so it induces a map 
\begin{displaymath}
\algmap{\cat{C}}{F}{\cat{D}}{\Cntr_{\cat{C}}}{\Cntr_{\cat{D}}}{\Cntr_F}
\end{displaymath}
with components functors $(\Cntr_F)_X\colon\Cntr_{\cat{C}}(X)\to\Cntr_{\cat{D}}(\cat{W}_FX)$
for any $X\in\WW{C}$ being just application of $F$ on the respective subobjects.

Now the natural transformation $\Img$ has components wiring diagram algebra maps $\Img_\cat{C}\colon\SpnS_\cat{C}\to\Cntr_\cat{C}$,
formed by the mappings
$$\SpnS_\cat{C}(X)=\cat{C}/(\Xin\times\Xout)\to\Sub_{\cat{C}}(\Xin\times\Xout)=\Cntr_\cat{C}(X)$$
which take the image of each $S\to\Xin\times\Xout$; these are functorial and symmetric lax monoidal. Finally, naturality of $\Img$
\[
  \begin{tikzcd}
    \SpnS_\cat{C}\ar[r,"{\SpnS_F}"]\ar[d,"\Img_{\cat{C}}"']&\SpnS_\cat{D}\ar[d,"\Img_{\cat{D}}"]\\
    \Cntr_{\cat{C}}\ar[r,"\Cntr_{F}"']&\Cntr_{\cat{D}}
  \end{tikzcd}
\]
is verified when we write the above commutativity inside $\WDalg$, i.e. arrows as in \cref{eq:morphWDalg}:
\[
\begin{tikzcd}[row sep=.3in, column sep=.6in]
\WW{\cat{C}}\ar[d,"\WW{F}"']\ar[drr, "\SpnS_\cat{C}", bend left]
\ar[drr, phantom,"\Downarrow{\scriptstyle\SpnS_F}"]&&&
\WW{\cat{C}}\ar[drr,"\SpnS_\cat{C}", bend left]\ar[d,equal]\ar[drr, phantom, "\Downarrow{\scriptstyle\Img_\cat{C}}"description]
\\
\WW{\cat{D}}\ar[d,equal]\ar[rr,"\SpnS_\cat{D}"description]\ar[drr, phantom,"\Downarrow{\scriptstyle\Img_\cat{D}}"]
&& \Cat\ar[r, phantom,"="description]
&
\WW{\cat{C}}\ar[rr,"\Cntr_\cat{C}"description]\ar[d,"\WW{F}"']\ar[drr, phantom, "\Downarrow{\scriptstyle\Cntr_F}"description]
&& \Cat
\\
\WW{\cat{D}}\ar[urr,bend right,"\Cntr_\cat{D}"'] && \phantom{A} &
\WW{\cat{D}}\ar[urr,bend right,"\Cntr_\cat{D}"'] && \phantom{A}
\end{tikzcd}
\]
\end{proof}
This last algebra map has the significant effect of translating machines of any kind to a safety contract consisting of all its `valid' behaviors,
i.e. lists of all possible inputs and outputs through time; see also \cref{def:safety_contract}.
As a result, whenever we have an interconnection \cref{picture_wiring_diagram} of arbitrary systems, we could 
directly reason about the valid behaviors of the composite system without composing the machines first. 

\appendix
\section{Discrete \texorpdfstring{\Conduche}{Conduche} fibrations}\label{appendix:DSF}

In this appendix, we discuss an equivalent view of interval sheaves
from \cref{sec:Int_sheaves} in terms of \emph{discrete \Conduche fibrations}, elsewhere \cite{Bunge.Fiore:2000a} called \emph{unique factorization lifting functors}. This view largely follows material found in \cite{Johnstone:1999a},
and allows us to think of continuous or discrete interval sheaves as categories equipped with `length' functors,
called \emph{durations} in \cite{Lawvere:1986a}, into the monoids $(\RRnn,0,+)$ and $(\NN,0,+)$ viewed as single-object categories.

\subsection{Discrete fibrations, opfibrations, and \Conduche fibrations}
	\label{sec:various_fibrations}
For any category $\cat{C}$, consider the diagram of sets and functions
\[
\begin{tikzcd}
	C_2\ar[r,shift left=5pt]\ar[r,"\circ"description]\ar[r,shift right=5pt]&C_1\ar[r,shift left=2pt,"s"]\ar[r,shift right=2pt,"t"']&C_0
\end{tikzcd}
\]
where $C_0$ is the set of objects, $C_1$ is the set of morphisms, and $C_2$ is the set of composable morphisms in $C$. The two functions $C_1\to C_0$ send a morphism $f$ to its source (domain) and its
target (codomain); the three functions $C_2\to C_1$ send a composable pair $(g,f)$ to $g$, $g\circ f$, and $f$. We have left out of our
diagram the $(i+1)$ functions $C_i\to C_{i+1}$ induced by identity morphisms.
For any functor $F\colon\cat{C}\to\cat{D}$, there is a function $F_i\colon C_i\to D_i$ for each $i\in\{0,1,2\}$, and each is induced by the
function $F_0$ on objects and the function $F_1$ on morphisms.

\begin{definition}\label{def:dopf_df_dcf}
For a functor $F\colon\cat{C}\to\cat{D}$, consider the commutative diagrams
\[
\begin{tikzcd}
	C_1\ar[r,"t"]\ar[d,"F_1"']&C_0\ar[d,"F_0"]\\
	D_1\ar[r,"t"']&D_0
\end{tikzcd}
\qquad
\begin{tikzcd}
	C_1\ar[r,"s"]\ar[d,"F_1"']&C_0\ar[d,"F_0"]\\
	D_1\ar[r,"s"']&D_0
\end{tikzcd}
\qquad
\begin{tikzcd}
	C_2\ar[r,"\circ"]\ar[d,"F_2"']&C_1\ar[d,"F_1"]\\
	D_2\ar[r,"\circ"']&D_1
\end{tikzcd}
\]
Then $F$ is called a \emph{discrete fibration} (resp.\ \emph{a discrete opfibration, a discrete \Conduche
fibration}) if the first (resp.\ second, third) square is a pullback.
\end{definition}

The first two conditions clearly correspond to the well-known definitions of discrete (op)fibrations.
For example, the first one says that for every morphism $h\colon x\to y$ in the base category $\cat{D}$ and every
object $d\in\cat{C}$ above $y$, there exists a unique morphism $c\to d$ which maps to $h$ via $F$. 
The third one, on the other hand, says that given a morphism $f\colon c\to d$
in the domain category $\cat{C}$ such that $Ff=v\circ u$ factorizes in the base category $\cat{D}$, there exists a unique
factorization $f=h\circ g$ with $Fh=v$ and $Fg=u$. This is also equivalent, \cite{Lawvere:1986a},
to the isomorphism of the factorization categories $\Fact(f)\cong\Fact(Ff)$ from \cref{def:factorization_linear}.

\begin{lemma}
If $F$ is a discrete opfibration (resp.\ a discrete fibration), then it is a discrete \Conduche fibration.
\end{lemma}
\begin{proof}
The top, bottom, and front of the following cube are pullbacks, so the back is too:
\[
\begin{tikzcd}[sep=small]
	C_2\ar[rr]\ar[rd]\ar[dd]&&C_1\ar[dd]\ar[rd]\\
	&C_1\ar[rr,crossing over]&&C_0\ar[dd]\\
	D_2\ar[rr]\ar[rd]&&D_1\ar[rd]\\
	&D_1\ar[from=uu,crossing over]\ar[rr]&&D_0
\end{tikzcd}
\]\end{proof}

\begin{lemma}\label{lem:DCF_composite}
If $F\colon\cat{C}\to\cat{D}$ and $G\colon\cat{D}\to\cat{E}$ are functors and $G$ is a discrete \Conduche fibration,
then $G\circ F$ is discrete \Conduche if and only if $F$ is.
The same holds if we replace discrete \Conduche with discrete (op)fibrations.
\end{lemma}
\begin{proof}
This is just the pasting lemma for pullback squares.\end{proof}

Small discrete \Conduche fibrations form a wide subcategory of the category of small categories, $\DCF\ss\Cat$ .
In particular, $\DCF/\cat{A}$ denotes the slice category of discrete
\Conduche fibrations over any $\cat{A}\in\Cat$. Our main case of interest is the case $\cat{A}=\cat{R}$, the additive monoid of reals,
which is also the primary example in the development of \cite{Lawvere:1986a}.
Hence a discrete \Conduche fibration $\len\colon\cat{C}\to\cat{R}$ for any category $\cat{C}$ 
amounts to a commutative diagram as to the left
\begin{equation}\label{length_functor}
\begin{tikzcd}[sep=large]
	C_2\ar[d,"\len\times\len"']\ar[r,shift left=6pt,"\pi_1"]\ar[r,"\circ" description]\ar[r,shift right=6pt,"\pi_2"']&C_1
	\ar[d,"\len"]\ar[r,shift left=2pt,"s"]\ar[r,shift right=2pt,"t"']&C_0\ar[d,"!"]\\
	\RRnn\times\RRnn\ar[r,shift left=6pt,"\pi_1"]\ar[r,"+" description]\ar[r,shift right=6pt,"\pi_2"']
	&\RRnn\ar[r,shift left=2pt,"!"]\ar[r,shift right=2pt,"!"']&\singleton
\end{tikzcd}
\qquad\qquad
\begin{tikzcd}[sep=large]
	C_2\ar[r,"\circ"]\ar[d,"\len\times\len"']\ar[dr,phantom,very near start, "\lrcorner" description]&C_1\ar[d,"\len"]\\
	\RRnn\times\RRnn\ar[r,"+"']&\RRnn
\end{tikzcd}
\end{equation}
for which the sub-diagram extracted to its right is a pullback in $\Set$. In other words, every morphism $f$ in $\cat{C}$ has a length
$\len(f)\in\RRnn$ and, for any way to write $\len(f)=\ti_1+\ti_2$ as a sum of nonnegative numbers,
there is a unique pair of composable morphisms $f=f_1\circ f_2$ in $\cat{C}$ having those lengths $\len(f_1)=\ti_1$
and $\len(f_2)=\ti_2$. 

\begin{proposition}\label{prop:mapping_DCF}
For any small category $\cat{C}$, the slice category $\DCF/\cat{C}$ is reflective in $\Cat/\cat{C}$. Moreover,
if $F\colon\cat{C}\to\cat{\cat{D}}$ is any functor, there is a diagram
\[
\begin{tikzcd}[sep=.6in]
\DCF/\cat{C}\ar[d,shift left=1.7,"U_\cat{C}","{\bbot}"']\ar[r,shift left=1.7,"\Sigma_F","{\bot}"']
&
\DCF/\cat{D}\ar[d,shift left=1.7,"U_\cat{D}","{\bbot}"']\ar[l,shift left=1.7,"\Delta_F"]
\\
\Cat/\cat{C}\ar[u,shift left=1.7,"L_\cat{C}"]\ar[r,shift left=1.7,"F\scirc(-)","{\bot}"']
&
\Cat/\cat{D}\ar[u,shift left=1.7,"L_\cat{D}"]\ar[l,shift left=1.7,"F^*"]
\end{tikzcd}
\]
which commutes for the right (and hence left) adjoints, where $F^*$ is given by pullback along $F$.
\end{proposition}
\begin{proof}
The existence of a left adjoint $L_\cat{C}$ to the inclusion $U_\cat{C}\colon\DCF/\cat{C}\to\Cat/\cat{C}$ is proven in
\cite[Prop.~1.3]{Johnstone:1999a},
by showing that $U_\cat{C}$ preserves all limits and satisfies a solution set condition.

If $\cat{A}\to\cat{D}$ is a discrete \Conduche fibration, then its pullback
$F^*(\cat{A})\to\cat{C}$ is also discrete \Conduche by the pasting lemma for pullbacks.
Hence $\Delta_F$ is defined as the restriction of $F^*$ on $\DCF/\cat{D}$, and $F^*\circ U_\cat{D}=U_\cat{C}\circ\Delta_F$.

Let $F_!=F\scirc(-)\colon\Cat/\cat{C}\to\Cat/\cat{D}$ be the left adjoint of $F^*$. We define
$\Sigma_F\colon\DCF/\cat{C}\to\DCF/\cat{D}$ to be the composite $L_\cat{D}\circ F_!\circ U_{\cat{C}}$.
Using the fact that $U_\cat{C}$ is fully faithful,
a calculation shows that $\Sigma_F$ is indeed left adjoint to $\Delta_F$:
\begin{equation*}
[\Sigma_F\cat{A},\cat{B}]=[(L_\cat{D}\scirc F_!\scirc U_\cat{C})\cat{A},\cat{B}]
\cong[U_\cat{C}\cat{A},(F^*\scirc U_\cat{D})\cat{B}]\cong[U_\cat{C}\cat{A},(U_\cat{C}\scirc \Delta_F)\cat{B}]
=[\cat{A},\Delta_F\cat{B}].\qedhere
\end{equation*}\end{proof}

\subsection{The equivalence \texorpdfstring{$\Shv{\Int}\cong\DCF/\cat{R}$}{IntDCFR}}
	\label{sec:equiv_shf_conduche}

Having discussed discrete \Conduche fibrations and their properties, we are ready
to show that the topos of interval sheaves is a special case.
The notion of a factorization-linear category $\cat{C}$ (\cref{def:factorization_linear})  turns out to capture all the necessary structure for the slice category $\DCF/\cat{C}$ to be a sheaf topos.

\begin{theorem}\cite[Prop.~3.6]{Johnstone:1999a}\label{thm:DCFR_Int}
Suppose that $\cat{C}$ is a factorization-linear category, let $\tw{\cat{C}}$ be its
twisted arrow category with its Johnstone coverage (\cref{def:Johnstone_coverage}),
and let $\Shv{\tw{\cat{C}}}$ be the associated sheaf topos. Then there is an equivalence of categories
\begin{displaymath}
\DCF/\cat{C}\simeq\Shv{\tw{\cat{C}}}.
\end{displaymath}
\end{theorem}

\begin{corollary}\label{cor:DCF_R=Int}
There is an equivalence between the topos of continuous sheaves (resp.\ discrete-interval sheaves) and discrete \Conduche fibrations over $\cat{R}$ (resp.\ over $\cat{N}$):
\begin{equation}\label{eqn:DCF_Int}
\DCF/\cat{R}\simeq\Shv{\Int}\quad\mathrm{and}\quad
\DCF/\cat{N}\simeq\Shv{\Int}_N
\end{equation}
\end{corollary}

\begin{remark}\label{rem:conduche_vs_sheaves_square}
Note that \cref{prop:mapping_DCF} does not follow from \cref{prop:int_sheaves_presheaves},
even though $\DCF/\cat{N}\simeq\Shv{\Int}_N$ and $\DCF/\cat{R}\simeq\Shv{\Int}$
and the upper adjunctions are essentially the same;
this is because discrete \Conduche fibrations and sheaves give different perspectives.
Categories emphasize composition, and adding
the \Conduche condition enforces that morphisms can be factorized. Conversely, presheaves
emphasize restriction, and adding the sheaf condition enforces that sections can be glued.

These two perspectives compare as follows.
Let $\WGrph$ denote the category of weighted graphs (with nonnegative edge weights),
i.e.\ objects are $G=\{E\rightrightarrows V,E\to\RRnn\}$.
If we define $\cat{I}$ to be the category with objects $\RRnn\sqcup\{v\}$ and two morphisms $s_\ti,t_\ti\colon v\to \ti$ for each $\ti\in\RRnn$, then we have $\WGrph\cong\Psh{\cat{I}}$.

There is a functor $\cat{I}\to\Int$ sending $v\mapsto 0$
and $\ti\mapsto \ti$ for all $\ti\in\RRnn$, and sending $s_\ti\mapsto\Tr_{0}$ and $t_\ti\mapsto \Tr_{\ti}$, the left and right endpoints. This
induces a left Kan extension between the presheaf categories, $\WGrph\to\Psh{\Int}$.
We also have a left adjoint $\WGrph\to\Cat/\cat{R}$, given by the free category construction, whose functor
to $\cat{R}$ sends a path to the sum of its weights. The diagram of left adjoints commutes:
\[
\begin{tikzcd}[sep=5pt]
	\Grph\ar[rr]\ar[dr]\ar[dd]&&\Psh{\Int_N}\ar[dr]\ar[dd]\\
	&\WGrph\ar[rr, crossing over]&&\Psh{\Int}\ar[dd]\\
	\Cat/\cat{N}\ar[rr]\ar[dr]&&\DCF/\cat{N}\simeq\Shv{\Int}_N\ar[dr]\\
	&\Cat/\cat{R}\ar[rr]\ar[from=uu, crossing over]&&\DCF/\cat{R}\simeq\Shv{\Int}
\end{tikzcd}
\]
\end{remark}

\subsection{The \Conduche perspective on interval sheaves and machines}\label{sec:B3}

The equivalence \cref{eqn:DCF_Int} between interval sheaves and discrete \Conduche fibrations
is in particular expressed as follows.
To every $\Int$-sheaf $A$, we may associate a category $\ascat{A}$
called its \emph{associated category}, as well as a functor $\len\colon\ascat{A}\to\cat{R}$
called its \emph{length} functor, similarly to \cref{length_functor}.
Explicitly, the object set of $\ascat{A}$ is the set $\Ob\ascat{A}\coloneqq A(0)$
of germs in $A$; morphisms in $\ascat{A}$ are sections $a\in A(\ti)$
of arbitrary length; composition is given by gluing sections. The functor $\len$ assigns
to each morphism $a$ its length $\len(a)\coloneqq\ti$.
Moreover, sheaf morphisms $F\colon A\to B$ correspond to
length-preserving functors $\ascat{F}\colon\ascat{A}\to\ascat{B}$ over $\cat{R}$.

Under the above correspondence, we can view continuous machines, \cref{def:ADS}, as
\begin{displaymath}
\ADS(A,B)=\Shv{\Int}/(A\times B)\simeq\big(\DCF/\cat{R}\big)\big/{\scriptstyle{\left(\ascat{A\times B}\xrightarrow{\len}\cat{R}\right)}}
\cong\DCF/\ascat{A}\times\ascat{B},
\end{displaymath}
namely themselves as discrete \Conduche fibrations over the product of the associated categories of the input
and output interval sheaves.

Moreover, the notions of totality and determinism for sheaf morphisms defined by \cref{prop:total_det_via_extensions}
also have equivalent expressions in the language of discrete \Conduche fibrations. Let $p\colon S\to A$ be a sheaf morphism,
let $\ascat{p}\colon\ascat{S}\to\ascat{A}$ be the associated functor and $s$ the source map.
Then $p$ is total (resp.\ deterministic) if and only if the induced function $h$
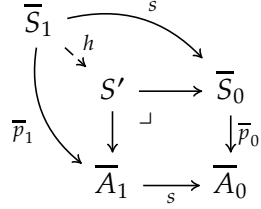
\begin{displaymath}
\begin{tikzcd}[column sep=.15in,row sep=.2in]
\ascat{S}_1 \arrow[drrr, bend left, "s"] \arrow[ddr, bend right, "\ascat{p}_1"'] \arrow[dr, dashed, "h"] &[-5pt] & & \\[-10pt]
& S' \ar[drr,phantom,very near start,"\lrcorner"]\arrow[rr]\arrow[d]
&& \ascat{S}_0\arrow[d, "\ascat{p}_0"] \\
& \ascat{A}_1 \arrow[rr,"s"'] && \ascat{A}_0
\end{tikzcd}
\end{displaymath}
is surjective (resp.\ injective); this is precisely condition (\cref{itm_fn_surj_all}).

Notice that the above conditions of $\ascat{p}$ could in fact be defined for an arbitrary functor $F\colon\cat{C}\to\cat{D}$. In the spirit of \cref{def:dopf_df_dcf},
one could define $F$ to be a (discrete) \emph{epi-opfibration} (resp.\ \emph{mono-opfibration})
if $h\colon\cat{C}_1\to \cat{D}_1\times_{\cat{D}_0}\cat{C}_0$ is surjective (resp.\ injective).
These express whether for each morphism in the base category and object above, say, the target,
there exists at least one, or maximum one, appropriate lifting in the domain category.
In the case $h$ is bijective, we recover the notion of a discrete opfibration.

Thus, $p$ is total (resp. deterministic) in the sense of \cref{def:total_det_map}
if and only if $\ascat{p}$ is an epi-opfibration (resp.\ mono-opfibration).
More informally, if and only if for all functors $a,b$ as shown in the diagram on the left (resp.\ right),
there exists a dotted lift:
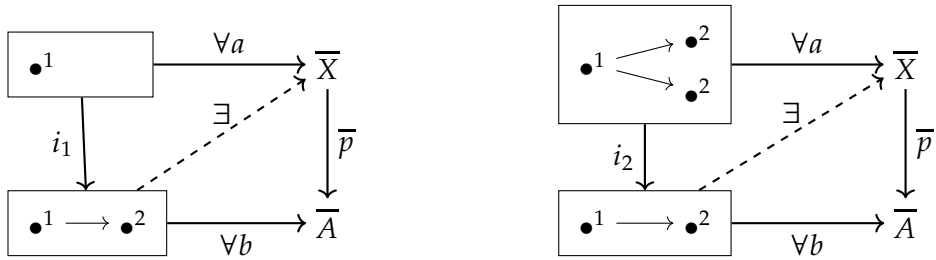
\[
\begin{tikzpicture}[baseline]
	\node (b1) {$\bullet^1$};
	\node[right=.75 of b1] (b2) {};
	\node[below=1.5 of b1] (b3) {$\bullet^1$};
	\node at (b2 |- b3) (b4) {$\bullet^2$};
	\node[rectangle, draw, fit=(b1) (b2)] (c1) {};
	\node[rectangle, draw, fit=(b3) (b4)] (c2) {};	
	\draw[->] (b3) -- (b4);
	\node[right=2 of c1] (c3) {$\ol{X}$};
	\node at (c3 |-c2) (c4) {$\ol{A}$};
	\draw[->, thick] (c1) to node[left] {$i_1$} (c2);
	\draw[->, thick] (c3) to node[right] {$\ol{p}$} (c4);
	\draw[->, thick] (c1) to node[above] {$\forall a$} (c3);
	\draw[->, thick] (c2) to node[below] {$\forall b$} (c4);
	\draw[->, thick, dashed] (c2) to node[above] {$\exists$} (c3);
\end{tikzpicture}
\hspace{1in}
\begin{tikzpicture}[baseline]
	\node (b1) {$\bullet^1$};
	\node[above right=-.25 and .75 of b1] (b21) {$\bullet^2$};
	\node[below right=-.25 and .75 of b1] (b22) {$\bullet^2$};
	\node[below=1.5 of b1] (b3) {$\bullet^1$};
	\node at (b21 |- b3) (b4) {$\bullet^2$};
	\node[rectangle, draw, fit=(b1) (b21) (b22)] (c1) {};
	\node[rectangle, draw, fit=(b3) (b4)] (c2) {};	
	\draw[->] (b1) -- (b21);
	\draw[->] (b1) -- (b22);
	\draw[->] (b3) -- (b4);
	\node[right=2 of c1] (c3) {$\ol{X}$};
	\node at (c3 |-c2) (c4) {$\ol{A}$};
	\draw[->, thick] (c1) to node[left] {$i_2$} (c2);
	\draw[->, thick] (c3) to node[right] {$\ol{p}$} (c4);
	\draw[->, thick] (c1) to node[above] {$\forall a$} (c3);
	\draw[->, thick] (c2) to node[below] {$\forall b$} (c4);
	\draw[->, thick, dashed] (c2) to node[above] {$\exists$} (c3);
\end{tikzpicture}
\]
where $i_1$ and $i_2$ are the obvious functors preserving object labels.
Finally, $p$ is total and deterministic if and only if $\ascat{p}$ is a discrete opfibration.

\bibliographystyle{alpha}
\bibliography{Library}

\end{document}